\documentstyle{article}
\newcommand{\version}{%
(May 18, 2007)}
\input xy
\xyoption{all}
\begin{document}
%
\newcommand{\refequ}[1]{$(\ref{#1})$}
\newtheorem{thm}     {Theorem}[section]
\newtheorem{defin}[thm]   {Definition}
\newtheorem{prop}[thm]    {Proposition}
\newtheorem{propdef}[thm]    {Proposition-Definition}
\newtheorem{cor}[thm]     {Corollary}
\newtheorem{lemma}[thm]   {Lemma}
\newtheorem{conj}[thm]    {Conjecture}
\newtheorem{question}[thm]    {Question}
\newtheorem{notation}[thm]{Notation}
\newcommand{\nothm}{\thethm\ \addtocounter{thm}{1}}
\newcommand{\nothmlbl}[1]{\newcounter{#1}\setcounter{#1}{\thethm}\nothm}
\newcommand{\rmk}{{\bf Remark }}
\newcommand{\rmks}{{\bf Remarks }}
\newcommand{\finrmk}{\medbreak\noindent}
\newcommand{\xcite}[1]{\cite{{#1}}}
\newcommand{\plcomment}[1]{{\bf COMMENTAIRE DE PASCAL:\{\{}{#1}{\bf\}\}}}
\newcommand{\plfootnote}[1]{{\footnote{COMMENTAIRE DE PASCAL: {#1}}}}
%
%
\newcommand{\proof}{{\sl Proof. }\par}
\newcommand{\Proof}[1]{{\bf Proof {#1}. }\par}
%
\newcommand{\cqfd} {\mbox{}\nolinebreak\hfill\rule{2mm}{2mm}\medbreak\par}
%
\newcommand{\biindice}[3]%
{
\renewcommand{\arraystretch}{0.5}
\begin{array}[t]{c}
#1\\ {\scriptstyle #2}\\ {\scriptstyle #3}
\end{array}
\renewcommand{\arraystretch}{1}
}
%
\newcommand{\Bbb}[1]{{\bf{#1}}}
\newcommand{\BR}{\Bbb R}         
\newcommand{\BQ}{\Bbb Q}         
\newcommand{\BC}{\Bbb C}         
\newcommand{\BN}{\Bbb N}         
\newcommand{\BZ}{\Bbb Z}         
\newcommand{\BL}{\Bbb L}         
\newcommand{\Bk}{\Bbb k}         
\newcommand{\calL}{{\cal L}}     
\newcommand{\calU}{{\cal U}}     
\newcommand{\calB}{{\cal B}}
\newcommand{\calC}{{\cal C}}
\newcommand{\calS}{{\cal S}}
\newcommand{\CP}{{\BC P}}
\newcommand{\SSets}{\mathrm{SSets}} 
\newcommand{\Top}{\mathrm{Top}} 
\newcommand{\quism}{\stackrel{\simeq}{\rightarrow}}  
\newcommand{\backquism}{\stackrel{\simeq}{\leftarrow}}
\newcommand{\iso}{\stackrel{\cong}{\rightarrow}}    
\newcommand{\cofarrow}{\to}
\newcommand{\Apl}{A_{PL}}                            
\newcommand{\Aplsimpl}{A_{PL}^{simpl}}
\newcommand{\Tor}{\mathrm{Tor}}
\newcommand{\Ext}{\mathrm{Ext}}
\newcommand{\im}{\mathrm{im}}
\newcommand{\pr}{\mathrm{pr}}
\newcommand{\Hom}{\mathrm{Hom}}
\newcommand{\coker}{\mathrm{coker}}
\newcommand{\gldim}{\mathrm{gldim}}
\newcommand{\fdim}{\mathrm{fdim}}
\newcommand{\cat}{{\mathrm{cat}}}
\newcommand{\nil}{{\mathrm{nil}}}
\newcommand{\Cl}{{\mathrm{Cl}}}
\newcommand{\id}{{\mathrm{id}}}
\newcommand{\Sing}{{\mathrm{Sing}}}
\newcommand{\rank}{\mathrm{rank}}
\newcommand{\mod}{\mathrm{mod}}
\newcommand{\del}{\partial}
\newcommand{\circlesimeq}{\bigcirc\,\,\,\simeq}
\newcommand{\decodiagpo}{(\bigcirc\!\!\!\!*\,)}
\newcommand{\bd}{\partial}
\newcommand{\isom}{\stackrel{\cong}{\rightarrow}}    
%
\begin{center}%
{\Large{\bf The rational homotopy type of a blow-up in the
stable case. }
\\ {\normalsize\version}}
\medbreak

{\large\bf Pascal Lambrechts\footnote{P. L. is a Chercheur Qualifi\'e at FNRS.}
and Don Stanley\footnote{Support of the Institut Mittag-Leffler is gratefully acknowledged}}%
\\Universit\'e de Louvain%
\\University of Regina
\end{center}
\vskip5mm
\begin{center}{\bf Abstract}\end{center}
Suppose $f:V\rightarrow W$ is an embedding of closed
oriented manifolds whose normal bundle has the structure of a complex vector bundle.
It is well known in both complex and symplectic geometry
that one can then construct a
manifold $\widetilde W$ which is the blow-up of $W$ along $V$.
Assume that $\dim W\geq 2\dim V +3$ and that $H^1(f)$ is injective. We construct an algebraic
model of the rational homotopy type of the blow-up $\widetilde{W}$
from an algebraic model of the embedding and the Chern
classes of the normal bundle. This implies that if the space $W$ is simply connected then the rational
homotopy type of $\widetilde W$ depends only on the rational homotopy class of $f$ and on the Chern
classes of the normal bundle.

{AMS-classification (2000)} : 55P62 -14F35 - 53C15 - 53D05
\vskip1cm
\newpage
%
\section{Introduction}
The blow-up construction comes from complex algebraic geometry;
Gromov \xcite{Gromov} and McDuff \xcite{McDuff} constructed the
blow-up for symplectic manifolds. McDuff used it to construct the
first examples of simply connected non-K\"ahler symplectic
manifolds. In this paper we study the rational homotopy type of
the blow-up construction.
\par
It is well known that all K\"ahler manifolds are symplectic and
fundamental problem is to find closed symplectic manifolds which cannot
be given a K\"ahler structure. A non-simply connected example was first found by
Thurston \cite{Thurston}, but to find simply connected examples proved
more difficult. In fact McDuff introduced the blow-up construction to
resolve this problem. She showed that the blow-up of Thurston's examaple
in a complex projective space has an odd third Betti number and is thus
not K\"ahler. This leads to more precise structural questions.
\par
A space is formal if its rational cohomology algebra serves as a rational model of it.
In particular this implies its cohomology algebra determines its rational homotopy
type and that it has no nontrivial Massey products in its cohomology. It has been shown
\cite{DGMS} that all closed K\"ahler manfolds are formal. Thurston's example was known
to be non formal. By constructing non-trivial Massey products Babenko and Taimanov
\cite{BabenkoTaimanov1} showed that the McDuff example is not
formal. Rudyak and Tralle \cite{RudyakTralle}
extended these results for the blow-ups along many other manifolds.
One of the applications of our model is to prove \cite{LS1} that the
blow-up along a manifold $M$ symplectically embedded in a large
enough $\CP(n)$ is formal if and only if $M$ is formal.
\par
Although McDuff was interested in the blow-up for symplectic manifolds,
the construction is more general.
Suppose $f:V\rightarrow W$ is an embedding of smooth
closed oriented manifolds
and that the normal bundle of the embedding has been given the structure of a
complex vector bundle.
This is enough data to construct the blow-up
$\widetilde W$ of $W$ along $V$ as described by McDuff
\cite{McDuff} and outlined in Section \ref{sec-outline} below.
If $f$ is an embedding of complex manifolds
with the canonical complex structure on the normal bundle then
this blow-up
is homeomorphic to the classical
blow-up of $W$ along $V$.
\par
Under some restrictions
we give a complete description of the rational homotopy type of
the blow-up $\widetilde W$ using only the rational homotopy class of
the map $f:V\rightarrow W$ and the Chern classes of the normal
bundle of $V$.
The rational homotopy type of a space is defined as long as all
the spaces involved are nilpotent (see Section \ref{func}).
Note that simply connected closed manifolds are always
nilpotent.
\begin{thm}
\label{blorgial}
(Theorem \ref{blorgy})
Let $f:V \rightarrow W$ be an embedding of closed orientable smooth manifolds and
suppose that the normal bundle $\nu$ is equipped with the structure of a complex
vector bundle. Assume that $\dim W\geq 2\dim V +3$,
that $W$ is simply connected and that $V$ is nilpotent.
Then the rational homotopy type of the blow-up
of $W$ along $V$, $\widetilde W$ can be explicitly determined
from the rational homotopy type
of $f$ and from the Chern classes $c_i(\nu)\in H^{2i}(V;\BQ)$.
\end{thm}
More generally our description holds as long as $V$, $W$
and $\widetilde W$ are
nilpotent, $H^1(f;\BQ)$ is injective and $\dim (W)\geq 2\dim(V) +3$
(Corollary \ref{UBC}). It is not so easy to determine if $\widetilde W$ is nilpotent, so the usual
case is $W$ simply connected.
We also note that without the dimension restriction, $\dim (W)\geq 2\dim(V) +3$,
there exists manifolds $V$ and $W$, and
homotopic embeddings $V \rightarrow W$ with isomorphic complex
normal bundles whose blow-ups are not rationally equivalent
\cite{Lambrechts}. Thus our dimension restriction, $\dim(W)\geq
2\dim(V)+3$, cannot be discarded.
\par
Sullivan  \cite{Sullivan} studies rational homotopy using
his piecewise linear forms $\Apl(\_)$ functor,
which is analagous to
the de Rham differential forms functor, $\Omega^*(\_)$.
In particular it is contravariant and
for a topological space $X$, $\Apl(X)$ is a commutative differential
graded algebra or CDGA (defined in Section \ref{sec-homo}).
Another way to think of $\Apl(X)$ is as a commutative version of
the cochains on $X$ with coefficients in $\BQ$, $C^*(X;\BQ)$.
A {\em model} of $X$ is any CDGA
weakly equivalent
to $\Apl(X)$ (see Section \ref{sec-cat}).
Similarly a model of a map $f:X \rightarrow Y$ is any map
of CDGAs weakly equivalent to $\Apl(f):\Apl(Y)\rightarrow \Apl(X)$.
Under some finiteness conditions any model of $X$
completely determines the rational homotopy type of $X$ (see Section \ref{func}). It is in
this sense that we determine the rational homotopy type of
$\widetilde W$. In fact we determine the homotopy type of
the CDGA $\Apl(\widetilde W)$ without any finiteness restrictions.
\par
In Section \ref{desc of model} and Theorem \ref{lemma-modelblowup}
using only a model for $f$ and the Chern
classes of the normal bundle we construct an explicit model for $\widetilde W$.
This is our main result and
Theorem \ref{blorgial} follows directly from it.
There are a number of interesting byproducts produced along the way to the proof.
For example we derive algebraic models for the compliment of the
embedding (see also \cite{LS}) and for the projectivization of a complex bundle.
\par
In the last section of the paper we give a few applications of the
model of the blow-up. First we study the special case of a blow-up
of $\CP(n)$, and looking at McDuff's example of the blow-up of
$\CP(n)$ along the Kodaira-Thurston manifold we
prove the existence of non-trivial Massey products by direct
calculation.  Our next
application is to calculate the cohomology algebra of the blow-up
along $f:V \rightarrow W$ under our dimension restrictions
(Section \ref{coho alg}). This is complementary to work of Gitler
\cite{Gitler} who gave a different description of this algebra when $H^*(f)$ is surjective.
In Section \ref{vlast} we use this calculation to show that there are infinitely many distinct rational
homotopy types of symplectic manifolds that can be constructed as the
blow-up of $\CP(5)$ along $\CP(1)$.

\subsection{Contents}
\noindent
1) Introduction
\par\noindent
2) Modelling the blow-up
\par\noindent
3) Background and notation
\par\noindent
4) Thom class and the shriek map
\par\noindent
5) Model of the complement of a submanifold
\par\noindent
6) Model of the projectivization of a complex bundle
\par\noindent
7) The model of the blow-up
\par\noindent
8) Applications
\vskip 0.4cm
A more detailed list of contents appears at the beginning of each section.
\subsection{Acknowledgements}
The authors would like to thank Professor Hans Baues for help with
homological algebra and for pointing out that there is a subtle
difficulty in finding the model of a pushout from the pullback of
the models. We would also like to thank Yves F\'elix, Francois
Lescure and Yuli Rudyak for enlightening conversations regarding
this work. We are very grateful to Barry Jessup for reading a
draft of the paper and suggesting many improvements. The second
author would like to thank Terry Gannon for his encouragement and
support. A large portion of this work was done while the authors
were at the Max Plank Institute, Universit\'e de Lille-1,
Universit\'e de Louvain and the Mittag-Leffler Institut. We would like to thank them for providing
great environments for research. Both authors also benefited from
the funding of FNRS, CNRS and NSERC without whose generosity this
work could never have been completed.
\section{Modelling the blow-up}\label{sec-outline}

In this section we first describe the topology of the blow-up construction and then
describe the model of the blow-up. The precise statement may be found
in Theorem \ref{lemma-modelblowup}.


Again suppose $f \colon V\rightarrow W$ is an
embedding of connected closed oriented manifolds and suppose
the normal
bundle $\nu$ of $f$ has been given a complex structure. Let $T$ be
a tubular neighborhood of $V$ in $W$. Let $\bd T$ be the boundary
of $T$ and $B=\overline{W\setminus T}$. Then $T\cup B=W$ and
$T\cap B=\bd T$. Hence we have a pushout
\begin{equation}\label{love}
\xymatrix
{
\bd T \ar[r]^k \ar[d] & B \ar[d] \\
T \ar[r] & W.
}
\end{equation}
By the Tubular Neighborhood Theorem \cite[Theorem 11.1]{MS} there is a
diffeomorphism between $T$ and the disc
bundle $D\nu$ that sends $V$ to the zero section of $D\nu$ and sends $\bd T$ to the sphere bundle
$S\nu$. Since $\nu$ is a complex bundle we can quotient by the
$S^1\subset \BC^*$-action on $S\nu\cong \bd T$. We obtain a complex
projective bundle $P\nu$ over $V$ and a commutative diagram
$$
\xymatrix
{
\bd T \ar[r]^q \ar[d] & P\nu \ar[dl] \\
V
}
$$
Next we can remove $T$ from $W$ and instead of putting it back
as in (\ref{love}) we can replace it by $P\nu$. This gives us a pushout
\begin{equation}
\label{lover}
\xymatrix
{
\bd T \ar[r]^k \ar[d]_{q} & B \ar[d] \\
P\nu \ar[r] & \widetilde W.
}
\end{equation}
The space $\widetilde W$ is called the {\em blow-up of $W$ along $V$}. This actually
only gives us the homeomorphism type of the blow-up but with
slightly more care we
can get the diffeomorphism type.
Since we are only studying the rational homotopy type, homeomorphism
type is more than enough. An important point is that (\ref{lover})
is also a homotopy pushout.
\par
In the next few paragraphs we analyze Diagram (\ref{lover}) and
give an idea of how we construct the model ${\cal B}(R,Q)$ of
$\Apl(\widetilde W)$. We begin with the fact that $\Apl(\_)$
takes homotopy pushouts to homotopy
pullbacks. So we are left with analyzing the homotopy pullback of
the diagram
$$
\xymatrix
{
\Apl(\bd T)
& \Apl(B) \ar[l]_{\Apl(k)} \\
\Apl(P\nu) \ar[u]^{\Apl(q)}
}
$$
which we get from applying $\Apl(\_)$ to Diagram (\ref{lover}).
To do this we construct models of the two ``legs" of the diagram,
the maps $\Apl(B)\rightarrow \Apl(\bd T)$ and $\Apl(P\nu)
\rightarrow \Apl(\bd T)$. The pullback of these models should be a
model of the pullback of the diagram. However we
need to be careful how the models are ``glued" together because different
gluings correspond to homotopy automorphisms of $\Apl(\bd T)$ and
can lead to pullbacks that are not homotopy equivalent.
\par
To construct a model of $k$ we make use of our cochain level version
of the classical shriek
map $f^!\colon H^{*-r}(V) \rightarrow H^*(W)$, where $r$ is the codimension
of $V$ in $W$.
Under our hypothesis that
$\dim(W)\geq 2\dim(V)+3$,
a model of $k$ can be constructed from any
model of $f$ and so the rational homotopy type of $\bd T$ and $B$
and the rational homotopy class of $k$ depend only on the rational homotopy class of
$f$. This is to be expected since we are in the stable range where
homotopic maps are isotopic. It is shown in \cite{LS} that
without our dimension restriction the rational homotopy class of $k$ can depend
on more than just the rational homotopy type of $f$, and so it is at this point
that the dimension restriction for the model of the blow-up arises.
In \cite{Lambrechts} it is shown that the blow-up
along homotopic embeddings with the same Chern classes
on their normal bundles can have different homotopy types.
\par
Next we more fully describe the model of $B$.
Suppose
$$
\phi\colon R\rightarrow Q
$$
is a model of $\Apl(f)\colon \Apl(W)\rightarrow \Apl(T)$.
We can consider $Q$ as a differential graded $R$-module (or $R$-dgmodule,
see Section \ref{sec-homo}). Suppose
$u_V\in H^*(Q)$ and $u_W\in H^*(R)$ are orientation classes, where
$m=\dim(V)$ and $n=\dim(W)$.
Let
$$
\phi^!:s^{-r} Q \rightarrow R
$$
be a map of $R$-dgmodules such that
$H(\phi^!_*)(u_V)=u_W$, where $r=n-m$. (Recall that the suspension
operator gives $(s^{-r}Q)^{j+r}=Q^j$, so $u_V\in H^n(s^{-r}Q)$.)
Such a map will be called a {\em shriek map}.
A shriek map always exists (Proposition \ref{prop-shriekexists})
and is unique up to homotopy (Proposition \ref{unish}).
Let $m=\dim V$, $n=\dim W$ and $r=n-m$.
We describe a model of $k: \bd T \rightarrow  B$, and refer the reader to
Lemma \ref{igog} for the details.
\begin{lemma}

Assume $n\geq 2m+3$.
If $R^{\geq n+1}=0$ and $Q^{\geq m+2}=0$ then there exists
explicit CDGA structures on $R\oplus ss^{-r} Q$ and
$Q\oplus ss^{-r} Q$ determined by the CDGA structures on $R$ and
$Q$ and by the shriek map $\phi^!$ such that
the map
\begin{equation}\label{near2}
\phi\oplus id\colon R\oplus ss^{-r} Q\rightarrow Q\oplus  ss^{-r} Q
\end{equation}
is a CDGA model of $\Apl(k)\colon\Apl(B)\rightarrow \Apl(\bd T)$.
\end{lemma}
From any model $\phi\colon R\rightarrow Q$, one satisfying
the degree restrictions can always be constructed
(Proposition \ref{prop-shriekexists}). The differential
on $R\oplus ss^{-r} Q$ comes from the fact that it is actually
the mapping cone on $\phi^!$ (see Section \ref{sec-mapcon}) and the CDGA structure is
what we call
the {\em semi-trivial} CDGA structure (see Definition \ref{def-semitrivCDGA}).
\par
Constructing a model of $q:\bd T\rightarrow P\nu$
is more straightforward.
The cohomology algebra of the projective bundle can be used
to define the Chern classes $c_i(\nu)$ of the
normal bundle $\nu$ \cite[IV.20]{BottTu}. This description allows us to
construct a model of $P\nu$. The free graded commutative algebra on the graded
generators $a_i$ is
denoted by $\Lambda(a_1, \dots , a_n)$.
Next we describe a model of $q:\bd T\rightarrow P\nu$,
and refer to Theorem \ref{prop-modelprojxi}
for the details.
\begin{thm}
Assume $n\geq 2m+3$.Let $2k=n-m$
and let $\gamma_0=1$ and
$\gamma_i\in Q$ be representatives of the Chern classes
$c_i(\nu)$. Suppose $|x|=2$ and $|z|=2k-1$.
Define CDGAs $(Q \otimes \Lambda(x,z);D)$ by $Dx=0$ and
$Dz= \Sigma_{i=0}^{k-1}\gamma_i x^{k-i}$ with  and $(Q\otimes \Lambda z;\overline D)$ by
$\overline D(z)=0$.
Then the projection map
\begin{equation}\label{loop}
{\rm proj}\colon Q \otimes \Lambda(x,z) \rightarrow Q\otimes \Lambda z
\end{equation}
sending $x$ to $0$
is a CDGA model of $\Apl(q):\Apl(P\nu)\rightarrow \Apl(\bd T)$.
\end{thm}
\par
In the last theorem the dimension restriction $n\geq 2m+3$ arises since
we require that $\overline{D}z=0$, which represents the Euler class of the
bundle, but other than that are not really needed and probably
a similar theorem without the dimension restrictions is true.
Having constructed models for both
$k$ and $q$ of Diagram (\ref{lover}) we can use
them to construct a model ${\cal B}(R,Q)$ of the pushout $\widetilde W$
by taking the pullback of these two
models. It is described precisely in Section \ref{desc of model}.
The form of ${\cal B}(R,Q)$ is:
$$
{\cal B}(R,Q)=\left(R\oplus Q \otimes \Lambda^+(x,z),D\right).
$$
As mentioned before, to get ${\cal B}(R,Q)$ we have to be careful how the
models of $k$ and $q$ fit
together. In the model (\ref{near2}) of $k$ the model of $\Apl(\bd T)$ is
$Q \oplus ss^{-r} Q$ whereas in the model (\ref{loop}) of $q$ it is
$Q\otimes \Lambda(z)$. These are always isomorphic as $Q$-modules and our dimension
restrictions imply that they are isomorphic as CDGAs. The problem is to
make sure we pick the correct isomorphism. While we are constructing our models
of $k$ and $q$, we keep track of the isomorphism with the help of orientation
classes on both of our manifolds and on the normal bundle.
In our special situation, this orientation information together with
the $Q$-dgmodule structures is enough to determine the isomorphism.
Once we have this nailed down it is straightforward to construct
a model of the blow-up from our models of $k$ and $q$ using a pullback.

\section{Background and notation}\label{b-and-n}
For this paper the ground field and the coefficient ring
will be the rational numbers,
unless otherwise stated. For a topological space $X$,
$H^*(X)$ refers to its singular cohomology.
We denote also by $H(\_)$ the functor taking homology of a
differential complex. We will use $H^n(\_)$ to refer both
to the $n$-th singular cohomology group and to the homology in
degree
$n$ of any differential graded object.
\subsection{Categorical preliminaries}\label{sec-cat}
For a small category $I$ and any category $\cal D$, let
${\cal D}^I$ be the {\em diagram category} defined as follows:
the objects of ${\cal D}^I$
are the functors $I \rightarrow {\cal D}$ and the morphisms are
natural transformations. We often refer to the objects in such a
category as diagrams.
For example our category $I$
can consist of two objects with exactly one non-identity map joining them
(this category can be depicted as $\bullet \rightarrow \bullet$)
and each diagram in ${\cal D}^I$ corresponds to a single map in $\cal D$.
Similarly $I$ could be the category
$\xymatrix{ \bullet & \bullet \ar[l] \ar[r] & \bullet}$
which corresponds to the data of a pushout in $\cal D$ or its dual
$\xymatrix{ \bullet\ar[r] & \bullet & \bullet\ar[l]}$
(data for a pullback). We also use squares of objects which correspond
to the category
$$
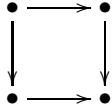
\xymatrix
{
\bullet \ar[r] \ar[d] & \bullet \ar[d] \\
\bullet \ar[r] & \bullet
}
$$
If
$$
{\mathbf A}=
\xymatrix
{
A_1 \ar[r]\ar[d] & A_2 \ar[d]\\
A_3 \ar[r] & A_4
}
$$
and
$$
{\mathbf B}=
\xymatrix
{
B_1 \ar[r]\ar[d] & B_2 \ar[d]\\
B_3 \ar[r] & B_4
}
$$
are squares of objects and
$$\alpha_i\colon A_i\rightarrow B_i$$ are maps, then
$$
{\left(
\begin{tabular}{cc}$\alpha_1$&$\alpha_2$\\$\alpha_3$&$\alpha_4$\end{tabular}
\right)}\colon
{\mathbf A}\rightarrow {\mathbf B}
$$
denotes the map between the squares that
would traditionally be denoted as a commutative cube:
$$
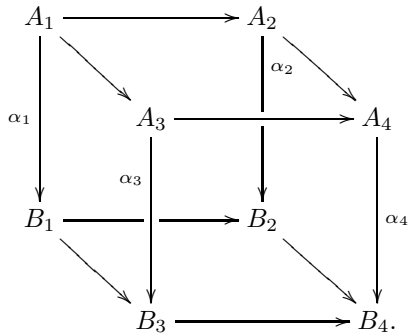
\xymatrix
{
A_1 \ar[dd]_{\alpha_1} \ar[rr]
\ar[dr] && A_2
\ar'[d]^{\alpha_2}  [dd]\ar[dr] \\
& A_3 \ar[dd]_(0.3){\alpha_3} \ar[rr] && A_4 \ar[dd]^{\alpha_4} \\
B_1 \ar'[r] [rr] \ar[dr]
&& B_2 \ar[dr] \\
& B_3 \ar[rr] && B_4.}
$$

\par

If $\cal D$ is the category of topological
spaces and $F, G\in {\cal D}^I$ then a morphism
$\eta\colon F\rightarrow G$ is a {\em homeomorphism} if for each
$i\in I$, $\eta(i)\colon F(i) \rightarrow G(i)$
is a homeomorphism for each element.
If $\cal D$ has weak equivalences then $\eta$
is called a {\em weak equivalence}
if $\eta(i)$ is a weak equivalence for each $i\in I$.
Two diagrams are weakly equivalent if they are connected by a
chain of weak equivalences. So $F,G\in {\cal D}^I$ would be weakly equivalent if there exists
a diagram in ${\cal D}^I$ as below with all the maps weak equivalences.
$$
\xymatrix
{
F & F_1 \ar[r]^{\simeq}\ar[l]_{\simeq} & F_2 & \cdots \ar[l]_{\simeq} &
F_n \ar[r]^{\simeq}\ar[l]_{\simeq} & G
}
$$
Notice that there may not be a direct map between $F$ and $G$.

\subsection{Homotopy theory of CDGA and $R$-DGMod}
\label{sec-homo}

A good reference for the categories
CDGA and $R$-DGmod is the book of Felix-Halperin-Thomas
\cite{FHT-book}.
Next we review the notion of CDGA. Let
$A=\oplus_{i=0}^{\infty} A^i$ be a graded vector space together
with an associative multiplication $\mu: A \otimes A \rightarrow A$,
a unit $1\in A^0$ and a
linear map of degree $+1$, $d:A \rightarrow A$, called the differential. We denote
$\mu(a,b)$ by $ab$ or sometimes $a\cdot b$.
 If $a\in A^n$ we write $|a|=n$ and say that
 $a$ has {\rm degree} $n$. We require
that the multiplication $\mu$ is graded commutative,
that $d^2=0$ and that $d$ is a derivation . In other
words  we require that $|ab|=|a|+|b|$, that $ab=(-1)^{|a||b|}ba$ and that the
Leibnitz law $d(ab)=(da)b + (-1)^{|a|}a(db)$ holds. Such an $(A,d)$ is
called a commutative differential graded algebra or CDGA.
Notice that the multiplication is suppressed. Sometimes we just
write $A$ also suppressing the differential.
We call $A$ {\em connected} if $A^0=\BQ \cdot 1$. Maps
between CDGAs are graded vector space homomorphisms that
commute with the multiplication and the differential. Clearly
we get a category which we will also denote CDGA.
A special case of CDGA are Sullivan algebras $(\Lambda V,d)$ where
$\Lambda V$ is the free graded commutative algebra on a graded vector space
$V$ and $d$ is a differential on $\Lambda V$ satisfying a certain nilpotence condition
as described in \cite[Part II]{FHT-book}. However we will relax
their condition that $V$ should be concentrated
in positive degrees
and just demand that $V$ be concentrated in non-negative degrees.
If $V=\BQ\langle a_i,b_j\rangle$ is a graded vector space with basis
$\{ a_i \}\cup \{ b_j\}$ where the $a_i$ are homogeneous in even degrees and the $b_j$
are homogeneous in odd degree then as a graded algebra
$\Lambda V\cong P(a_i)\otimes E(b_i)$ where $P(a_i)$ is the free
graded polynomial algebra on the $a_i$ and $E(b_i)$ is the free
graded exterior algebra on the $b_i$.

\par

We can also define a category of differential graded modules over
some fixed CDGA. Let $(R,d)$ be a CDGA. Let
$M=\oplus_{i=-\infty}^{\infty}M^i$ be a differential graded
$R$-module with structure map $\mu:R\otimes M \rightarrow M$ and
differential of degree $+1$, $d': M \rightarrow M$. We require
that $d'(rm)=(dr)m+(-1)^{|r|}d'm$. Of course $M$ is an $R$-module
in the ungraded sense as well and $|rm|=|r|+|m|$. We call such an
object a differential graded $R$-module or an $R$-dgmodule. Maps
are $R$-module maps that preserve the grading and
commute with the differential.
The category of $R$-dgmodules and $R$-dgmodule maps is denoted by
$R$-DGmod.
\par
Weakly equivalent diagrams of CDGA or $R$-dgmodules are
 referred to as {\em models} of each other.
For example let $f:X \rightarrow Y$ be a map of spaces. If a diagram
$\phi:A\rightarrow B$ of CDGA is weakly equivalent to the diagram
$\Apl(f):\Apl(Y) \rightarrow \Apl(X)$, then we will
say that $\phi:A \rightarrow B$ is a model of $\Apl(f)$.
In this case, as is common,
we also call $\phi$ a model of $f$. Similarly if $A$ is a model of $\Apl(X)$
we call $A$ a model of $X$.
The concept of a CDGA model of a space or a map is well
established.
Note however that many authors reserve the term model to refer to free
CGDA whereas we use it more generally.
For our purposes we also need models
of more general diagrams. In fact since we will be gluing diagrams
together we will sometimes need that the equivalences between
our diagrams preserve certain extra structure.
\par
Both CDGA and $R$-DGMod are closed model
categories (see \cite{DwyerSpalinski} for a
review of model categories and \cite{BG} for a proof of the fact that
CDGA satisfies the axioms of a closed model category).\
It is not necessary for the reader to be familiar with closed
model categories since all of the relevant results can be proved
directly in these categories.
In both categories the
closed model structure is determined by the following families of maps:
\begin{itemize}
\item fibrations are surjections,
\item weak equivalences are quasi-isomorphisms.
\end{itemize}
We use the terms fibration and surjection interchangeably when working
in CDGA or $R$-DGMod.
By a {\em cellular cofibration} we mean:
\begin{itemize}
\item in CDGA, a relative Sullivan algebra
$B\to B\otimes \Lambda V$
as defined in \cite[Chapter 14]{FHT-book} except that we allow
non-negatively graded $V$ instead of just positively graded $V$.
\item in $R$-DGMod, a semi-free extension
$M\to M\oplus R\otimes V$
as defined in \cite[\S 2]{FHT-Jameshandbook}
\end{itemize}
We will generally deal with cofibrations that are
cellular ones. Note that not all cofibrations are of this form
(see \cite[\S 4.4]{BG}) but at least all cofibrations between connected CDGA are.
A map that is both a cofibration and a weak equivalence
is called an acyclic cofibration. Similarly a map that is both a
fibration and a weak equivalence is called an acyclic fibration.
If $\emptyset$ denotes the initial object
of the category (it is $\BQ$ concentrated
in dimension $0$ in CDGA and $0$ in $R$-DGMod)
then a {\em cellular cofibrant object} is an object $X$ such that
the map $\emptyset\to X$ is a cellular cofibration.
In CDGA cellular cofibrant objects are the Sullivan algebras \cite[\S 12]{FHT-book}
and in $R$-DGMod they are the semi-free $R$-dgmodules \cite[\S 2]{FHT-Jameshandbook}.
Dually if $*$ denotes the terminal object then an object
$X$ is called {\em fibrant}
if $X\rightarrow *$ is a fibration.
Note that all objects in CDGA and $R$-DGMod are fibrant.
\par
The following lemma is a slight modification of one of the axioms of a closed
model category.
\begin{lemma}
Suppose $A$ and $B$ are CDGAs such that $H^0(A)=H^0(B)=0$.
Any CDGA map $f:A\rightarrow B$ can be factored into a cellular cofibration follow
by an acyclic fibration.
\end{lemma}
\proof
This is proved as part of \cite[Theorem 6.1]{Halperin} with the added condition that
$A$ is augmented. However the augmentation is not used to get
the result of the lemma.
\cqfd
The analogous result in $R$-DGmod is also true and follows by standard arguments.
\begin{lemma}
Any map in $R$-DGmod can be factored into an acyclic cellular cofibration followed
by a fibration and also into a cellular cofibration followed by an acyclic fibration.
\end{lemma}
\par
Homotopies in a closed model category can be defined with the help of a
cylinder object which we describe now.
Let $\_ \coprod \_$ denote the coproduct, and $\nabla:
A\coprod A \rightarrow A$ the fold map. Factor $\nabla$ into a
cellular cofibration $i_0 +i_1: A \coprod A \rightarrow Cyl A$ followed by an
acyclic fibration. This implies in particular that
the maps $i_0,i_1\colon A \rightarrow Cyl A$ are weak equivalences.
The object $Cyl A$ is called the {\em cylinder object} of $A$,
it is unique up to homotopy.
   Two morphisms $f_0,f_1:A \rightarrow X$
are {\em homotopic} if there exists a map $H:Cyl A \rightarrow X$ such
that $Hi_0=f_0$ and $Hi_1=f_1$.
 We write $f_0\simeq f_1$, or
$f_0\simeq_R f_1$ if we wish to emphasize the fact that the
homotopy is in the category $R$-DGMod.
Then it is easy to check from the
definition of a homotopy in $R$-DGMod:
\begin{lemma}\label{hom-char}

Two morphisms $f_0,f_1\colon A\to X$ in $R$-DGMod with $A$ cofibrant
are homotopic if there exists an $R$-module degree $-1$ morphism
$$h\colon A\to X$$
such that $d_Xh+hd_A=f_0-f_1$.
Note that such a homotopy can also be seen as an $R$-module
degree $0$ morphism $h\colon sA\to X$ where $s$ is the suspension
(see Definition \ref{def-suspension}).
\end{lemma}
In CDGA our notion of homotopy is also equivalent to the
more traditional one (\cite[Chapter 12(b)]{FHT-book}).
\par
We recall the notion of sets of homotopy classes of maps in CDGA
and $R$-DGmod.
Let $X,Y\in {\rm CDGA\ or \ } R{\rm -DGMod}$. Factor $\emptyset\rightarrow X$ as a cofibration
followed by an acyclic fibration
$\emptyset\rightarrow \hat X \rightarrow X$.
We define the set of homotopy classes of maps from $X$
to $Y$, $[X,Y]$ as the set of equivalence classes of $\Hom(\hat X,Y)$ under the homotopy
relation in either CDGA or $R$-DGMod.
$$
[X,Y]=\Hom(\hat X,Y)/\simeq.
$$
We will write
$[X,Y]_R$ if we wish to emphasize we are looking at homotopy classes
of $R$-dgmodule maps. The Lifting Lemma
(Lemma \ref{lemma-tk-lifting} below) implies that $[X,Y]$ is independent of the choice of
$\hat X$ (any two choices give naturally isomorphic sets) and there is
an obvious map $\Hom(X,Y)\rightarrow [X,Y]$.
Weak equivalences in the range and domain of $[\_,\_]$ induce
bijections of sets. To those already familiar with closed model
categories we point out that we don't have to replace the
range by a fibrant object since all
objects in our category are already fibrant.
\subsection{The functors $\Apl(\_)$ and $|\_|$}
\label{func}
A connected space $X$ is {\em nilpotent} if $\pi_1(X)$ is a nilpotent group and $\pi_n(X)$ is a nilpotent
$\pi_1(X)$-module for each $n\geq 2$ (See \cite[II.2]{HMR}).
Note that all simply connected spaces are nilpotent.
A nilpotent space $X$ is {\em rational}  if $\pi_n(X)$ is a rational vector space for
each $n\geq 1$. A nilpotent space $X$ has
{\em finite $\BQ$-type} if $H_n(X;\BQ)$ is finite dimensional
for each $n\geq 1$.
\par
Let $Top$ be the category of topological spaces.
There are adjoint functors
$$
\Apl:Top \rightarrow CDGA
$$
$$
|\_ |:CDGA \rightarrow Top
$$
introduced by Sullivan
\cite{Sullivan}. We refer the reader
to \cite[II 10 and II 17]{FHT-book}  for their definitions (see also
\cite[Section 8]{BG} and \cite[Section 7 and 8]{Sullivan}).
These functors induce an equivalence
between the homotopy types of finite type CDGA and of
nilpotent rational topological spaces of finite $\BQ$-type
(see \cite[9.4]{BG} for a more precise statement).
It is in this sense that the homotopy type of $\Apl(X)$ as a CDGA
determines the rational homotopy type of $X$. A key property of $\Apl$
that we will use repeatedly is the $\BQ$-version of the natural algebra de Rham isomorphism
\begin{equation}\label{equ-isocan}
H\Apl (\_)\cong H^*(\_;\BQ)
\end{equation}
This isomorphism holds for all spaces and will be fixed throughout the paper.
\par
In order to reduce clutter in our equations, for a map $f$ of spaces we will often
write $f^*$ for $\Apl(f)$.
\subsection{Standard abuses of notation}\label{sec-abu}
Because of the isomorphism (\ref{equ-isocan}),
to any cohomology class
$\mu\in H^n(X;\BZ)$, using the change of coefficients
morphism $H^*(\_, \BZ) \rightarrow H^*(\_, \BQ)$,
we can associate a cohomology class in
$H^n(\Apl(X))$. Abusing notation we also denote this class by
$\mu\in H^n(\Apl(X))$.
\par
Also if $R$ is a CDGA model of another CDGA $R'$ through a fixed chain of weak equivalences

$$
\xymatrix
{
R \ar[r]^{\simeq}
& \dots & R' \ar[l]_{\simeq}
}
$$
and if $\mu\in H^n(R)$ is a cohomology class we will denote the corresponding
class $\mu\in H^n(R')$ by the same name.
Note that which classes correspond only depends on which element of
$[R,R']$ is represented by the chain of weak equivalences.
\par
We will often abuse ``inc" to denote any inclusion  map and ``proj" to denote any
projection map. The abuse is twofold since not only will we often not explicitly define
our inclusions but also we will have many different maps referred to by the same notation.
We will also sometimes use
``$\xymatrix{\ar@{^{(}->}[r] &}$" to denote cofibrations or inclusions.
\subsection{Closed model category facts}
\label{sec-CMCfactoids}
The following result is standard. Let $\mathcal C$ be a closed model category.
\begin{lemma}[Lifting Lemma]
\label{lemma-tk-lifting}
Suppose given the following solid arrow diagram in $\calC$
$$
\xymatrix{
A\ar[r]^f\ar[d]_g&B\ar[d]^{v}\\
C\ar[r]_{u}\ar@{-->}[ru]^l&D.
}
$$
(i) If the diagram is commutative, $g$ is a cofibration,
 $v$ is a fibration and
 $g$ or $v$ is a weak equivalence then
there exists a  lift $l$ making each triangle commutative.
Moreover such a lift $l$ is unique up to homotopy.
\par
\noindent
(ii) If the diagram is commutative
(respectively, commutative up to homotopy), $g$
is a cofibration and $v$ is a weak equivalence
then there exists a  lift $l$ such that $l\circ g=f$
(respectively, $l\circ g\simeq f$) and  $v\circ l\simeq u$;
that is, the upper triangle is commutative and the lower triangle
is commutative up to homotopy (respectively, both triangles are
commutative up to homotopy.)
Moreover such a lift $l$ is unique up to homotopy.
\end{lemma}
\proof
See \cite[Chapter 14]{FHT-book} for the CDGA case, and
\cite[Lemma A.3]{FHT-gorenstein} for the special case of $A=0$ in
$R$-DGmod. The case for general $A$ follows by standard techniques.
\cqfd
The following next two lemmas allow us to convert certain
homotopy commutative diagrams into strictly commutative diagrams.

\begin{lemma}
\label{lemma-tk-nose}
Let $\mathcal C$ be a closed model category and suppose
$$
\xymatrix{
&\hat A\ar[ld]_f\ar[d]^{\overline f}\ar[rd]^{f'}\\
B&\ar[l]^\beta \hat B\ar[r]_{\beta'}&B'
}
$$
is a homotopy commutative diagram in $\mathcal C$.
If $\hat A$ is a  cofibrant object and if
$(\beta,\beta')\colon\hat B\to B\oplus B'$ is a fibration (in other words a surjection in
CDGA or $R$-DGMod),
then there exists a morphism
$\hat f\colon\hat A\to\hat B$ such that $\overline f\simeq \hat f$ and
the diagram
$$
\xymatrix{
&\hat A\ar[ld]_f\ar[d]^{\hat f}\ar[rd]^{f'}\\
B&\ar[l]^\beta \hat B\ar[r]_{\beta'}&B'
}
$$
is strictly commutative.
\end{lemma}
\proof
Consider the following solid arrow diagram
$$
\xymatrix{
\hat A \ar[d]_{i_0} \ar[r]^{\overline f} & \hat B \ar[d]^{(\beta,\beta')} \\
Cyl\hat A \ar[r]_H \ar@{-->}[ur]^G & B \oplus B' }
$$
where $H$ is a homotopy between $(\beta,\beta')\overline f$ and
$(f,f')$. The map $i_0$ is always a weak equivalence and here is a
cofibration since $\hat A$ is cofibrant. Also $(\beta,\beta')$ is a
fibration by hypothesis and hence there is a lift $G$ by Lemma \ref{lemma-tk-lifting}(i).
Let $\hat f=Gi_1$. This makes the diagram of the conclusion commute since
$Hi_1=(f,f')$ and $Gi_1\simeq Gi_0=\overline f$.
\cqfd

\begin{lemma}
\label{lemma-tk-commonmodel}
Suppose in CDGA or $R$-DGMod
that $f\colon A\to B$ is a model of
$f'\colon A'\to B'$. If we are in CDGA then
also assume that
$H^0(A)=H^0(B)=H^0(A')=H^0(B')=\BQ$.
Then there exists a cellularly cofibrant $\hat A$,
a cellular cofibration
$\hat f\colon\xymatrix@1{\hat A\,\,\ar[r]&\hat B}$,
and weak equivalences $\alpha\colon \hat A \quism A$,
$\alpha'\colon \hat A \quism A'$
$\beta\colon\hat B\quism B$
and  $\beta'\colon\hat B\quism B'$  such that
$(\beta,\beta')\colon\hat B\to B\oplus B'$ and
$(\alpha,\alpha')\colon \hat A \to A \oplus A'$ are surjective
and the following  diagram is strictly commutative
$$
\xymatrix{
\ar[d]_f A&
\ar[d]^{\hat f}
\hat A\ar[l]^\simeq_{\alpha}\ar[r]_\simeq^{\alpha'}&
\ar[d]^{f'} A'\\
B&
\hat B\ar[l]^\simeq_{\beta}\ar[r]_\simeq^{\beta'}&
B'.
}
$$
In addition the isomorphisms $H^*(\alpha')H^*(\alpha)^{-1}\in \Hom(H^*(A),H^*(A'))$
and $H^*(\beta')H^*(\beta)^{-1}\in \Hom(H^*(B),H^*(B'))$ are the same
as the one determined by the original string of weak equivalences
making $f$ a model of $f'$.
\end{lemma}
\proof
If we let $I$ be the category $(\bullet \rightarrow \bullet)$ with two objects and
one non-identity map, then $\cal C^I$ is the category of maps in $\cal C$.
According to \cite[Section 5.1]{Hovey} for any model category $\cal C$,
$\cal C^I$ can be given a model structure
such that if $A_1 \rightarrow B_1$ and $A_2\rightarrow B_2$ are in $\cal C^I$ then a map
between them
$$
\xymatrix
{
A_1 \ar[d]\ar[r]^{\alpha} & A_2 \ar[d] \\
B_1 \ar[r]_{\beta} & B_2
}
$$
is a weak equivalence if $\alpha$ and $\beta$ are both weak equivalences and
a fibration if both $\alpha$ and $\beta$ are fibrations.
\par
Since $f$ is a model of $f'$ they are connected by a sequence of weak equivalences,
$$
\xymatrix
{
f=f_0 & g_0 \ar[r]^{\simeq} \ar[l]_-{\simeq} & f_1 &
\cdots \ar[l]_{\simeq} & g_{n-1} \ar[r]^-{\simeq} \ar[l]_{\simeq} & f_n=f'
}
$$
with each $g_i,f_i\in {\cal C}$. We can factor
$g_i\rightarrow f_i\oplus f_{i+1}$ as an acyclic cofibration followed by a fibration.
and therefore we can assume that all of the maps in the sequence are fibrations
(since they are surjections when evaluated at each object of $I$). Also if $g_i$ was cofibrant we
replaced it with something cofibrant.
\par
Now let $\hat f\stackrel{\simeq}{\rightarrow} f_0$ be a weak equivalence such that
$\hat f$ is cofibrant. Since they are fibrations, we can lift along the
maps of the sequence to get a diagram
$$
\xymatrix
{
f & \hat f \ar[r]^{\simeq} \ar[l]_{\simeq} & f'
}
$$
As before we can assume that the map $\hat f \rightarrow f\oplus f'$ is a fibration.
So we get a diagram
$$
\xymatrix
{
\ar[d]_f A&
\ar[d]^{\hat f}
\hat A\ar[l]^\simeq_{\alpha}\ar[r]_\simeq^{\alpha'}&
\ar[d]^{f'} A'\\
B&
\hat B\ar[l]^\simeq_{\beta}\ar[r]_\simeq^{\beta'}&
B'.
}
$$
Since ${\hat f} \rightarrow f\oplus f'$ is a fibration,
$(\alpha,\alpha')\colon \hat A \rightarrow A\oplus A'$ and
$(\beta,\beta')\colon \hat B \rightarrow B\oplus B'$
are surjections. Finally
it follows by a diagram chase that
the isomorphisms $H^*(\alpha')H^*(\alpha)^{-1}\in \Hom(H^*(A),H^*(A'))$
and $H^*(\beta')H^*(\beta)^{-1}\in \Hom(H^*(B),H^*(B'))$ are the same
as the one determined by the original string of weak equivalences
making $f$ a model of $f'$.
\cqfd
\subsection{Homotopy pullbacks}
\label{sec-hompull}
Homotopy pullbacks are homotopy invariant versions of pullbacks.
They exist in any closed model category (\cite[Section 10]{DwyerSpalinski}).
We now give a description
that is adapted to our
categories CDGA and $R$-DGMod.
Let
$f:D\rightarrow C$ be a map in
CDGA or $R$-DGMod and
$
\xymatrix
{
D\ar[r]^i & D' \ar[r]^{f'} & C
}
$
be a factorization into an acyclic cofibration followed by a fibration.
Then the {\em homotopy pullback} of the diagram
$$
\xymatrix
{
B\ar[r]^g & C & D \ar[l]_f
}
$$
is the pullback of
$$
\xymatrix
{
B\ar[r]^g & C & D' \ar[l]_{f'}.
}
$$
There is a map induced by $id_B$, $id_C$ and $i$ from the pullback to the homotopy
pullback. Both the homotopy pullback and this induced map are unique up to
homotopy.
A homotopy pullback in CGDA or RDG-Mod gives rise to a Mayer-Vietoris sequence, corresponding
to the fact that pushouts of spaces have such sequences and $\Apl(\_)$ takes homotopy
pushouts to homotopy pullbacks.
The following lemma states some standard facts about
homotopy pullbacks.
\begin{lemma}
\label{hopo}
Let
\begin{equation}
\label{poho}
\xymatrix
{
B_1\ar[r]\ar[d] & C_1\ar[d] & D_1 \ar[l]\ar[d]\\
B_2\ar[r] & C_2 & D_2 \ar[l]
}
\end{equation}
be a commutative diagram in CDGA or $R$-DGMod.
Let $E_i$ be the pullback of
$
\xymatrix
{
B_i\ar[r] & C_i & D_i \ar[l]
}
$
and $E'_i$ its homotopy pullback.
\begin{itemize}
\item[(i)] There is an induced map $E'_1\rightarrow E'_2$ such that
the following diagram of induced maps commutes
$$
\xymatrix
{
E_1 \ar[d]\ar[r] & E_1' \ar[d] \\
E_2 \ar[r] & E_2'.
}
$$
\item[(ii)] If the vertical maps in (\ref{poho}) are weak equivalences then the induced map
$E_1' \rightarrow E_2'$ is a weak equivalence.
\end{itemize}
\end{lemma}
\proof
Part $(i)$ follows from the Lifting Lemma and
Part $(ii)$ follows from Mayer-Vietoris.
\cqfd
The following lemma follows directly from Lemma \ref{hopo}.
\begin{lemma}
\label{inducted}
Suppose that in CDGA or $R$-DGMod, the diagram
$$
A_1 \rightarrow B_1 \rightarrow C_1 \leftarrow D_1
$$
is a model of the diagram
$$
A_2 \rightarrow B_2 \rightarrow C_2 \leftarrow D_2.
$$
Let $E_i$, $F_i$ be the homotopy pullbacks of
$A_i \rightarrow C_i \leftarrow D_i$
and
$B_i \rightarrow C_i \leftarrow D_i$
respectively. Then
the induced map $E_1 \rightarrow F_1$ is a model of
the induced map $E_2 \rightarrow F_2$.
\end{lemma}
\subsection{Making a homotopy commutative diagram strictly commutative}
\label{sec-strict}
In general if we have a homotopy commutative diagram we cannot always replace it with a strictly
commuting one. However in the following particularly simple case we can.
Recall that weakly equivalent diagrams are connected by a sequence of weak equivalences and weak
equivalences between diagrams are strictly commuting diagrams (see Section \ref{sec-cat}).
\begin{lemma}\label{gorgle}
Assume that we have a homotopy commutative diagram in CDGA or $R$-DGMod
\begin{equation}\label{cat}
\xymatrix
{
A'\ar[d]^{\simeq} \ar[r]^{f'} & B' \ar[d]^{\simeq} &
C' \ar[d]^{\simeq} \ar[l]_{g'} & D' \ar[l]_{h'}\ar[d]^{\simeq}\\
A \ar[r]_{f} & B & C \ar[l]^{g} & D \ar[l]^{h}
}
\end{equation}
such that all vertical arrows are weak equivalences. Then the
diagrams that make up the top and bottom row of (\ref{cat}) are
weakly equivalent.
\end{lemma}
\proof
Replacing the top row by something weakly equivalent we can assume
that the objects in the top
row are cofibrant and the middle two vertical arrows are fibrations.
Consider first the left square of Diagram (\ref{cat}):

$$
\xymatrix
{
A' \ar[d]_{\alpha} \ar[r]^{f'} & B' \ar[d]^{\beta} \\
A \ar[r]_f & B.
}
$$
Since $\beta$ is a surjection and $A'$ is cofibrant the map
$f\alpha$ can be lifted through $\beta$ to get a new map
$\overline f:A' \rightarrow B'$. Also by Lemma \ref{lemma-tk-lifting} there is a homotopy
$H\colon Cyl(A') \rightarrow B'$ between $f'$ and $\overline f$.
So the following diagram is strictly commutative with all
horizontal arrows weak equivalences:
\begin{equation}\label{equ-1}
\xymatrix
{
A' \ar[d]\ar[r]^{f'} & B'\ar[d]^= \\
Cyl(A') \ar[r]_H & B'\\
A' \ar[u] \ar[r]_{\overline f} \ar[d]_{\alpha} & B' \ar[u]_= \ar[d]^{\beta}\\
A \ar[r]_{f} & B.
}
\end{equation}
Similarly, we can replace the center and right squares by a
strictly commuting diagram
\begin{equation}\label{equ-2}
\xymatrix
{
B' \ar[d]_= & C' \ar[l]_{g'} \ar[d] & D' \ar[l]_{h'} \ar[d]\\
B' & Cyl(C')\ar[l]_{H'} & Cyl D' \ar[l]_{Cyl h'}\\
B' \ar[u]^= \ar[d]_= & C' \ar[u] \ar[d]_=
\ar[l]_{\overline g} & D' \ar[u] \ar[d] \ar[l]^{h'}\\
B' & C' \ar[l]_{\overline g} & Cyl D' \ar[l]_{H''}\\
B' \ar[u]^= \ar[d]_{\beta} & C'\ar[u]^= \ar[d] \ar[l]_{\overline g} &
D' \ar[l]_{\overline h} \ar[u] \ar[d]\\
B & C\ar[l]^{g} & D \ar[l]^{h}
}
\end{equation}
where $Cyl h'$ is any lift in the following solid arrow diagram
$$
\xymatrix
{
D' \coprod D' \ar[d] \ar[r] & Cyl (C') \ar[d] \\
Cyl D' \ar@{-->}[ur]^{Cyl h'} \ar[r]_{H''} & C'
}
$$
Notice that  a lift exists since in this diagram the left map is a cofibration and the
right an acyclic fibration.
Now glue the bottom of Diagram (\ref{equ-1}) to the top
of Diagram (\ref{equ-2}) to get a sequence of weak equivalences
connecting the left and right columns of the new diagram. This
completes the proof of the lemma.
\cqfd
\subsection{Mapping cones}
\label{sec-mapcon}
\begin{defin}
\label{def-suspension}
Let $k\in \BZ$.
The {\em $k$-th suspension} of an $R$-dgmodule $M$ is the $R$-dgmodule
$s^kM$ defined by
\begin{itemize}
\item $(s^kM)^j\cong M^{k+j}$ as vector spaces for $j\in\BZ$ and this isomorphism
is denoted by $s^k$,
\item $r\cdot (s^kx)=(-1)^{|r||k|}s^k(r\cdot x)$ for $x\in M$ and $r\in R$,
\item $d(s^kx)=(-1)^ks^k(dx)$ for $x\in M$
\end{itemize}
\end{defin}

\begin{defin}
\label{def-mappingcone}
Let $R$ be a CDGA.
The {\em mapping cone} of an $R$-dgmodule morphism
$f\colon A\to B$ is the $R$-dgmodule $(B\oplus_f sA,d)$ defined as follows:
\begin{itemize}
\item $B\oplus_f sA=B\oplus sA$ as graded $R$-modules,
\item $d(b,sa)=(d_B(b)+f(a),-sd_A(a))$
for $a\in A$ and $b\in B$.
\end{itemize}
\end{defin}
Recall that a mapping cone gives rise to a long exact cohomology sequence, therefore
the following two lemmas follow easily from the five lemma.
\begin{lemma}
\label{lemma-hmtpy-cofibre}
Let
$$
\xymatrix@1{0\ar[r]&A\ar[r]^i&B\ar[r]^p&C\ar[r]&0
}
$$
be a short exact sequence of $R$-dgmodules.
Then there is a quasi-isomorphism of $R$-dgmodules
$$
p\oplus0\colon B\oplus_i sA\quism C.
$$
\end{lemma}
\begin{lemma}
\label{lemma-tk-conehmtpy}
Let
\begin{equation}
\label{diag-tk-conehmtpy}
\xymatrix{
A\ar[r]^{f}\ar[d]_\alpha&B\ar[d]_\beta\\
A'\ar[r]^{f'}&B'
}
\end{equation}
be a homotopy commutative diagram of $R$-dgmodules and let
$$
h\colon sA\to B'
$$
be an $R$-dgmodule homotopy between $\beta f$ and $f'\alpha$ as in Lemma \ref{hom-char}.
Then there exists a commutative diagram
$$
\xymatrix{ 0\ar[r]&B\ar@{^(->}[r]\ar[d]^{\beta}&
B\oplus_fsA\ar[r]\ar[d]^{\beta\oplus_hs\alpha}&
sA\ar[r]\ar[d]^{s\alpha}&0\\
0\ar[r]&B'\ar@{^(->}[r]& B'\oplus_{f'}sA'\ar[r]& sA'\ar[r]&0 }
$$
of $R$-dgmodules
in which each line is a short exact sequence and
where ${\beta\oplus_hs\alpha}$ is the $R$-dgmodule morphism defined by
$$
{(\beta\oplus_hs\alpha)}(b,sa)=(\beta(b)+h(sa),s\alpha(a)).
$$
Moreover if $\alpha$ and $\beta$ are quasi-isomorphisms then so is
$\beta\oplus_hs\alpha$. If the diagram $(\ref{diag-tk-conehmtpy})$
is strictly commutative then we can take $h=0$ and
$\beta\oplus_hs\alpha=\beta\oplus s\alpha$.

\end{lemma}
\subsection{Sets of homotopy classes in $R$-DGMod}
\label{sec-homclass}
Let $R$ be a CDGA and  $M$ and $N$ be $R$-dgmodules with differentials
$d_M$ and $d_N$ respectively.
Forgetting about the differentials we denote by
$\Hom^i_R(M,N)$
the vector space of $R$-module maps from
$M$ to $N$ raising degree by $i$. We can put a differential $D$
on $\Hom_R(M,N)=\oplus_{i\in \BZ} \Hom^i_R(M,N)$ as follows.
For $f\in \Hom^i_R(M,N)$ and $m\in M$ define
$$
(Df)(m)=d_N(f(m))+(-1)^i f(d_M(m)).
$$
$D$ is of degree $+1$ and turns $\Hom_R(M,N)$ into
a chain complex.
It is easy to check that the cycles are exactly the chain maps.
Let
\begin{equation}\label{not_fun}
\rho:H^0(\Hom_R(M,N))\rightarrow[M,N]_R
\end{equation}
be any map which when
restricted to cycles in $\Hom^0(M,N)$ gives their equivalence class
in $[M,N]_R$.
Since a map is a cycle if and only if it is
a chain map and a boundary if and only if it is homotopic to $0$,
if $M$ is cofibrant then $\rho$ is an isomorphism.
\par
We consider $R$ and $H^n(R)$ as chain complexes with $0$ differential, with
$R$ concentrated in degree $0$ and $H^n(R)$ concentrated in degree $n$.
Let $\overline \epsilon:R \rightarrow s^{-n}H^n(R)$ be any chain
map that induces the identity isomorphism after taking $H^n$. In
other words such that $\overline \epsilon$ maps dimension $n$
cocycles to their equivalence classes in cohomology. The condition
that $\overline \epsilon$ is a chain map of course implies that
all other dimensions go to $0$. Any such $\overline \epsilon$
gives us a map
\begin{equation}\label{equ-mape}
\epsilon: R \rightarrow \Hom_{\BQ}(R,s^{-n}H^n(R)),\ r \mapsto
\overline \epsilon(r\cdot \_)
\end{equation}
It is straightforward to check that $\epsilon$ is a degree $0$ map
of $R$-dgmodules. Note that in the general case
where we don't assume commutativity, for $s\in R$ we would have
$\epsilon(r)(s)=(-1)^{|r||s|}\overline\epsilon(sr)$. Since we
are working in the graded commutative situation
$(-1)^{|r||s|}sr=rs$. Thus we get the same formula in our
situation but with fewer signs.
\begin{defin}
A connected graded algebra $R$ is a {\em Poincar\'e duality
algebra} of formal dimension $n$ if $\epsilon$ is a
quasi-isomorphism.
\end{defin}
For $M$ a left $R$-dgmodule, $N$ a right $R$-dgmodule and $L$ a
chain complex, define
\begin{equation}\label{fun}
\phi:\Hom_R(M, \Hom_{\BQ}(N,L)) \rightarrow \Hom_{\BQ}(N\otimes_R M,
L)
\end{equation}
by $\phi(f)(n\otimes m)=(-1)^{|n||m|}f(m)(n)$. Clearly $\phi$ is a degree $0$
isomorphism of chain complexes.
\begin{lemma}
\label{lemma-tk-PD} Let $R$ be a CDGA such that
$H^*(R)$ is a Poincar\'e duality algebra of formal dimension $n$.
Let $P$ be an $R$-dgmodule. Then the map
$$
H^n:[P,R]_R\rightarrow \Hom(H^n(P),H^n(R)),\ [f]\mapsto H^n(f).
$$
is an isomorphism of $\BQ$-modules.
\end{lemma}
\proof By the definition of $[P,R]_R$ we may assume that $P$ is
cofibrant. Since $R$ is a Poincar\'e duality algebra
and $P$ is cofibrant we have a quasi-isormorphism
$\epsilon_*\colon H^0(\Hom_R(P,R))
\rightarrow H^0(\Hom_R(P,\Hom_{\BQ}(R,s^{-n}H^n(R))$
induced by Equation (\ref{equ-mape}).
Also $\phi$ and $\rho$ from equations (\ref{not_fun}) and (\ref{fun}) are isomorphisms.
So we get a string of isomorphisms.

\begin{eqnarray*}
[P,R]_R & \stackrel{\rho^{-1}}{\rightarrow} & H^0(\Hom_R(P,R)) \\
& \stackrel{\epsilon_*}{\rightarrow} & H^0(\Hom_R(P,
\Hom_{\BQ}(R,s^{-n}H^n(R))\\
& \stackrel{\phi_*}{\rightarrow} & H^0(\Hom_{\BQ}(R\otimes_R P,s^{-n}H^n(R)))\\
& = & H^0(\Hom_{\BQ}(P,s^{-n}H^n(R))) \\
& = & \Hom_{\BQ}(H^n(P), H^n(R))
\end{eqnarray*}
 It is straightforward to check that the composition of these
maps is the same as the map induced by taking $H^n$.
\cqfd
\begin{lemma}
\label{lemma-tk-nullhmtpy} Let $P$ and $X$
be $R$-dgmodules. If $H^{<r}(P)=0$ and $H^{\geq r}(X)=0$ for
some $r\in\BZ$ then $[P,X]_R=0$.
\end{lemma}
\proof This is a straightforward obstruction argument. See Lemma \ref{hom-char}.
\cqfd

\subsection{Turning $R$-dgmodule structure into CDGA structure}
\label{sec-turn}
\begin{defin}
Let $R$ be a commutative graded algebra (CGA)
and let $X$ be a right $R$-module.
Then, the {\em semi-trivial CGA} structure on $R\oplus X$ is the multiplication
$$
\mu\colon (R\oplus X)\otimes(R\oplus X)\to (R\oplus X)
$$
defined, for homogeneous elements $r,r'\in R$ and $x,x'$ in $X$, by
\begin{itemize}
\item $\mu(r\otimes r')=r.r'$
\item $\mu(x\otimes r')=(-1)^{|r'||x|}r'\cdot x$
\item $\mu(r\otimes x')=r\cdot x'$
\item $\mu(x\otimes x')=0$.
\end{itemize}
\end{defin}
Note that if $d$ is a differential of $R-$modules on $R\oplus X$,
it is in general not true that the multiplication $\mu$ defined above
defines a CDGA structure on $(R\oplus X,d)$ because the Leibnitz
rule is not necessarily satisfied. However, we have the following:
\begin{prop}
\label{propdef-semitrivCDGA}
Let $R$ be a CDGA, let $Q$ be an $R$-dgmodule and let $f\colon Q\to R$
be a morphism of $R$-dgmodules.
If either
\begin{enumerate}
\item $f\colon Q\to R$ is the inclusion of an ideal with the $R$ module structure
on $Q$ given by multiplication in $R$, or
\item there exist $p\in\BN$ such that $Q^i=0$ for $i<p$ and $i\geq 2p$
\end{enumerate}
hold then the semi-trivial CGA structure on the mapping cone
$R\oplus_fsQ$ defines a CDGA structure.
\end{prop}
\proof
To see that the semi-trivial CGA structure defines a CGDA structure we only have to check that
the Leibniz law holds. By definition
$(r,sq)\cdot (r',sq')=(rr',s(rq'+(-1)^{|q||r'|}r'q))$. Since the product is bilinear
and the differential is linear, we
can check Leibniz on terms of the form $(r,0)\cdot (r',0)$,
$(r,0)\cdot(0,sq')$, $(0,sq')\cdot(r',0)$ and $(0,sq)\cdot(0,sq')$.
Since $R$ satisfies Leibniz so do terms of the first type. Since $Q$ is an $R$-module and $f$ is
an $R$-dgmodule map, it follows directly from the definition of the differential that Leibniz is
satisfied for terms of the second and third type.
\par
Next consider terms of the fourth type. On one side of the Leibniz equation we have
$d[(0,sq)\cdot(0,sq')]=d(0)=0$.
On the other side we have,
$$
d(0,sq)\cdot(0,sq')+(-1)^{|sq|}(0,sq)\cdot d(0,sq')
$$
$$
=(f(q),-sdq)\cdot(0,sq')+(-1)^{|q|-1}(0,sq)\cdot(f(q'),-sdq')
$$
$$
= f(q)sq'+(-1)^{|q|-1}(-1)^{|q'|(|q|-1)}f(q')sq
$$
In Case $2$ this term is zero for degree reasons. In Case $1$ we continue using the fact that
the multiplication is given by that in $R$,
$$
= (-1)^{|q|}s(qq')+(-1)^{|q'|}(-1)^{|q|-1}(-1)^{|q'|(|q|-1)}s(q'q)
$$
$$
=(-1)^{|q|}s(qq')+(-1)^{|q||q'|}(-1)^{|q'|}(-1)^{|q|-1}(-1)^{|q'|(|q|-1)}s(qq')
$$
$$
=(-1)^{|q|}s(qq')-(-1)^{|q|}s(qq')=0
$$
So in both cases the product is zero and Leibniz holds.
\cqfd
\begin{defin}\label{def-semitrivCDGA}
The CDGA structure of the last proposition is called the {\em semi-trivial CDGA structure}
on the mapping cone $R\oplus_fsQ$.
\end{defin}

\section{Thom class and the shriek map}
\label{section-modelcompl}
Throughout this section we fix a connected
oriented smooth manifold $W$ of dimension $n$
and an oriented connected closed smooth submanifold $V\subset W$ of dimension $m$
and codimension $r=n-m$.
We denote this embedding by
$$
f\colon V\hookrightarrow W,
$$
and the
orientation classes by $u_V\in H^m(V;\BZ)$ and $u_W\in
H^n(W;\BZ)$.
\par
In this section we introduce the ingredients needed for
describing a model of $W\setminus V$ and we build a first version of such a model.
In Section \ref{sec-Thom} we review the description of the Thom isomorphism and
the normal bundle of $V$ in $W$ as a tubular neighborhood $T$
(\cite[Chapters 10 and 11]{MS}). We also describe $W$ as a pushout
(Diagram \ref{decodiagpo}), for which we find a rational model
in Section \ref{section-modelblowup}.
We associate a Thom class $\overline\theta\in H^*(T,\bd T;\BZ)$ to the
normal bundle of the embedding $f$ compatible in a certain way with the orientations
of $V$ and $W$ (Lemma \ref{shape}).
In Section \ref{sec-shr} and Section \ref{sec-unish}
we introduce our chain level version of the shriek map and
prove its existence and uniqueness.
Finally in Section \ref{sec-preap} we
give a first CDGA model of Diagram (\ref{decodiagpo}) below using
mapping cones (Lemma \ref{lemma-Aplconemodel}). This is a very specific model but
an intermediate step in constructing the CDGA model
of this diagram based on any model of $f$, which is done
in the Lemma \ref{igog}.

\subsection{Thom class and orientations}
\label{sec-Thom}
For this subsection alone our coefficients for cohomology are $\BZ$.
By the tubular neighborhood theorem \cite[Theorem 11.1]{MS}
the normal bundle $\nu$ of the embedding $f$ is diffeomorphic to some open
neighborhood $T'$ of $V$ in $W$.
The associated normal disk bundle $D\nu$ can be identified with a compact
manifold $T\subset T'$ such that the inclusion makes $V$ a
strong deformation retract of $T$.
Let $B:=\overline{W\setminus T}$ be the closure of
the complement of $T$. Then $B$ is a compact manifold with boundary
$\del B=\del T=B\cap T$ and $B\cup T=W$.
In other words we have the following
push-out
\begin{equation}\label{decodiagpo}
\xymatrix{
\del T\ \ar@{^(->} [r]^k\ar@{^{(}->}[d]_i&
B\ar@{^(->}[d]_l\\
T\ \ar@{^(->}[r]^j&W
}
\end{equation}

This pushout is also a homotopy pushout since the inclusion
$\del T \rightarrow T$ is a cofibration \cite[Theorem 13.1.10 and Proposition 13.5.4]{Hir}.
It is clear that the projection of the normal bundle $D\nu$
on $V$ defines a homotopy inverse
$$
\pi\colon T\stackrel{\simeq}{\rightarrow} V
$$
of the inclusion
such that the composition
$$
V\subset T\stackrel{j}{\hookrightarrow}W
$$
is equal to $f$. Also $B\simeq W\setminus V$.

The projection map $(T,\del T)\to V$
is an orientable $(D^r,S^{r-1})$-bundle.
Let $*\in V$ be some point and consider
$(D^r,S^{r-1})=(T,\bd T)|_*$
as the restriction of the bundle
over that point.
We thus have an associated inclusion
${\rm inc}\colon(D^r, S^{r-1}) \rightarrow (T, \bd T)$.
Let $\tau_W$ and $\tau_V$ be the tangent bundles of $W$ and
$V$ respectively. We can write
$$
\tau_W|_V\cong \tau_V\oplus \nu
$$
Because $W$ and $V$ are orientable so is $\nu$.
Let $D^m\subset V$ be a closed neighborhood of $*$ with boundary
$S^{m-1}$ and $D^n$ be a neighborhood of $*$ in $W$ with
boundary $S^{n-1}$.
The given orientations $u_V$ and $u_W$ induce orientations
$$
u_{D^m}\in H^m(D^m,S^{m-1}) \cong H^m(V,V\setminus *)\cong H^m(V)
$$
and
$$
u_{D^n}\in H^n(D^n,S^{n-1}) \cong H^n(W,W\setminus *)\cong H^n(W).
$$
By \cite[Lemma 11.6]{MS} these correspond to orientations on
$\tau_W$ and $\tau_V$.
Suppose an orientation
\begin{equation}\label{equ-ornormal}
u_{D^r}\in H^r(D^r,S^{r-1})
\end{equation}
on $\nu$ has also been given.
We say that $u_V$, $u_W$ and $u_{D^r}$ are {\em compatible} if they satisfy the formula
\begin{equation}\label{equ-compat}
u_{D^n}=u_{D^m}\times u_{D^r}\in
H^{m+r}((D^m,S^{m-1})\times (D^r,S^{r-1})) =H^n(D^n,S^{n-1}).
\end{equation}
By \cite[Theorem 9.1]{MS} there
exists a Thom class $\bar\theta\in H^r(T,\del T)$ such that
$\overline \theta|_*=u_{D^r}$, in other words such that
\begin{equation}\label{equ-incT}
{\rm inc}^*(\overline{\theta})=u_{D^r}.
\end{equation}
It also
induces a Thom isomorphism:
\begin{equation}\label{equ-Thom}
\_\cup\bar\theta\colon H^*(T)\iso H^{*+r}(T,\del T).
\end{equation}

The following formula summarizes the link
between the Thom class $\overline \theta$ and the fixed orientation
classes $u_V$ and $u_W$.
\begin{lemma}\label{shape}
Suppose that $u_v$, $u_W$ and $u_{D^r}$ are compatible.
Then the following sequence of arrows sends $u_V$ to $u_W$
$$
\xymatrix
{
H^*(V)\ar[r]^{\cong}_{\pi^*} &
H^*(T) \ar[r]^-{\_\cup \overline \theta} & H^{*+r}(T,\bd T)
 & H^{*+r}(W,B) \ar[l]_{\cong}^{\rm restriction} \ar[r] & H^{*+r}(W).
}
$$
\end{lemma}
\proof
The lemma follows from Equation (\ref{equ-compat}) and the commutativity of the following diagram
$$
\xymatrix
{
H^*(T) \ar[d] \ar[r]^{\_\cup \overline \theta} & H^{*+r}(T,\bd T) \ar[d]
 & H^{*+r}(W,B) \ar[d] \ar[r]\ar[l]_{\cong} & H^{*+r}(W) \ar[d] \\
H^*(D^m,S^{m-1})\ar[r]_{\_\times u_{D^r}} & H^{*+r}(D^n,S^{n-1})\ar@{=}[r] &
H^{*+r}(D^n,S^{n-1})\ar@{=}[r] & H^{*+r}(D^n,S^{n-1}).
}
$$
\cqfd
\subsection{The shriek map}
\label{sec-shr}
To construct a CDGA-model of $W\setminus V$ we will need an analogue
of the shriek map (or Gysin map, or transfer map) which we recall now.
Let $\mu_V\in H_m(V,\BZ)$ and $\mu_W\in H_n(W,\BZ)$ be homology classes dual
to $u_V$ and $u_W$.
The classical {\em cohomological shriek map} (see \cite[VI.11.2]{Bredon},
also called Gysin map, pushforward map or umkehr map) is a map
$$
f^!\colon s^{-r}H^*(V,\BZ) \cong H^{*-r}(V,\BZ) \to H^*(W,\BZ)
$$
By \cite[VI.14.1]{Bredon} $f^!$ is a morphism of $H^*(W,\BZ)$-modules
and induces an isomorphism in degree $n$ sending
$s^{-r} u_V$ to $u_W$. These two facts characterize $f^!$ as can also be seen
using Lemma \ref{lemma-tk-PD}.

An important ingredient in our model of $W\setminus V$ will be the
following analogue of the shriek map at the level of models.
Following our standard abuse of notation (Section \ref{sec-abu})
denote also by $u_V\in H\Apl(V)$ and $u_W\in H\Apl(W)$ the images
of the orientation classes through the isomorphism between
$H\Apl(\_)$ and $H^*(\_)$ (Equation (\ref{equ-isocan})).
\begin{defin}
\label{def-shriek} Let $\phi\colon R\to Q$ be a CDGA-model of
$\Apl(f)\colon \Apl(W)\to\Apl(V)$. A {\em shriek map} associated
to $\phi$ is any $R$-dgmodule morphism
$$
\phi^!\colon s^{-r}Q\to R
$$
such that $H^n(\phi^!)(s^{-r}u_V)=u_W$ .
\end{defin}
\begin{prop}
Let $\phi\colon R\to Q$ be a CDGA-model of $\Apl(f)\colon
\Apl(W)\to\Apl(V)$. An  $R$-dgmodule morphism $\tilde\phi\colon s^{-r}Q\to R$
is a shriek map if and only if
$$
H(\tilde \phi)\colon s^{-r}H^*(V)\cong H(s^{-r}Q) \rightarrow H(R)\cong H^*(W)
$$
is the classical cohomological shriek map $f^!$.
\end{prop}
\proof
Assume that $H(\tilde \phi)$ is the cohomological shriek map $f^!$.
We know that
$H_*(f)(u_V\cap \mu_V)=H_*(f)(1)=u_W\cap \mu_W$. Thus by
\cite[VI.14.1]{Bredon}
$f^!(s^{-r}u_V)=u_W$ and so $\tilde \phi$ is a shriek map.

Next assume that $\tilde \phi$ is a shriek map, then $H^n(\phi^!)(s^{-r}u_V)=u_W$
by definition. We also know that the cohomological shriek map $f^!$ satisfies
$f^!(s^{-r}u_V)=u_W$. Letting $(H(R),0)$ be the CGDA and
$H^*(Q)$ the $H^*(R)$-dgmodule in Lemma \ref{lemma-tk-PD} we
see that $H^r(\phi^!)\simeq f^!$ in the category $H^*(R)$-DGMod. Hence they must be equal
on homology, which completes the proof of the other direction.
\cqfd
\subsection{Uniqueness and existence of the shriek map}
\label{sec-unish}
\begin{prop}\label{unish}
The shriek map is unique up to homotopy. More precisely let
$\phi\colon R \rightarrow Q$ be a CDGA model of $\Apl(f)\colon
\Apl(W)\to\Apl(V)$. Let $\phi^!, \overline \phi^!\colon  s^{-r}Q
\rightarrow R$ be shriek maps associated to $\phi$, then
$\phi^!=\overline \phi^!$ in $[s^{-r}Q, R]_R$.
\end{prop}
\proof
We know that $H^n(\phi^!)(s^{-r} u_V)=H^n(\overline \phi^!)(s^{-r} u_V)$
since both $\phi^!$ and $\overline \phi^!$ are shriek maps. Thus as
$H^n(s^{-r}Q)$ is one dimensional, $H^n(\phi^!)=H^n(\overline \phi^!)$ in
$\Hom(H^n(s^{-r}Q),H^n(R))$ and so Lemma \ref{lemma-tk-PD} implies that
$\phi^!=\overline \phi^!$ in $[s^{-r}Q, R]_R$.
\cqfd
The next proposition and its proof show how to associate to any CDGA model of an
embedding $f:V\rightarrow W$, a suitable CDGA model together with a shriek map which
will be used in our model for $W\setminus V$ and for the blow-up.
\begin{prop}
\label{prop-shriekexists}
Assume that $H^1(f)$ is injective and that
$n\geq m+2$. Suppose $\phi'\colon R'\to Q'$ is a CDGA-model
of the embedding $V\stackrel{f}{\hookrightarrow}W$.
Then we can construct from $\phi'$ another CDGA model $\phi\colon R\to Q$
of $f$ together with a shriek map $\phi^!\colon s^{-r}Q\to R$, such that
\begin{itemize}
\item $R$ and $Q$ are connected, that is, $R^0=Q^0=\BQ$,
\item $R^{\geq n+1}=0$, and
\item $Q^{\geq m+2}=0$.
\end{itemize}
\end{prop}
\proof By taking a minimal Sullivan model of $R'$,
$\beta\colon\tilde R=\Lambda X\quism R'$ and a minimal relative
Sullivan model $\tilde \phi\colon\tilde R\cofarrow \tilde Q=\tilde
R\otimes \Lambda Y$ of the composition $\phi'\beta$ we get a new
model $\tilde\phi$ of $f$. Since
$H^0(V)=H^0(W)=\BQ$ and $H^1(f)$ is injective, $\tilde R$ and $\tilde Q$
 are connected.

By \cite[Lemma 14.1]{FHT-book}, $\tilde Q$ is a semi-free $\tilde R$-dgmodule
and so is $s^{-r}\tilde Q$, and hence they are cofibrant.
Since $H^*(\tilde R)\cong H^*(W;\BQ)$
is a connected Poincar\'e duality algebra of formal dimension $n$,
Lemma \ref{lemma-tk-PD} implies that
$$
\left[s^{-r}\tilde Q,\tilde R\right]_{\tilde R}\cong
\Hom(H^n(s^{-r}\tilde Q),H^n(\tilde R)).
$$
From the definition of $[\_,\_]$, $\Hom(A,B)\rightarrow [A,B]$ is surjective if
$A$ is cofibrant, so
since $s^{-r}\tilde Q$ is cofibrant we can
take a representative $\tilde\phi^!\colon s^{-r}\tilde Q\to\tilde R$
of a homotopy class such that
$H^n(\tilde \phi^!)(u_V)=u_W$.
Then $\tilde\phi^!$ is a shriek map associated to $\tilde\phi$.
\par
We next adjust this shriek map so that the dimension conditions
at the end of the proposition are satisfied.
Set
\begin{eqnarray*}
I&=&\tilde R^{\geq n+1}\oplus(\mbox{a complement of
the cocycles in $\tilde R^n)$}\\
J&=&\tilde Q^{\geq m+2}\oplus(\mbox{a complement of
the cocycles in $\tilde Q^{m+1})$}.
\end{eqnarray*}
Since $\tilde R$ is connected, $I$ is an ideal and it is acyclic
because $H^{>n}(\tilde R)=0$. Similarly $J$ is an acyclic ideal in
$\tilde Q$. Define $R=\tilde R/I$ and $Q=\tilde Q/J$. Since $n\geq
m+2$, $\tilde\phi$ induces a map $\phi\colon R\to Q$ which is a
model of $f$ and since $m+1+r\geq n+1$, $\tilde\phi^!$ induces an
$R$-dgmodule morphism $\phi^!\colon s^{-r}Q\to R$ which is a
shriek map. Obviously $R$ and $Q$ are connected and $R^{\geq
n+1}=Q^{\geq m+2}=0$.
\cqfd

\subsection{Preliminaries with $\Apl$}
\label{sec-preap}
Recall the notation of Diagram (\ref{decodiagpo}) at the beginning
of the section. Consider the ladder of maps between short exact sequences of pairs:
$$
\xymatrix
{
\Apl(W,B) \ar[r]^{\iota} \ar[d]^{\rho} &
\Apl(W) \ar[r]^{l^*} \ar[d]^{j^*} & \Apl(B) \ar[d]^{k^*} \\
\Apl(T,\del T) \ar[r]^{\iota'} & \Apl(T) \ar[r]^{i^*} & \Apl(\del T)
}
$$
where
$\iota$
is the kernel of $l^*$ and $\iota'$
is the kernel of $i^*$.
Also $\rho$
is the induced
map which is well defined since Diagram (\ref{decodiagpo})
commutes and $\Apl$ is a functor. As each is a kernel of a CDGA map,
each is a differential ideal -
$\Apl(W,B)$ inherits an $\Apl (W)$-dgmodule structure and
$\Apl(T,\bd T)$ inherits an $\Apl (T)$-dgmodule structure.
Since $j^*$ is an algebra map, $\rho$ is an $\Apl(W)$-dgmodule map.
Note that $l:B\rightarrow W$ and
$i:\bd T \rightarrow T$ are cofibrations
and so $l^*$ and $i^*$ are surjections
\cite[Proposition 10.4 and Lemma 10.7]{FHT-book}. This will be used later.

\begin{lemma}
\label{lemma-AplWB=Tdt}
The restriction map
$$
\rho\colon\Apl(W,B)\rightarrow\Apl(T,\del T)
$$
is a surjective weak equivalence of $\Apl(W)$-dgmodules.
\end{lemma}
\proof
As observed above $\rho$ is an $\Apl(W)$-dgmodule map. It is
a weak equivalence by excision. To see that it is surjective,
let $\alpha\in \Apl(T,\bd T)$. Since $j$ is a cofibration we know
by \cite[Proposition 10.4 and Lemma 10.7]{FHT-book}
that $j^*$ is a surjection. Let $\beta\in \Apl(W)$ be such that
$j^*(\beta)=\iota'(\alpha)$. Since $i^*\iota'=0$, $k^*l^*(\beta)=0$.
So we can extend $l^*(\beta)$ by $0$ on simplices
contained in $T$ and arbitrarily
to the rest of the singular simplices in $W$ (the ones whose image
is contained in neither $B$ or $T$) to get $\beta'\in \Apl(W)$. We know that
$l^*(\beta')=l^*(\beta)$ and $j^*(\beta')=0$. Thus $l^*(\beta-\beta')=0$ and
$j^*(\beta-\beta')=j^*(\beta)=\iota'(\alpha)$. So $\beta-\beta'$ lifts to
$\tilde \beta \in \Apl(W,B)$ such that $\rho(\tilde\beta)=\alpha$. Since
$\alpha$ was arbitrary $\rho$ is surjective.
\cqfd
See \ref{sec-cat} for notation concerning squares.
\begin{lemma}
\label{lemma-Aplconemodel}
Consider the following squares:
$$
{\mathbf F}'=\xymatrix{
\Apl(W)\ar[r]^{j^*}\ar[d]_{\rm inc}&
\Apl(T)\ar[d]_{\rm inc}\\
\Apl(W)\oplus_{\iota}s\Apl(W,B)\ar[r]^-{j^*\oplus s\rho}&
\Apl(T)\oplus_{\iota'}s\Apl(T,\bd T)
}
$$
and
$$
{\mathbf F}=\xymatrix{
\Apl(W)\ar[r]^{j^*}\ar[d]^{l^*}&
\Apl(T)\ar[d]^{i^*}
\\
\Apl(B)\ar[r]^{k^*}&
\Apl(\del T) }
$$
(i) With
the semi-trivial CDGA structure (see Definition \ref{def-semitrivCDGA})
on the mapping cones, $\mathbf F'$ is a CDGA squares.
\par
\noindent
(ii)
$
\Theta_6=\left(
\begin{tabular}{cc}$\id$&$\id$\\$l^* + 0$&$i^* + 0$\end{tabular}
\right)\colon {\mathbf F}'\rightarrow {\mathbf F}
$
is a weak equivalence of CDGA squares.
\end{lemma}
\proof
(i) By Proposition \ref{propdef-semitrivCDGA} (i) there are
semi-trivial CDGA-structures on the mapping cones:
$$
\Apl(W)\oplus_\iota s\Apl(W,B) \quad {\rm and} \quad
\Apl(T)\oplus_{\iota'}s\Apl(T,\del T).
$$

It is then easy to see that $\mathbf F'$
is a CDGA squares.
\par
(ii)
It is immediate that the following diagram is one in CDGA
and is commutative:
\begin{equation}
\label{diag-Aplconemodel}
\xymatrix{
\Apl(W)\oplus_{\iota}s\Apl(W,B)
\ar[r]^-{j^*\oplus s\rho}\ar[d]_{l^* + 0}&
\Apl(T)\oplus_{\iota'}s\Apl(T,\del T)\ar[d]_{i^* + 0}\\
\Apl(B)\ar[r]^{k^*}&
\Apl(\del T).
}
\end{equation}
We need only to prove that the  vertical arrows in
Diagram $(\ref{diag-Aplconemodel})$
are quasi-isomorphisms.
As noted above Lemma \ref{lemma-AplWB=Tdt},
$l^*$ is surjective. Thus
we have a short exact sequence
$$
\xymatrix{
0\ar[r]&\Apl(W,B)\ar[r]^{\iota}&\Apl(W)\ar[r]^{l^*}&
\Apl(B)\ar[r]&0.
}
$$
By Lemma \ref{lemma-hmtpy-cofibre} this implies that the vertical
map $l^* + 0$ is a quasi-isomorphism. The proof that
$i^* + 0$ is a quasi-isomorphism is similar. \cqfd

\begin{lemma}
\label{Thom-qi}
Let $\theta\in \Apl^r(T,\del T)\cap\ker d$ be a representative of
the Thom class
$\bar\theta$ chosen in Section \ref{sec-Thom}. Then multiplication
by $\theta$,
$$
\_\cdot\theta\colon s^{-r}\Apl(T)\quism\Apl(T,\del T)
$$
is a quasi-isomorphism of left $\Apl(T)$-dgmodules.
\end{lemma}
\proof
This follows from the Thom isomorphism of Equation (\ref{equ-Thom}).
\cqfd

\begin{lemma}\label{chirp}
The following sequence of arrows
$$
\xymatrix{
s^{-r}\Apl(V) \ar[r] &
s^{-r}\Apl(T)\ar[r]^{\_\cdot\theta}_{\simeq}&
\Apl(T,\del T)&
\ar[l]_{\rho}^{\simeq}\Apl(W,B)\ar[r]^{\iota}&
\Apl(W)&
}
$$
induces an isomorphism
$$
H^{n-r}(V)\cong H^{n}(W).
$$
that takes $s^{-r}u_V$ to $u_W$.
\end{lemma}
\proof
Apply Lemma \ref{shape} and the natural equivalence
of Equation (\ref{equ-isocan})
between
$H^*(\_)$ and $H(\Apl(\_))$.
\cqfd

%
\section{A model of the complement of a submanifold}
\label{sec-compl}
%
As in Section \ref{section-modelcompl} we suppose we are
given an embedding of closed manifolds $f:V \hookrightarrow W$ as well as orientation
classes $u_V\in H^m(V;\BZ)$ and $u_W\in H^n(W;\BZ)$. Again
$W$ is of dimension $n$ and $V$ of dimension $m$ and codimension $r=n-m$.
We also use the notation from Diagram
(\ref{decodiagpo}), the CDGA maps $\iota$,$\iota'$, and $\rho$ from Section \ref{sec-preap} and
the representative $\theta$ of the Thom class from  Lemma \ref{Thom-qi}.
\par
Assume that $H^1(f)$ is injective and
$\dim W\geq 2\dim V+3$. We fix a
CDGA-model
$$
\phi\colon R\to Q
$$
of $j^*\colon\Apl(W)\to\Apl(T)$ such that $R$ and $Q$ are connected, $R^{\geq n+1}=0$ and
$Q^{\geq m+2}=0$. Notice that $\phi$ is also a model of $\Apl(f)$.
By our standard abuse of notation (see Section \ref{sec-abu}) this
determines orientation classes
$u_W\in H(R)$ and $u_V\in H(Q)$.
We suppose also that we have been given an
associated shriek map
$$
\phi^!\colon s^{-r}Q\to R.
$$
Notice that by Proposition \ref{prop-shriekexists} we can always build such a CDGA model $\phi$ and
shriek map $\phi^!$.
\par
Our aim in this section is to describe a CDGA-model of the map
 $k:\bd T\rightarrow B$
using a CDGA model of the embedding $f$ under the hypotheses that
$n\geq 2m+3$.
In fact we will give a CDGA model of the Diagram (\ref{decodiagpo}) of the last section
(Lemma \ref{igog})
using only the model $\phi$ and the shriek map $\phi^!$.
This extra precision
is necessary to get the model of the blow-up.
\par
The section is organized as follows. In Section \ref{sec-supps}
we fix a common model
$\hat \phi:\hat R\rightarrow \hat Q$
 of $\phi$ and $j^*$. In Section \ref{sec-comsh}
we construct a common model $\hat \phi^!$ of $\phi^!$
and $\iota$. The common model $\hat \phi\hat \phi^!$ of
$\phi\phi^!$ and $j^*\iota$ comes with a $\hat R$-dgmodule structure.
In Section \ref{sect-give} we show that it is homotopic to a
$\hat Q$-dgmodule map $\chi$.
This extension of structure would not be necessary if we just wanted
a model of $\Apl(B)$, but it is necessary to obtain a model of
$\Apl(k)$.
In Section \ref{sec-diagmod} (Lemma \ref{lemma-RQmodel-k}) we show that the
cone on $\phi^!$ is a $\hat R$-dgmodule model of the cone on
$\iota$ and that the cone on $\phi\phi^!$ is a $\hat Q$-dgmodule model
of the cone on $j^*\iota$. We already know from the last section
(Lemma \ref{lemma-Aplconemodel}) that the cones on
$\iota$ and $j^*\iota$ are CDGA models of
$\Apl(B)$ and $\Apl(\bd T)$ respectively. We also construct the maps between
our models that we will need. In Section \ref{sec-extend}
(Lemma \ref{kuif}) we give conditions under which a map between diagrams with certain
dgmodule structure can be extended to a CDGA map.
Next we put everything together and construct a CDGA model
of Diagram (\ref{decodiagpo}), (Lemma \ref{igog}).
Since we are constructing the model of the blow-up we also keep track of the weak
equivalences connecting the models.
\subsection{A common model of $\phi$ and $j^*$}
\label{sec-supps}
Solely for purposes of the proof we use Lemma \ref{lemma-tk-commonmodel}
to fix a
commuting diagram of CDGAs
\begin{equation}
\label{diag-phihat}
\xymatrix{
R\ar[d]_{\phi} &
\hat R\ar[l]_{\alpha}^{\simeq}\ar[r]^-{\alpha'}_-{\simeq}\ar[d]_{\hat\phi}&
\Apl(W)\ar[d]_{j^*}\\
Q &
\hat Q\ar[l]_{\beta}^{\simeq}\ar[r]^-{\beta'}_-{\simeq}&
\Apl(T)
}
\end{equation}
such that $\alpha,\alpha',\beta$ and $\beta'$ are quasi-isomorphisms,
$\hat R$ is cellular, $\hat \phi$ is a cellular cofibration and the maps
\begin{eqnarray*}
(\alpha,\alpha')&\colon&\hat R\to R\oplus\Apl(W)\\
(\beta,\beta')&\colon&\hat Q\to Q\oplus\Apl(T)
\end{eqnarray*}
are surjective.
As a consequence $R$ and
$\Apl(W)$ are $\hat R$-dgmodules and $Q$ and $\Apl(T)$ are $\hat
Q$-dgmodules. Notice also that since $\hat\phi$ is a cellular cofibration
and $\hat R$ and $\hat Q$ are connected, $\hat Q$ is a semifree
$\hat R$-dgmodule.
By the second part of Lemma \ref{lemma-tk-commonmodel} the homology
classes $u_W\in H(R)$ and $u_V\in H(Q)$ corresponding to the
orientations $u_W\in H^*(W)=H(\Apl(W))$ and $u_V\in
H^*(V)\stackrel{\pi^*}{\cong} H^*(T)=H(\Apl(T))$ using the quasi-isomorphisms of
Diagram (\ref{diag-phihat}) are the same as those given
by the original string of weak equivalences that made
$\phi$ a model of $j^*$.
\subsection{A common model of $\phi^!$ and $\iota$}
\label{sec-comsh}
Next we
construct a common model $\hat\phi^!$ of $\phi^!$ and of
$\iota\colon\Apl(W,B)\rightarrow \Apl(W)$ defined in Section \ref{sec-preap}.
Recall also from Section \ref{sec-preap} the definition $\rho$
and the cocycle $\theta$ from Lemma \ref{Thom-qi}.
We consider the following
commutative solid arrow diagram of
$\hat R$-dgmodules
\begin{equation}\label{long}
\xymatrix{
& &  s^{-r} \Apl(T) \ar[dr]_{\simeq}^{\_ \cdot \theta} \\
s^{-r}Q\ar[d]_{\phi^!}&
\ar[l]_{s^{-r}\beta}^{\simeq}s^{-r}\hat
Q\ar@{-->}[r]^{\gamma'}_{\simeq}\ar@{-->}[d]_{\hat\phi^!}
\ar[ur]^{s^{-r}\beta'}_{\simeq} &
\Apl(W,B)\ar[d]_{\iota} \ar[r]_{\rho}^{\simeq} & \Apl(T,\del T) \\
R&
\hat R\ar[l]_{\alpha}^{\simeq}\ar[r]^-{\alpha'}_-{\simeq}&
\Apl(W).
}
\end{equation}
In the next lemma we construct $\hat R$-dgmodule maps $\gamma'$ and
$\hat \phi'$ making the diagram commute. Notice that all
horizontal and diagonal maps are weak equivalences.
Some of the maps above are more than just $\hat R$-dgmodule
maps and this extra structure will be used when we enhance our
$\hat R$-dgmodule structure to get a CDGA structure in Section \ref{sec-extend}.
\begin{lemma}\label{still}
There exists an $\hat R$-dgmodule weak equivalence $\gamma':s^{-r}
\hat Q \rightarrow \Apl(W,B)$ and a $\hat R$-dgmodule map $\hat
\phi':s^{-r} \hat Q \rightarrow \hat R$ making Diagram
(\ref{long}) above commute.
\end{lemma}
\proof
The quasi-isomorphism $\gamma'$ is just any lift in the following solid arrow
diagram
$$
\xymatrix
{
& & \Apl(W,B) \ar[d]^{\rho}_{\simeq} \\
s^{-r} \hat Q \ar@{-->}[urr]^{\gamma'}
\ar[rr]_{(\_\cdot \theta)s^{-r} \beta'}^{\simeq} & &
\Apl(T,\del T).
}
$$
The Lifting Lemma (Lemma \ref{lemma-tk-lifting}) implies the lift $\gamma'$
exists since by Lemma \ref{lemma-AplWB=Tdt} $\rho$ is an
acyclic fibration and $s^{-r}\hat Q$ is a
cofibrant $\hat R$-dgmodule.
\par
By the same argument there exists a lift $\tilde \phi^!$ in the solid arrow diagram
$$
\xymatrix
{
& \hat R \ar[d]^{\alpha'} \\
s^{-r}\hat Q \ar[r]_{\iota\gamma'}
\ar@{-->}[ur]^{\tilde\phi^!} & \Apl(W)
}
$$
It follows from Lemma
\ref{chirp} that
$H(\alpha)H(\tilde \phi^!)H(s^{-r}\beta)^{-1}(s^{-r}u_V)=u_W$. So
since by Proposition \ref{unish} $\phi^!$ represents the unique homotopy class of
maps with that property, we have that $\phi^!s^{-r}\beta\simeq
\alpha \tilde\phi^!$. Since $(\alpha,\alpha'):\hat R \rightarrow
R\oplus\Apl(W)$ is a surjection we can use
Lemma \ref{lemma-tk-nose} to replace
$\tilde\phi'$ by a map $\hat \phi^!:s^{-r}\hat Q\rightarrow \hat
R$ making Diagram (\ref{long}) commute on the nose.
\cqfd

\subsection{Replacing $\hat\phi\hat\phi^!$ by a $\hat Q$-dgmodule morphism}
\label{sect-give}
In this subsection we suppose fixed Diagrams (\ref{diag-phihat})
and (\ref{long}) with the maps $\gamma'$ and $\hat \phi^!$ constructed in Lemma
\ref{still}.
Here we show that the $\hat R$-dgmodule map $\hat\phi\hat\phi^!$ can be replaced by a
$\hat Q$-dgmodule map $\chi$ which is homotopic to
$\hat\phi\hat\phi^!$ in a controlled way.
\begin{lemma}
\label{lemma-chi}
\begin{itemize}
\item[(i)]
$\rho\gamma'\colon s^{-r}\hat Q \to\Apl(T,\bd T)$ is a morphism of
$\hat Q$-dgmodules.
\item[(ii)]
$\phi\phi^!=0$ and is thus a morphism of $\hat Q$-dgmodules.
\end{itemize}
\end{lemma}
\proof (i) Lemma \ref{still} says that $\gamma'$
makes Diagram (\ref{long}) commute so we have that
$$
\rho\gamma'=(\_ \cdot \theta) s^{-r}(\beta')
$$
Of course $(\_\cdot\theta)$ is an $\Apl(T)$-dgmodule
map and hence a $\hat Q$-dgmodule map. Also $\beta'$ and thus
$s^{-r}\beta'$ are $\hat Q$-dgmodule maps. Hence the composition
$$
(\_\cdot \theta)s^{-r}\beta'=\rho\gamma'
$$
is a $\hat Q$ dgmodule map.

(ii) We have assumed that $n\geq 2m+3$, hence $r=n-m\geq m+3$. Since $Q^{\geq m+2}=0$ the statement
follows.
\cqfd
In general $\hat\phi\hat\phi^!:s^{-r} \hat Q\rightarrow \hat Q$ is
not a $\hat Q$-dgmodule map but it can be adjusted as described in the following lemma.
\begin{lemma}\label{longer}
There exists a $\hat Q$-dgmodule
map $\chi:s^{-r}\hat Q \rightarrow \hat Q$ and a $\hat R$-dgmodule
homotopy $h:ss^{-r}\hat Q \rightarrow \hat Q$ from $\chi$ to $\hat
\phi\hat \phi^!$ such that $\beta h=\beta'h=0$. In particular
$(\beta,\beta')\chi=(\beta,\beta')\hat\phi\hat\phi^!$.
\end{lemma}
\proof
Since $n\geq 2m+3$, $r=n-m\geq m+1$, so
we know that $H^{<r}(s^{-r}\hat Q)=0$ and $H^{\geq
r}(T;\BQ)\subset H^{\geq m+1}(T;\BQ)=H^{\geq m+1}(V;\BQ)=0$. Thus Lemma
\ref{lemma-tk-nullhmtpy} implies that
$[s^{-r}\hat Q,\Apl(T)]_{\hat Q}=0$. Similarly $[s^{-r}\hat Q,Q]_{\hat Q}=0$.
Because
$s^{-r}\hat Q$ is a free $\hat Q$-dgmodule on one generator, it is
semi-free as a $\hat Q$-dgmodule. Thus using Lemma \ref{lemma-chi}
we see that the following diagram is homotopy commutative
in the category of $\hat Q$-dgmodules
$$
\xymatrix{
s^{-r}Q\ar[d]_{\phi\phi^!}&
\ar[l]^{\simeq}_{s^{-r}\beta} \ar[r]^{\rho \gamma'}_{\simeq} s^{-r}\hat Q
\ar[d]_0&
\Apl(T,\bd T)\ar[d]_{\iota'}\\
Q&
\ar[l]^{\simeq}_{\beta} \ar[r]^{\beta'}_{\simeq} \hat Q&
\Apl(T).
}
$$
Since $(\beta,\beta')$ is surjective, Lemma \ref{lemma-tk-nose}
asserts that we can replace the zero-map in the previous diagram
by a $\hat Q$-dgmodule morphism $\chi\colon s^{-r}\hat Q\to\hat Q$
making the diagram strictly commutative, in other words such that
$(\beta,\beta')\chi=(\beta,\beta')\hat\phi\hat\phi^!$.
\par
Next we construct the homotopy $h$.
Let $H'\colon Cyl (s^{-r}\hat Q) \rightarrow Q\oplus \Apl(T)$ be the map corresponding to the
constant
homotopy $h'=0\colon ss^{-r}\hat Q \rightarrow Q\oplus \Apl(T)$. We get the following
commutative solid arrow diagram
$$
\xymatrix
{
s^{-r}\hat Q \oplus s^{-r}\hat Q \ar[r]^-{\chi+\hat\phi\hat\phi^!}
\ar[d]_{i_0+i_1} & \hat Q \ar[d]^{(\beta,\beta')} \\
Cyl s^{-r}\hat Q \ar[r]_{H'} \ar@{-->}[ur]^H & Q\oplus \Apl(T).
}
$$
The map $i_0+i_1$ is a cellular $\hat R$-dgmodule cofibration such that
$H^{<r}(i_0+i_1)$ is an isomorphism. Also
since $n\geq 2m+3$, $r-1=n-m-1\geq m+1>\dim(V)$ and so
$H^{\geq r-1}(Q\oplus \Apl (T))=0=H^{\geq r-1}(\hat Q)$. Thus a lift $H$ exists. Our
desired $h\colon ss^{-r} \hat Q \rightarrow \hat Q$ is the homotopy corresponding to $H$.
\cqfd

\subsection{A dgmodule model of $j^*\oplus s\rho$.}
\label{sec-diagmod}
Recall from Lemma \ref{lemma-Aplconemodel} that the map of mapping cones
$$
j^*\oplus s\rho\colon \Apl(W) \oplus_{\iota} s\Apl(W,B)
\rightarrow \Apl(T) \oplus_{\iota'}s\Apl(T,\bd T)
$$
is a CDGA model of $k^*\colon \Apl(B)\rightarrow \Apl(\bd T)$.
Our aim in this subsection is to give another model (as dgmodules) of that map.
\par
We consider Diagrams (\ref{diag-phihat}) and (\ref{long}) with the maps
$\hat \phi^!$ and $\gamma'$ from Lemma \ref{still} as well as the
$\hat Q$-dgmodule map $\chi\colon s^{-r}Q \rightarrow \hat Q$ and the
$\hat R$-dgmodule homotopy $h:\chi\simeq \hat \phi\hat\phi$ from Lemma \ref{longer}.
Recall the notation from Lemma \ref{lemma-tk-conehmtpy}.

\begin{lemma}
\label{lemma-RQmodel-k}
The map
$$
\hat \phi\oplus_h \id\colon
\hat R\oplus_{\hat\phi^!}ss^{-r}\hat Q
\to
\hat Q\oplus_{\chi}ss^{-r}\hat Q
$$
is a morphism of $\hat R$-dgmodules making the following diagram commutative
in the category of $\hat R$-dgmodules
\begin{equation}\label{struck}
\xymatrix{
R\oplus_{\phi^!}ss^{-r} Q
\ar[d]_{\phi\oplus\id}&
\hat R\oplus_{\hat\phi^!}ss^{-r}\hat Q
\ar[l]_{\alpha\oplus ss^{-r}\beta}^{\simeq}
\ar[r]^-{\alpha'\oplus s\gamma'}_-{\simeq}
\ar[d]^{\hat \phi\oplus_h \id}&
\Apl(W)\oplus_{\iota}s\Apl(W,B)
\ar[d]^{j^*\oplus s\rho}
\\
Q\oplus_{\phi\phi^!}ss^{-r} Q&
\hat Q\oplus_{\chi}ss^{-r}\hat Q
\ar[l]_{\beta\oplus ss^{-r}\beta}^{\simeq}
\ar[r]^-{\beta'\oplus s\rho\gamma'}_-{\simeq}&
\Apl(T)\oplus_{\iota'}s\Apl(T,\bd T)
}
\end{equation}
Moreover the horizontal maps in this diagram are quasi-isomorphisms and the
bottom row consists of $\hat Q$-dgmodule maps.
\end{lemma}
\proof
Since Lemma \ref{longer} says that
$\beta h=\beta' h=0$, the
diagram in the statement of the lemma is commutative.

We see that
the horizontal maps of the diagram are
 quasi-isomorphisms by Lemma \ref{lemma-tk-conehmtpy}.
The fact that the bottom line of the diagram is
of $\hat Q$-dgmodules follows from the construction of $\beta$,
$\beta'$ and $\gamma'$.
\cqfd
\subsection{Extending dgmodule structure to CDGA structure}
\label{sec-extend}
For this subsection
we suppose fixed Diagrams (\ref{diag-phihat})
and (\ref{long}) with the maps $\gamma'$ and $\hat \phi^!$ constructed in Lemma
\ref{still} and
$\chi$ and $h$ from Lemma \ref{longer}.
\par
Here we use CDGA squares. At first sight
this may seem clumsy but it keeps track of the $\hat R$ and $\hat Q$ module
structures in a convenient way. Notation concerning squares is described in
Section \ref{sec-cat}.
\begin{lemma}
\label{esc}
Assume that $r$ is even.
There exists a unique CDGA structure on $\hat Q\oplus_{\chi} ss^{-r} \hat Q$
(respectively, $Q\oplus ss^{-r}  Q$) extending its $\hat Q$-dgmodule
(respectively, $Q$-dgmodule) structure. Moreover we can find CDGA isomorphisms $\hat e$
and $e$,  such that $e(z)=ss^{-r} 1$,
$\hat e(z)=(\hat q, ss^{-r} 1)$ for some $q\in \hat Q^{r-1}$, $e|_Q=id$,
$\hat e|_{\hat Q}=id$,  and which make
the following diagram commute
$$
\xymatrix
{
\hat Q\otimes \Lambda (z) \ar[r]^{\hat e} \ar[d]_{\beta\otimes id}
& \hat Q\oplus_{\chi} ss^{-r} \hat Q \ar[d]^{\beta\oplus s s^{-r}\beta} \\
 Q\otimes \Lambda (z) \ar[r]_e & Q\oplus ss^{-r}  Q
}
$$
where $|z|=r-1$ and $dz=0$.
\end{lemma}
\proof
Recall that $\phi\phi^!=0$ for dimension reasons.
The map $e$ is determined by the conditions that
$e(z)=ss^{-r} 1$, $e|_Q=id$ and the fact that it is a $Q$-module map. It is clearly
an isomorphism.
Since $\chi(ss^{-r} 1)$ is a cocycle in $\hat Q$ and since $H^{\geq r}(\hat Q)=0$, there exists
$\hat q\in \hat Q^{r-1}$ such that $d(\hat q)=-\chi(ss^{-r} 1)$.
For degree reasons $\beta(\hat q)=0$. We see also that $(\hat q,ss^{-r} 1)$ is a cocycle in
$\hat Q\oplus_{\chi} ss^{-r} \hat Q$.
Then $\hat e$ is
determined by the formula $\hat e(z)=(\hat q, ss^{-r} 1)$
and is an isomorphism. Clearly the diagram commutes. The CDGA structures
are unique since all products are determined by $(ss^{-r}(1))^2$ and
the $\hat Q$ and $Q$ module structures and $ss^{-r}(1)$ must square to $0$ since
it is in odd dimension.
\cqfd
\begin{lemma}\label{kuif}
For $r$ even, let $\mathbf{D}$ denote the commutative square
$$
\mathbf{D}=
\xymatrix
{
\hat R \ar[r]^{\hat \phi} \ar@{^{(}->}[d] & \hat Q \ar@{^{(}->}[d] \\
\hat R\oplus_{\hat\phi^!}ss^{-r}\hat Q
\ar[r]_{\hat \phi\oplus_h id} & \hat Q\oplus_{\chi}ss^{-r}\hat Q.
}
$$
and let $(\hat Q \otimes \Lambda z,dz=0)$ be the relative
Sullivan algebra from Lemma \ref{esc} together with the CDGA isomorphism
$\hat e\colon\hat Q \otimes \Lambda z \rightarrow \hat Q\oplus_{\chi}ss^{-r}\hat Q$.
\par
(i) There exists a CDGA square
$$
\mathbf{D}'=
\xymatrix
{
\hat R \ar[r]^{\hat \phi} \ar[d] & \hat Q \ar@{^{(}->}[d] \\
\hat A \ar[r]^-{\psi} & \hat Q\otimes \Lambda z
}
$$
and a weak equivalence between $\hat R$-dgmodule squares
$
\mathbf{\Theta}=\left(
\begin{array}{cc}
id & id\\
\theta_3 & (\hat e)^{-1}
\end{array}
\right)
\colon\mathbf{D} \rightarrow \mathbf{D}'.
$
\par
(ii)
Suppose
$$
\mathbf{C}=\xymatrix
{
C_1 \ar[r]\ar[d] & C_2 \ar[d]\\
C_3 \ar[r] & C_4
}
$$
is another CDGA square
and there is a map
$
\mathbf{\Theta}'=\left(
\begin{array}{cc}
\theta'_1 & \theta'_2\\
\theta'_3 & \theta'_4
\end{array}
\right)\colon \mathbf{D} \rightarrow \mathbf{C} $ such that $\theta'_1$ and $\theta'_2$
are CGDA maps, $\theta'_3$ is a $\hat R$-dgmodule map and
$\theta'_4$ is a $\hat Q$-dgmodule maps where $C_3$ (respectively, $C_4$) has been
given the $\hat R$-dgmodule (respectively, $\hat Q$-dgmodule) induced by
$\theta_1'$ (respectively, $\theta_2'$).
If $H^i(C_3)=H^i(C_4)=0$
for $i\geq 2r-3$ and
$C_3 \rightarrow C_4$ is a fibration, then there exists a CDGA map $\overline \theta_3$
such that
$\mathbf{\overline \Theta}= \left(
\begin{array}{cc}
\theta'_1 & \theta'_2\\
\overline\theta_3 & \theta'_4\hat e
\end{array}
\right)
\colon \mathbf{D}' \rightarrow \mathbf{C}
$
is a map of CDGA squares
making the following diagram commute.
$$
\xymatrix{
\mathbf{D} \ar[r]^{\mathbf{\Theta}}
\ar[d]_{\mathbf{\Theta}'} & \mathbf{D}'
\ar[ld]^{\mathbf{\overline \Theta}} \\
\mathbf{C}
}
$$
\end{lemma}
\proof
In view of Lemma \ref{esc}
to construct $\mathbf{D}'$ and $\mathbf{\Theta}$,
and prove part (i),
it is enough to construct
a relative Sullivan model $\hat R \rightarrow \hat A$ (thus also giving
$\hat A$ a $\hat R$-dgmodule structure),
a CDGA map $\psi:\hat A\rightarrow \hat Q\otimes \Lambda (z)$ and
an $\hat R$-dgmodule map
$\theta_3:\hat R \oplus_{\hat \phi^!} ss^{-r}\hat Q \rightarrow \hat A$ so that the
following diagram commutes
$$
\xymatrix
{
\hat R \oplus_{\hat \phi^!} ss^{-r} \hat Q \ar[r]^-{\theta_3}
\ar[dr]_{\hat e^{-1} (\hat \phi\oplus_h id)}  & \hat A \ar[d]^{\psi} \\
& \hat Q \otimes \Lambda z
}
$$
and $\theta_3$ is an equivalence.
Since $\hat \phi$ is a cellular cofibration of CDGAs,
$\hat Q=\hat R \otimes \Lambda U$ and so
$ss^{-r} \hat Q=ss^{-r} \hat R \otimes \Lambda U$.
There is an isomorphism of $\hat R$-dgmodules
$\rho\colon ss^{-r}\hat R \otimes \Lambda U \rightarrow
\hat R \otimes ss^{-r} \Lambda U$ given by
$\rho(ss^{-r}\alpha\otimes \beta)=(-1)^{(1-r)|\alpha|}\alpha\otimes ss^{-r}\beta$
where $\alpha\in \hat R$ and $\beta\in \Lambda U$.
Thus $ss^{-r}\hat Q\cong \hat R \otimes ss^{-r}\Lambda U$ as
$\hat R$-dgmodules and so
$\hat R\oplus_{\hat \phi^!} ss^{-r} \hat Q\cong \hat R\oplus \hat R \otimes ss^{-r}\Lambda U$.
Setting $U_0=ss^{-r} \Lambda U$,
$\hat R \oplus_{\hat \phi^!} ss^{-r} \hat Q$ is isomorphic as an
$\hat R$-dgmodule to $\hat R\otimes (\BQ\oplus U_0)$.
The inclusion $\BQ \oplus U_0\rightarrow \Lambda U_0$ induces an inclusion
$R\otimes (\BQ \oplus U_0)\rightarrow R\otimes \Lambda U_0$.
There is a unique differential on $R\otimes \Lambda U_0$
satisfying the Liebniz law such that the inclusion
$f_0\colon\hat R \oplus_{\hat \phi^!} ss^{-r} \hat Q \cong R\otimes
(\BQ \oplus U_0)\rightarrow \hat R \otimes \Lambda U_0$
is an $\hat R$-dgmodule map. Clearly then we get a commuting diagram
$$
\xymatrix
{
\hat R \ar[dr]\ar[r] & \hat R \oplus_{\hat \phi^!} ss^{-r} \hat Q\ar[d]^{f_0} \\
& \hat R \otimes \Lambda U_0. }
$$
Moreover any $\hat R$-dgmodule map from $\hat R \oplus_{\hat \phi^!} ss^{-r} \hat Q$
into a CDGA extends uniquely to a CDGA map out of
$\hat R \otimes \Lambda U_0$.
In particular we get a CDGA map
$g_0\colon\hat R\otimes \Lambda U_0 \rightarrow \hat Q\oplus_{\chi}ss^{-r}\hat Q$ and a
commutative diagram
$$
\xymatrix { \hat R \oplus_{\hat \phi^!} ss^{-r}\hat Q\ar[r]^{f_0}
\ar[dr]_{\hat \phi\oplus_h id}&
\hat R\otimes \Lambda(U_0) \ar[d]^{g_0} \\
& \hat Q \oplus_{\chi} ss^{-r}\hat Q.
}
$$
The elements of the cokernel of
$f_0$ are elements in $\Lambda U_0$ of product length at least
two. So they are in degree at least
$2r-2$. We add a minimal set of generators
$U_1$ of degree at least $2r-3$
to kill the cohomology of the cokernel of $f_0$. We then  get a
CDGA extension $h_1\colon \hat R \otimes \Lambda U_0 \rightarrow \hat R \otimes
\Lambda U_0\otimes \Lambda U_1$ and an
$\hat R$-dgmodule map
$f_1=h_1f_0\colon\hat R \oplus_{\hat \phi^!} ss^{-r} \hat Q \rightarrow
\hat R \otimes \Lambda U_0\otimes \Lambda U_1$.
We have assumed that $n\geq 2m+3$ and $r=n-m$, so $2r-2=n-m+r-2\geq m+r+1 >m+r-1$.
Recalling that $\hat Q$ models $\Apl(V)$ and
$\hat Q \oplus_{\chi} ss^{-r} \hat Q$ models $\Apl(\bd T)$ we know that
$\hat Q \oplus_{\chi} ss^{-r} \hat Q$ has trivial cohomology in degrees greater
than $n-1=m+r-1$ and so also in degrees
greater than or equal to $2r-2$ where the boundaries of $U_1$ lie.
Therefore the map $g_0$ can be extended over these new elements to a CDGA map
$g_1\colon \hat R\otimes \Lambda U_0 \otimes \Lambda U_1
\rightarrow \hat Q \oplus_{\chi} ss^{-r} \hat Q$.
\par
Of course $f_1$ has a new cokernel but since $H^0(\hat R)=\BQ$ (of course the homology
is $0$ in negative degrees)
and $U_1$ was chosen minimally the cohomology of this cokernel
is still in degrees greater than $2r-3$.
So again we can add generators to kill the cohomology of the
cokernel. We continue this process countably many times to get
our CDGA, $\hat A=\hat R\otimes \Lambda(\oplus_{i\geq 0} U_i)$ and our map
$\theta_3=f_{\infty}:\hat R \oplus_{\hat \phi^!} ss^{-r} \hat Q \rightarrow \hat A$.
Since each element of the cohomology of the cokernel of
$\theta_3$ was killed at the next stage, $H(\theta_3)$ is surjective.
Since $2r-2>n$, $H^{\geq 2r-2}(\hat R\oplus_{\hat \phi^!} ss^{-r} \hat Q)=0$
and we have not killed anything in
$H(\hat R\oplus_{\hat \phi^!} ss^{-r} \hat Q)$; so
$H(\theta_3)$ is injective. Thus $\theta_3$ is a quasi-isomorphism
of $\hat R$-dgmodules.
The original map $g_0$ extends to a CDGA map
$g_{\infty}\colon \hat A \rightarrow \hat Q \oplus_{\chi} ss^{-r} \hat Q$ since
$H^{\geq 2r-2}(\hat Q \oplus_{\chi} ss^{-r} \hat Q)=0$
and so at each stage
all obstructions are trivial for degree reasons.
So we can set $\psi=\hat e^{-1} g_{\infty}$ and we have completed the
construction of $\mathbf{D}'$ and $\mathbf{\Theta}$.
\par
We now proceed to construct the map $\mathbf{\overline\Theta}$ and prove part (ii).
Any $\hat Q$-dgmodule map out of $\hat Q\oplus_{\chi} ss^{-r}\hat Q$ is automatically
a CDGA map since $(ss^{-r} 1)\cdot (ss^{-r} 1)=0$ because $ss^{-r}1$ is of
odd degree.
Thus $\theta'_4$ and $\theta'_4\hat e$ are CDGA maps.
So we only have
to construct a CDGA map $\overline \theta_3:\hat A \rightarrow C_3$
extending $\theta'_3$ and making the following diagram commute
$$
\xymatrix
{
\hat A \ar[r]^-{\hat e \psi} \ar[d]_{\overline \theta_3} &
\hat Q \oplus_{\chi} ss^{-r} \hat Q \ar[d]^{\theta'_4} \\
C_3 \ar[r] & C_4.
}
$$
Clearly the $\hat R$-dgmodule map $\theta'_3$ extends uniquely to a CDGA map
$h_0\colon\hat R\otimes \Lambda U_0 \rightarrow C_3$
making the following diagram commute
$$
\xymatrix
{
\hat R \otimes \Lambda U_0 \ar[d]_{h_0} \ar[r]^{g_0} &
\hat Q \oplus_{\chi} ss^{-r} \hat Q \ar[d]^{\theta_4} \\
C_3 \ar[r] & C_4.
}
$$
Since $H^i(C_3)=0$ for $i>2r-2$, we can further extend
to a CDGA map $h'_1:\hat R \otimes \Lambda U_0\otimes \Lambda U_1
\rightarrow C_3$. We show that the following CDGA diagram
$$
\xymatrix
{
\hat R \otimes \Lambda U_0 \otimes \Lambda U_1
\ar[d]_{h'_1} \ar[r]^{g_1} & \hat Q \oplus_{\chi} ss^{-r} \hat Q \ar[d]^{\theta_4} \\
C_3 \ar[r] & C_4
}
$$
commutes up to CDGA homotopy. Indeed it already commutes on
$\hat R \otimes \Lambda U_0$ and hence
the two ways of going around differ by an element
of $[\Lambda U_1,C_4]$ and this group is $0$ since
$H^i(C_4)=0$ for $i\geq 2r-3$ and $U_1$ has no elements in degrees $\leq 2r-3$.
Because $C_3 \rightarrow C_4$ is a
surjection we can replace $h'_1$ by a map
$h_1:\hat R \otimes \Lambda U_0 \otimes \Lambda U_1 \rightarrow C_3$
making the last diagram commute exactly. Using this same method at each
stage and taking the direct limit we get our desired CDGA map
$\overline \theta_3:\hat A \rightarrow C_3$ making the
diagram commute.
\cqfd
The approach taken in the last lemma gives a hint of how to
approach the problem of extending an $\hat R$-dgmodule structure to
a CDGA structure when there are no dimension restrictions. At each stage one would
have to choose which representatives of the cokernel to kill.
Now we describe a CDGA model of Diagram (\ref{decodiagpo})
\begin{lemma}\label{hate}
Assume that $r$ is even.
With the semi-trivial CDGA on the mapping cones, the square
$$
\mathbf{E}=
\xymatrix
{
R \ar[r]^{\phi} \ar@{^{(}->}[d] & Q \ar@{^{(}->}[d] \\
R\oplus_{\phi^!}ss^{-r} Q
\ar[r]_{\phi\oplus id} & Q\oplus_{\phi\phi^!}ss^{-r} Q.
}
$$
is a commuting square of CDGAs.
\par
Also there exists a CDGA square
$$
\mathbf{E}'=
\xymatrix
{
R \ar[r]^{\phi} \ar@{^{(}->}[d] & Q \ar@{^{(}->}[d] \\
A \ar[r]_-{\kappa} & Q\otimes \Lambda (z).
}
$$
such that $\kappa$ is a fibration and there exists a CDGA
quasi-isomorphism $\xi\colon R\oplus_{\phi^!} ss^{-r} Q \rightarrow A$ such that
$
\mathbf{\Theta}_1=\left(
\begin{array}{cc}
\id & \id\\
\xi & e^{-1}
\end{array}
\right)
\colon \mathbf{E} \rightarrow \mathbf{E}'$ is a weak equivalence of CDGA squares.

\end{lemma}
\proof
That
$\mathbf{E}$ is a CDGA square follows using Proposition \ref{propdef-semitrivCDGA}. To
construct $\mathbf{E}'$ and $\mathbf{\Theta}_1$
we factor the CDGA map
$\phi\oplus id\colon R \oplus_{\phi^!} ss^{-r} Q \rightarrow Q \oplus_{\phi\phi^!} ss^{-r} Q$
as an acyclic cofibration
$\xi\colon R \oplus_{\phi^!} ss^{-r} Q \rightarrow A$ followed by a fibration
$\kappa\colon A \rightarrow Q \oplus_{\phi\phi^!} ss^{-r} Q$. The lemma follows
easily.
\cqfd
In the next lemma $\mathbf{D}$ and $\mathbf{D}'$ are the diagrams of
Lemma \ref{kuif}, $\mathbf{E}$ and $\mathbf{E}'$ are the
diagrams of the previous lemma
and
$\mathbf{F}$, $\mathbf{F'}$ and $\mathbf{\Theta}_6$ are the diagrams
and map from Lemma \ref{lemma-Aplconemodel}.
\begin{lemma}\label{igog}
Assume that $r$ is even.
Let $\rho$ be the map from \ref{sec-preap}, $\gamma'$ from
Section \ref{sec-comsh}, $\hat e$ from Lemma \ref{esc} and
$\mathbf{\Theta}_1$ from Lemma \ref{hate}.
Set
$$
\mu'=(i^* + 0)(\beta' \oplus s\rho\gamma')\hat e:
\hat Q\otimes \Lambda z\rightarrow \Apl(\bd T)
$$
which is a quasi-isomorphism.
There exist CDGA maps $\zeta\colon \hat A\stackrel{\simeq}{\rightarrow} A$ and
$\zeta'\colon \hat A\stackrel{\simeq}{\rightarrow} \Apl(B)$ such that
if we set
$
\mathbf{\Theta}_2=\left(
\begin{array}{cc}
\alpha & \beta\\
\zeta & \beta\otimes id
\end{array}
\right)
$
 and
$
\mathbf{\Theta}_3=\left(
\begin{array}{cc}
\alpha' & \beta'\\
\zeta' & \mu'
\end{array}
\right)
$
then we have a chain of quasi-isomorphisms of CDGA squares
$$
\mathbf{E}
\stackrel{\mathbf{\Theta}_1}{\rightarrow}
\mathbf{E}'
\stackrel{\mathbf{\Theta}_2}{\leftarrow}
\mathbf{D}'
\stackrel{\mathbf{\Theta}_3}{\rightarrow}
\mathbf{F}.
$$

\end{lemma}
\proof
Set
$$
\mathbf{\Theta}_4=\left(
\begin{array}{cc}
\alpha & \beta \\
\alpha\oplus ss^{-r}\beta & \beta \oplus_h ss^{-r}\beta
\end{array}
\right)\colon {\mathbf D} \rightarrow {\mathbf E}
$$
and
$$
\mathbf{\Theta}_5=\left(
\begin{array}{cc}
\alpha' & \beta'\\
\alpha'\oplus s\gamma' & \beta'\oplus s\rho\gamma'
\end{array}
\right)\colon {\mathbf D} \rightarrow {\mathbf F'}.
$$
By Lemmas \ref{lemma-Aplconemodel}, \ref{lemma-RQmodel-k}  and
\ref{hate} we have a
weak equivalences of $\hat R$-dgmodule squares
$\mathbf{\Theta}_4\mathbf{\Theta}_1\colon {\mathbf D} \rightarrow {\mathbf E'}$ and
$\mathbf{\Theta}_6\mathbf{\Theta}_5\colon{\mathbf D} \rightarrow {\mathbf F}$.
Also the maps between
the objects in the right hand column of each square are $\hat Q$-dgmodule maps.
We can  apply Lemma \ref{kuif} to these maps and get
the string of quasi-isomorphisms of squares as stated in the lemma.
\cqfd

\section{Model of the projectivization of a complex bundle}
\label{section-modelblowup}
In this section we suppose that $f:V\rightarrow W$ is a smooth embedding of closed manifolds
of codimension $2k$.
We will assume that $k>1$.
We also suppose that the normal bundle $\nu$ of the embedding has some fixed structure
of a complex vector bundle,
$$
\nu\colon\BC^k\to E\stackrel{\pi}\to V.
$$
Let $T$ be a compact tubular neighborhood of $V$ in $W$. By the Tubular Neighborhood
Theorem, we can
identify
$T$ with the disk bundle $D\nu$ and $\bd T$ with its sphere bundle
$S\nu$ in such a way that the zero section of $D\nu$ corresponds to the inclusion
$$
\sigma: V \rightarrow T.
$$
We fix such an identification.
We also suppose we have been given some CDGA model $Q$ of $\Apl(T)$ and a common
model $\hat Q$
$$
\xymatrix
{
Q & \hat Q \ar[r]^-{\beta'} \ar[l]_{\beta} & \Apl(T)
}
$$
of $Q$ and $\Apl(V)$
such that $(\beta,\beta')$ is surjective and $\beta$
and $\beta'$ are quasi-isormorphisms as in Section \ref{sec-supps}.
Clearly $\sigma^*\beta'\colon \hat Q\rightarrow \Apl V$ is also
a quasi-isomorphism.
By Lemma \ref{lemma-tk-commonmodel} such a common model can
be constructed from any CDGA model $Q$ of $\Apl(T)$.
\par
The aim of this section is to describe the projective bundle
$P\nu$ associated to $\nu$ and give a CDGA model for this projective bundle.
In Section \ref{Chern}
we review the definition of the projective bundle and of Chern classes and prove
the triviality of a certain line bundle
(Lemma \ref{lemma-qa0}). Next in Sectoin \ref{orient} we consider the
pullback over a point of the sphere bundle and of its projectivization.
Using orientation information we
show that the models of these pullbacks and the model of
the sphere bundle $P\nu$ can be
chosen in a compatible way (Lemmas \ref{B1} and \ref{B3}). Then in
Section \ref{mod of proj} we construct a model of the projectivization of the
sphere bundle and use the results from
Section \ref{orient} to show this model is compatible with the model of
the boundary of the normal bundle
(Lemma \ref{(I)} and Proposition \ref{prop-modelprojxi}).

\subsection{The projective bundle and the Chern classes}\label{Chern}
Next we recall the definition of the Chern classes of a complex bundle using the
associated projective bundle as described in
\cite[IV.20]{BottTu}. Consider the universal complex line bundle $\gamma^1$
over $\CP(\infty)$.
$$
\BC \rightarrow E\gamma^1 \rightarrow \CP(\infty)
$$
where $\CP(\infty)=\{ l \colon l {\rm \ is \ a \ } \BC {\rm \ line \ in \ } \BC^{\infty}\}$
and $E\gamma^1=\{ (l,v) \colon l \in \CP(\infty), v\in l \}$.
This complex line bundle can also be viewed as an oriented real vector bundle
of rank 2. Therefore
we have an associated Euler class $e(\gamma^1)$ and we set
$$
a_{\infty}\colon = -e(\gamma^1)\in H^2(\CP(\infty), \BZ)
$$
This is our preferred generator of the cohomology algebra
$H^*(\CP(\infty),\BZ)$. Note that if $j\colon S^2=\CP(1)\rightarrow \CP(\infty)$ is the
obvious inclusion then $j^*(a_{\infty})\in H^2(S^2;\BZ)$ is the orientation class corresponding
to the orientation coming from its complex structure.

To any complex vector bundle $\nu$ we can associate its {\em
projective bundle} which is defined as follows
(see \cite[page 269]{BottTu}).
Set $E_0=E\setminus\{\mbox{zero section}\}$ and consider the bundle
$$
\nu_0\colon\BC^k\setminus\{0\}\to E_0\stackrel{\pi_0}\to V
$$
Then $\BC^*=\BC\setminus\{0\}$ acts on each fibre in $E_0$
by complex multiplication and we can define the orbit space $P\nu=E_0/\BC^*$.
In other words,
$$
P\nu=\{(v,\ell):v\in V,\ell\mbox{\,is a complex line in the fibre \,}E_v:=\pi^{-1}(v)\}
$$
and we have the projective bundle
$$
\CP(k-1)\stackrel{\rm inc}\to P\nu\stackrel{\pi'}\to V
$$
Denote by $q\colon E_0\to P\nu$ the quotient map and consider
the commutative diagram

\begin{equation}\label{quotie}
\xymatrix
{
E_0\ar[rr]^q\ar[rd]_{\pi_0}&&P\nu\ar[ld]^{\pi'}\\
& V.
}
\end{equation}
Since the inclusion $\bd T\cong S\nu\hookrightarrow E_0$ is a homotopy equivalence,
we also use $q$ to denote the composition
$
\xymatrix
{
\bd T\cong S\nu\ar[r] & E_0 \ar[r]^q & P\nu.
}
$

We now come to the definition
of the Chern classes $c_i(\nu)\in H^{2i}(V;\BZ)$ in terms of the
projective bundle.
The pullback of $\nu$ along $\pi'$ is a
$\BC^k-$bundle $\pi'^*(\nu)$ over $P\nu$ containing a
tautological line bundle defined by
\begin{equation}\label{equ-tautol}
\lambda=\{(v,\ell,x):v\in V,\ell\mbox{\,is a line in \,}E_v,x\in\ell\}
\end{equation}
The complex line bundle $\lambda$ is classified by some map
$\alpha:P\nu \rightarrow \CP(\infty)\simeq BU(1)$, that is
$\lambda\cong \alpha^*(\gamma^1)$. We define the
{\em canonical class} of $\nu$ as the cohomology class
\begin{equation}\label{equ-a}
a=\alpha^*(a_{\infty})\in H^2(P\nu,\BZ).
\end{equation}
Since
the restriction of that class to each fibre $\CP(k-1)$ of $\pi'$
is a generator of the cohomology (because the pullback of $\lambda$
to that fibre is the universal line bundle over $\CP(k-1)$), the
Leray-Hirsch Theorem gives an isomorphism of $H^*(V,\BZ)-$algebras
\begin{equation}\label{equ-chernclasses}
H^*(P\nu,\BZ)\cong H^*(V,\BZ)[a]/ \left(\sum_{i=0}^k c_i(\nu)a^{k-i}\right)
\end{equation}
where $c_0(\nu)=1$ and by definition the $c_i(\nu)\in H^{2i}(V,\BZ)$ for
$i>0$ are the Chern classes.
Notice a straightforward calculation shows that
\begin{equation}\label{equ-sup}
a=-c_1(\lambda).
\end{equation}
\begin{lemma}
\label{lemma-qa0}
$q^*(a)=0$.
\end{lemma}
\proof
Consider
the following diagram of vector bundles
$$\xymatrix{q^*\lambda\ar[d]\ar[r]
&\lambda\ar[rd]\ar@{^{(}->}[r]
&\pi'^*(\nu)\ar[r]\ar[d]\ar@{}[rd]|{\mbox{pullback}}
&\nu\ar[d]\\
E_0\ar[rr]^q&&P\nu\ar[r]^{\pi'}\ar[rd]_{\alpha} &V\\
&&&\CP(\infty)=BU(1)}$$
The cohomology class
$q^*(a)$ is the classifying class of the line bundle
$q^*(\lambda)$. We will prove that this bundle is trivial, so that
$q^*(a)=0$.

By (\ref{equ-tautol}) we have
\begin{eqnarray*}
q^*\lambda&=&\{(e,v,\ell,x):e\in E_0,v\in V,v=\pi_0(e),
\ell\mbox{\,is a line in \,}E_v,x\in\ell\}\\
&\cong&\{(e,x):e\in E_0,x\in \BC.e\subset E_{\pi_0(e)}\}
\end{eqnarray*}
Therefore we have the trivialization
$$E_0\times\BC\isom q^*\lambda\,\,,\,\,(e,z)\mapsto(e,e.z)$$
\cqfd
\subsection{Orientations}\label{orient}
As in Section \ref{sec-Thom} we suppose $*\in V$ is some fixed point.
Consider the map
$$
\eta:S^{2k-1} \rightarrow \CP(k-1)
$$
obtained as the restriction of the map
$$
q:S\nu \rightarrow P\nu
$$
to the fibres over our point $*\in V$.
The aim of this subsection is to give an explicit CDGA model of the map
$\eta$.

First we define some more cohomology classes. Recall the class
$a\in H^2(P\nu)$ defined in (\ref{equ-a}).
We also denote by $a\in H^2(\CP(k-1))$ its restriction to the
fibre. Notice that the restriction of $\lambda$ to $\CP(k-1)$ is the tautological
line bundle $\gamma_{k-1}$ and hence by Equation (\ref{equ-sup}) and the naturality of the Chern classes
\begin{equation}\label{equ-acpk-1}
a=-c_1(\gamma_{k-1}).
\end{equation}
Equivalently we can take Equation (\ref{equ-acpk-1}) to be the definition
of $a\in H^2(\CP(k-1))$.
The fibre $D^{2k}$ of the disk bundle $D\nu$
over our point $*\in V$
has a canonical orientation given by its complex structure
$D^{2k}\subset \BC^k$. This determines the cohomology class
\begin{equation}\label{equ-uD}
u_{D^{2k}}\in H^{2k}(D^{2k},S^{2k-1};\BZ)
\end{equation}
of Equation (\ref{equ-ornormal})
which through the connecting homomorphism
$$
\xymatrix
{
\delta\colon H^{2k-1}(S^{2k-1}; \BZ) \ar[r]^{\cong} &
H^{2k}(D^{2k},S^{2k-1}; \BZ)
}
$$
determines an orientation class
\begin{equation}\label{equ-uS}
u_S\in H^{2k-1}(S^{2k-1},\BZ).
\end{equation}
\par
Next we extend $\eta$ to some map
$\tilde \eta:D^{2k}\rightarrow \CP(k)$ and describe
$\CP(k)$ as a suitable pushout.
For $z=(z_1, \dots , z_k)\in \BC^k$ we set
$||z||_2=\left(\Sigma_{i=1}^k |z_i|^2\right)^{1/2}$.
Consider the disc and the sphere
$$
D^{2k}=\{ z\in \BC^k\colon ||z||_2\leq 1\}
$$
$$
S^{2k-1}=\{ z\in \BC^k\colon ||z||_2=1\}.
$$
Define a map
$$
\tilde\eta:D^{2k} \rightarrow \CP(k),
z \mapsto [z_1: \dots :z_k: 1-||z||_2]
$$
Then $\eta=\tilde\eta |_{S^{2k-1}}$ where $\CP(k-1)$ is considered as the hyperplane
with last coordinate $0$.
\par
We need to replace the map $\eta$ by an equivalent map that is a cofibration.
This is done using the standard mapping cylinder construction. We will describe it explicitly
so that it is easy to see how it is compatible with $\tilde\eta$.
Define the contraction
$$
\rho:D^{2k}\rightarrow D^{2k}, z\mapsto z/2.
$$
Set
$$
X=\{ \tilde\eta(z)\colon z\in D^{2k}, 1/2\leq||z||_2\leq 1\}\subset \CP(k)
$$
and denote by $l\colon S^{2k-1} \rightarrow D^{2k}$,
$\tilde l\colon X \rightarrow \CP(k)$ and
$s\colon \CP(k-1) \rightarrow X$ the obvious inclusions. Then the composite
$$
\tilde h=\tilde\eta\rho:D^{2k} \rightarrow \CP(k)
$$
is a homeomorphism onto its image and a cofibration.
Also we can think of $X$ as a tubular neighbourhood of $\CP(k-1)$ in
$\CP(k)$ although we will not use this fact.
\par
We summarize a few facts about these maps in the following lemma.
Its proof is straightforward.
\begin{lemma}\label{sapsap}

(i) The restriction of $\tilde h$ to $S^{2k-1}$ induces a cofibration
$$
h:S^{2k-1} \rightarrow X
$$
(ii) There is a pushout
$$
\xymatrix
{
S^{2k-1} \ar[r]^h \ar[d]_l & X \ar[d]^{\tilde l} \\
D^{2k} \ar[r]_{\tilde h} & \CP(k)
}
$$
(iii) The inclusion map $s:\CP(k-1)\rightarrow X$ is a homotopy
equivalence.
\par
\noindent
(iv) The following diagram commutes up to homotopy
$$
\xymatrix
{
S^{2k-1} \ar[r]^h \ar[dr]_\eta & X \\
& \CP(k-1). \ar[u]_{s}
}
$$
\end{lemma}
\par
By abuse of notation we also use $a\in H^2(X)$ to denote
the preimage through the isomorphism $s^*$ of the element $a\in H^2(\CP(k-1))$
defined in (\ref{equ-acpk-1}).
\begin{lemma}\label{sap}
Consider the Sullivan algebra $(\Lambda(x,z);dx=0,dz=x^k)$ with
$|x|=2$ and $|z|=2k-1$.
There exists a homotopy commutative diagram of CDGA
$$
\xymatrix
{
\Lambda(x,z) \ar[d]_{g_0} \ar[r]^{\rm proj} &
\Lambda(z) \ar[d]^{f_0} \\
\Apl (X)  \ar[r]_-{\Apl(h)}  & \Apl(S^{2k-1}).
}
$$
such that $[g_0(x)]=a$ (defined in (\ref{equ-acpk-1})) and $[f_0(z)]=-u_S$
(defined in (\ref{equ-uS})).
\end{lemma}
\proof
Since $k>1$, the restriction map
$$
\iota: H^2(X,S^{2k-1})\rightarrow H^2(X)
$$
is an isomorphism and so the class $a$ lifts to a unique class
$a_0\in H^2(X,S^{2k-1})$.
\par
Denote by $\tilde a\in H^2(\CP(k))$ the opposite of the first
Chern class of the tautological line bundle $\gamma_k$ over
$\CP(k)$. Since $\gamma_{k-1}$ is the restriction of
$\gamma_k$ over $\CP(k-1)$,
Equation (\ref{equ-acpk-1}) and the naturality of the Chern classes implies that
$\tilde l^*(\tilde a)=a$.
In view of the pushout of Lemma \ref{sapsap}(ii), $\tilde a$
lifts to some class
$\tilde a_0\in H^2(\CP(k),D^{2k})$
such that $(\tilde l,l)^*(\tilde a_0)=a_0$.
\par
Let $\alpha_0\in \Apl^2(X,S^{2k-1})\cap \ker d$ be a representative
of $a_0$ and $\alpha$ be its image in $\Apl^2(X)\cap \ker d$.
Since $a^k=0$ in $H^2(X)$, there exists $\zeta\in \Apl^{2k-1}(X)$
such that $d\zeta=\alpha^k$. Define
$$
g_0\colon (\Lambda(x,z);dz=x^k) \rightarrow \Apl(X)
$$
by $g_0(x)=\alpha$ and $g_0(z)=\zeta$.
\par
Since $h^*(\alpha)=0$, $h^*(\zeta)$ is a cocycle in $\Apl(S^{2k-1})$ and so
we can define
$$
f_0\colon \Lambda(z) \rightarrow \Apl(S^{2k-1}), \ z\mapsto h^*(\zeta).
$$
This definition makes the diagram of the lemma commutative.
\par
We proceed to prove that $[h^*(\zeta)]=-u_S$. This will imply that $[f_0(z)]=-u_S$
and thus complete the proof of the lemma.
\par
Let $\delta$ and $\delta'$ denote the cohomology connecting homomorphisms
of the pairs $(D^{2k},S^{2k-1})$ and $(X,S^{2k-1})$ respectively.
Consider the following diagram in cohomology
$$
\xymatrix
{
H^{2k-1}(S^{2k-1};\BZ)\ar@/^1.5pc/[rr]^{\delta'} \ar@/_3pc/[dd]_{\delta}
& H^{2k}(X;\BZ) &
H^{2k}(X,S^{2k-1};\BZ) \ar[l]_{\iota} \\
H^{2k}(D^{2k};\BZ) & H^{2k}(\CP(k);\BZ) \ar[l]_{\tilde h^*} \ar[u]_{\tilde l^*}
 & H^{2k}(\CP(k),D^{2k};\BZ) \ar[l]_{\tilde \iota} \ar[u]_{(\tilde l,l)^*}\\
H^{2k}(D^{2k},S^{2k-1};\BZ)\ar[u] &
H^{2k}(\CP(k),X;\BZ)\ar[u]_{\tilde\iota'}\ar[l]_-{(\tilde h,h)^*}.
}
$$
A diagram chase at the chain level gives us the formula
\begin{equation}\label{Lin1}
\tilde \iota'((\tilde h,h)^*)^{-1}\delta=-\tilde \iota((\tilde l,l)^*)^{-1}\delta',
\end{equation}
the minus sign corresponding to the fact that in the Mayer-Vietoris sequence
there is a minus sign on one of the maps.
\par
By the definition of the pushout of Lemma \ref{sapsap}(ii) and
the construction of $\zeta$ we have
$$
\delta'([h^*\zeta])=[\alpha^k]=a_0^k\in H^{2k}(X,S^{2k-1})
$$
Since $(\tilde l,l)^*$ and $\tilde \iota$ are multiplicative we get
\begin{equation}\label{Lin2}
\tilde\iota((\tilde l,l)^*)^{-1}\delta'[h^*\zeta]=\tilde a^k\in H^{2k}(\CP(k))
\end{equation}
On the other hand since $\tilde a=-c_1(\lambda)$,
\cite[page 170]{MS} says that the orientation of
$\CP(k)$ induced by its complex structure corresponds to the
class
$$
\tilde a^k\in H^{2k}(\CP(k)).
$$
Since $\tilde h$ preserves the given orientations, we have
$\tilde \iota'((\tilde h,h)^*)^{-1}u_D=\tilde a^k$, and hence
\begin{equation}\label{Lin3}
\tilde\iota'((\tilde h,h)^*)^{-1}\delta(u_S)=\tilde a^k.
\end{equation}
Equations (\ref{Lin1}), (\ref{Lin2}) and (\ref{Lin3}) imply that
$$
\tilde\iota'((\tilde h,h)^*)^{-1}\delta(u_S)=-\tilde\iota'((\tilde h,h)^*)^{-1}\delta [h^*\zeta]
$$
Thus since $\tilde\iota'((\tilde h,h)^*)^{-1}\delta\colon H^{2k-1}\rightarrow H^{2k}\CP(n)$
is an isomorphism
we deduce that
$$
[h^*\zeta]=-u_S
$$
This completes the proof of the lemma.
\cqfd
\begin{lemma}\label{B1}
Let $|x|=2$ and $|z|=2k-1$.
Suppose we are given a CDGA diagram
$$
\xymatrix
{
(\Lambda (x,z);dz=x^k) \ar[d]_g \ar[r]^-{\rm proj} &
\Lambda(z) \ar[d]^f \\
\Apl (\CP(k-1)) \ar[r]_-{\Apl(\eta)} & \Apl(S^{2k-1}).
}
$$
such that $[g(x)]=a$ and $[f(z)]=-u_S$. Then the diagram is homotopy commutative
and $f$ and $g$ are quasi-isomorphisms.
\end{lemma}
\proof
Observe that since $H^{2k-1}(\CP(k-1))=0$ the homotopy class of the map
$g$ is determined by the fact that $[g(x)]=a$. This equation also implies that
$g$ is a quasi-isomorphism. Similarly the equation $[f(z)]=-u_S$ determines
the homotopy class of the map and implies that it is a quasi-isomorphism.
The lemma then follows easily from Lemma \ref{sapsap} (iii) and (iv) and Lemma \ref{sap}.
\cqfd
\begin{lemma}\label{hunger}
Consider the map
$\mu':\hat Q\otimes \Lambda(z)\rightarrow \Apl(\bd T)$ from
Lemma \ref{igog}.
\par
(i) Consider the connecting homomorphism
$$
\delta:H^{2k-1}(\bd T)\rightarrow H^{2k}(T,\bd T)
$$
and the Thom class $\overline \theta$ from Section \ref{sec-Thom}.
Then
$$
\delta[\mu'(z)]=-\overline \theta.
$$
\par
(ii)
Let ${\rm inc}:S^{2k-1}\rightarrow \bd T$ be the inclusion. Then
$$
{\rm inc}^*[\mu'(z)]=-u_S.
$$
\end{lemma}
\proof
By the definition of $\mu'$ in Lemma \ref{igog}, the following diagram commutes
$$
\xymatrix
{
\hat Q\otimes \Lambda(z)\ar[ddrr]_{\mu'} \ar[rr]^-{\hat e} &&
\hat Q \oplus ss^{-2k} \hat Q
\ar[d]^{\beta'\oplus s\rho\gamma'} \\
&& \Apl(T)\oplus_{\iota'} s\Apl(T,\bd T) \ar[d]^{i^* + 0}\\
&& \Apl(\bd T).
}
$$
By the construction of $\hat e$ in Lemma \ref{esc} there exists $\hat q\in \hat Q^{2k-1}$
such that $\hat e(z)=(\hat q,ss^{-2k}1)$. By Lemma \ref{still}
$\rho\gamma'(ss^{2k}(1))=\theta$ and we all see that
$$
(\beta'\oplus s\rho\gamma')(\hat q,ss^{-2k}1)=(\beta'(\hat q),s\theta),
$$
so
\begin{equation}\label{star}
\mu'(z)=i^*(\beta'(\hat q)).
\end{equation}
Since $z$ is a cocycle, $(\beta'(\hat q),s\theta)\in \Apl(T)\oplus_{\iota'}s\Apl(T,\bd T)$
is also a cocycle which implies that in $\Apl(T)$
$$
d(\beta'(\hat q))=-\iota'(\theta).
$$
The definition of the connecting homomorphism then implies that
$$
\delta [i^*\beta'(\hat q)]=-[\theta]=-\overline\theta
$$
Combining this formula with Equation (\ref{star}) proves (i).
\par
(ii) This is a consequence of (i), of the fact that
$u_{D^{2k}}={\rm inc}^*(\overline \theta)$ (see Equation (\ref{equ-incT}) and
of the naturality of the connecting homomorphism in the following diagram
$$
\xymatrix{
H^{2k-1}(\bd T) \ar[r]^{\delta} \ar[d]_{{\rm inc}^*} & H^{2k}(T,\bd T)
\ar[d]^{{\rm inc}^*} \\
H^{2k-1}(S^{r-1})\ar[r]^{\cong}_{\delta} & H^{2k}(D^{2k},S^{2k-1}).
}
$$
\cqfd
\begin{lemma}\label{B3}
The inclusion of the point $*\in V$
determines an augmentation on $\hat Q$ using the composition
$$
\hat Q \stackrel{\sigma^*\beta'}{\rightarrow} \Apl(V) \rightarrow \Apl(*)=\BQ
$$
In turn the augmentation determines the projection map
${\rm proj}\colon \hat Q\otimes \Lambda (z) \rightarrow \Lambda(z)$.
Suppose $f\colon (\Lambda(z);dz=0)\rightarrow \Apl(S^{2k-1})$ is any CDGA map such
that $[f(z)]=-u_S$.
Then the following diagram commutes up to CDGA homotopy
$$
\xymatrix
{
\hat Q \otimes \Lambda(z) \ar[r]^{\rm proj} \ar[d]_{\mu'} &
\Lambda(z) \ar[d]^f  \\
\Apl (\bd T) \ar[r]_{{\rm inc}^*} & \Apl (S^{2k-1})
}
$$
\end{lemma}
\proof
In this proof we denote by $*$ an augmentation map followed by
a unit map. For fixed domain and range this composition
is unique up to homotopy since all of our CDGAs are homologically connected.
Since ${\rm proj}|_{\hat Q}=*$, $f {\rm proj}|_{\hat Q}=*$.
Also $\mu'|_{\hat Q}$ factors up to homotopy through $\Apl(V)$ so
${\rm inc}^*\mu'|_{\hat Q}\simeq *$. Because $dz=0$ in
$\hat Q \otimes \Lambda (z)$ (see Lemma \ref{esc})
the diagram commutes if and only if
${\rm inc}^*[\mu'(z)]=[f {\rm proj}(z)]$ in $H\Apl(S^{2n-1})$.
Lemma \ref{hunger} says that
${\rm inc}^*[\mu'(z)]=-u_S$ and by assumption
$[f(z)]=-u_S$. Thus the diagram homotopy commutes as required.
\cqfd
\subsection{A model of the projective and the sphere bundles}
\label{mod of proj}
Recall from the start of Section \ref{section-modelblowup} that
$V\hookrightarrow W$ has codimension $2k$ and so this is also the real rank of the normal
bundle $\nu$, which also has a complex structure of rank $k$. Recall also
the maps $\sigma$, $\beta$ and $\beta'$
from the start of Section \ref{section-modelblowup}.

\par

We describe the Sullivan models that we will prove are models for the
projective bundle $P\nu$.
Set $\gamma_0=1$ and for $1\leq i\leq k$ let $\gamma_i\in Q^{2i}\cap \ker d$ be
representatives of the images of the Chern classes
$$
c_i(\nu)\in H^{2i}(V,\BZ)\rightarrow
H^{2i}(V,\BQ) \stackrel{(\sigma\beta')^*}{\cong} H^{2i}(\hat Q)
\stackrel{\beta^*}{\cong} H^{2i}(Q).
$$
Since $\beta$ is a surjective quasi-isomorphism there exists
$\hat \gamma_i\in \hat Q^{2i}\cap \ker d$ such that
$\beta(\hat \gamma_i)=\gamma_i$.
Also we can
take $\hat \gamma_0=1$. Let $|x|=2$ and $|z|=2k-1$
and define relative Sullivan models
\begin{equation}\label{equ-Qxz}
(Q\otimes \Lambda (x,z);Dx=0, \ Dz=\Sigma_{i=0}^k \gamma_i x^{k-i})
\end{equation}
and similarly
\begin{equation}\label{equ-hatQxz}
(\hat Q\otimes \Lambda (x,z);\hat Dx=0, \ \hat Dz=\Sigma_{i=0}^k \hat \gamma_i x^{k-i})
\end{equation}
These models are motivated by Equation (\ref{equ-chernclasses}) for the
Chern classes.
For the next lemma
recall that $\sigma:V\rightarrow T$ corresponds to the inclusion of the zero section.
In addition recall the classes $a\in H^2(P\nu)$ defined in (\ref{equ-a}) and its restriction
to $\CP(k-1)$ also denoted by $a\in H^2(\CP(k-1))$.
\begin{lemma}\label{B2}
Consider the projective bundle associated to $\nu$
$$
\xymatrix
{
\CP(k-1) \ar[r]^-{\rm inc} & P\nu \ar[r]^{\pi'} & V.
}
$$
Let $\hat D$ be the differential given in Equation (\ref{equ-hatQxz}).
Suppose $g\colon(\Lambda(x,z);dz=x^k)\rightarrow \Apl(\CP(k-1))$ is
any map such that $[g(x)]=a$.
Then there exists a quasi-isomorphism
$\theta'\colon(\hat Q\otimes \Lambda(x,z);\hat D) \rightarrow \Apl(P\nu)$
making the following diagram commute up to homotopy
$$
\xymatrix
{
\hat Q \ar[d]_{\sigma^*\beta'} \ar[r] & (\hat Q\otimes
\Lambda(x,z),\hat D)
\ar[d]_{\theta'}^{\simeq} \ar[r]^-{\rm proj} & (\Lambda(x,z),dz=x^k) \ar[d]^g \\
\Apl (V) \ar[r]^{\pi'^*} & \Apl(P\nu) \ar[r]^-{{\rm inc}^*} & \Apl(\CP(k-1)).
}
$$
\end{lemma}
\proof
As observed at the start of the section $\sigma^*\beta'$ is a quasi-isomorphism.
Define $\theta'|_{\hat Q\otimes \Lambda (x)}$ so that
$\theta'|_{\hat Q}=(\sigma\pi')^*\beta'$ and $\theta'(x)$
is any representative of the image of $a$ under the map
$$
H^*(P\nu, \BZ) \rightarrow H^*(P\nu, \BQ)\cong H(\Apl(P\nu))
$$
Then we know from the definition of the Chern classes (see Section \ref{Chern}) that
\break
$\Sigma_{i=0}^k c_i(\nu) a^{k-i}=0$ in $H^*(P\nu, \BZ)$ and so
$[\Sigma_{i=0}^k {\pi'}^*\sigma^*\beta'(\hat \gamma_i) (\theta'(x))^{k-i}]=0$ in
$H(\Apl(P\nu))$. Thus an extension over $z$ of
$\theta'|_{\hat Q\otimes \Lambda(x)}$ exists. Let $\theta'$ be any
such extension. It is clear using Equation (\ref{equ-chernclasses}) that $\theta'$
is a quasi-isomorphism and makes the left hand square
of our diagram commute. Since $({\rm inc}^*\theta')|_{\hat Q}=0$ the right
hand square commutes when restricted to $\hat Q$. Since
$g(x)$ represents ${\rm inc}^*(a)$ we see that the right
hand square commutes up to homotopy when restricted to
$\Lambda(x)$. Thus it commutes up to homotopy when restricted to
$\hat Q\otimes \Lambda(x)$.
We know that $[\Lambda(z), \Apl\CP(k-1)]=H^{2k-1}(\CP(k-1))=0$.
This implies that, up to homotopy, there is a unique extension of
the map $\hat Q\otimes \Lambda(x) \rightarrow \Apl(\CP(k-1))$
over $\hat Q\otimes \Lambda(x,z)$. Thus
the right hand square in the diagram commutes up
to homotopy and the lemma has been proven.
\cqfd
In this section $\theta'$ yielded a model of
$\Apl(P\nu)$. Lemma \ref{(I)} will show that the two are compatible.
This lemma
controls the automorphism of the model of $\Apl(\bd T)$.
If $Q$ is our model of $\Apl(T)$ then the model of $\Apl(\bd
T)$ is $Q\otimes \Lambda z$. An automorphism $\phi\in [\Apl(\bd
T),\Apl(\bd T)]$ can be considered as an element of $[Q\otimes
\Lambda z, Q\otimes \Lambda z]\cong [Q, Q\otimes \Lambda z]\times
[\Lambda z, Q\otimes \Lambda z]$. We control the first factor by
working in the category of $Q$-dgmodules. Another way to think of
this is that we work with objects together with maps from $Q$. To
handle the second factor we can observe that for dimension reasons
$[\Lambda z, Q \otimes \Lambda z]\cong[\Lambda z, \Lambda z]$ and in
turn we have the general isomorphism $[\Lambda z, \Lambda z]\cong
\Hom(H^{2k-1}\bd T,H^{2k-1}\bd T)$. So we only have to control a
single element in homology and this is the reason we have been keeping track of
orientation classes. In the last section $\mu'$ gave us a model of
$\Apl(\bd T)$.
\par
Recall $\beta$ from the start of the section,
$\mu'$ from Lemma \ref{igog}, $\theta'$ from Lemma \ref{B2},
$q\colon\bd T\rightarrow P\nu$ defined in Section \ref{Chern} and the
CDGAs $Q\otimes \Lambda(x,z)$ and $\hat Q \otimes \Lambda(x,z)$ defined in
equations (\ref{equ-Qxz}) and (\ref{equ-hatQxz}).
\begin{lemma}\label{(I)}
Assume that $2k\geq \dim(V)+2$.

The following diagram commutes up to homotopy
\begin{equation}
\label{hard}
\xymatrix
{
Q \otimes \Lambda(x,z) \ar[r]^{\rm proj}
&  Q \otimes \Lambda(z) \\
\hat Q \otimes \Lambda(x,z) \ar[u]^{\beta\otimes id} \ar[r]^{\rm proj}
\ar[d]_{\theta'} & \hat Q \otimes \Lambda(z) \ar[u]_{\beta\otimes id} \ar[d]^{\mu'} \\
\Apl(P\nu) \ar[r]_{q^*} & \Apl(\bd T)
}
\end{equation}
\end{lemma}
\proof
The top square clearly commutes.
Since we will be using different inclusions,
to avoid confusion for the rest of the proof
set $I={\rm inc}^*:\Apl(\bd T) \rightarrow \Apl(S^{2k-1})$.
To avoid too much clutter in the equations we will also sometimes write $I$ instead of
$I_*$ for induced maps.
Recall the classes $a\in H^2(P\nu)$ defined in (\ref{equ-a}) and also its restriction
to $\CP(k-1)$ also denoted by $a\in H^2(\CP(k-1))$.
Let $g\colon \Lambda(x,z) \rightarrow \Apl(\CP(k-1))$ be any map such that
$[g(x)]=a$ and $f\colon\Lambda (z)\rightarrow \Apl(S^{2k-1})$ be any map
such that $[f(z)]=-u_S$. Note that these equations determine $g$ and $f$ up to
homotopy.
To see that the bottom square of (\ref{hard}) commutes
consider the following cube:
$$
\xymatrix
{
\hat Q \otimes \Lambda(x,z) \ar[dd]_{\theta'} \ar[rr]^{\rm proj}
\ar[dr]_{\rm proj} && \hat Q \otimes \Lambda(z)
\ar'[d]^{\mu'}  [dd]\ar[dr]^{\rm proj} \\
& \Lambda(x,z) \ar[dd]_(0.3){g} \ar[rr]^(0.3){\rm proj} && \Lambda(z) \ar[dd]^f \\
\Apl(P\nu) \ar'[r] [rr]_{q^*} \ar[dr]_{{\rm inc}^*}
&& \Apl(\bd T) \ar[dr]^I \\
& \Apl(\CP(k-1)) \ar[rr]_{\Apl(\eta)} && \Apl(S^{2k-1}). }
$$
The back face is the one we wish to show is homotopy commutative.
The top and bottom faces clearly commute. The front, right and left faces
are homotopy commutative by Lemmas
\ref{B1}, \ref{B3} and \ref{B2}. So we get that $I\mu'({\rm proj}) \simeq
Iq^*\theta'$. Next consider the coaction sequence
(\cite[Chapter I Section 3 Proposition 4]{Qui})
associated to the map from
the cofibration sequence of CDGA
$$
\Lambda (s^{-1}z) \rightarrow \hat Q \otimes \Lambda(x)
\rightarrow \hat Q \otimes \Lambda(x,z) \rightarrow \Lambda (z)
$$
into $\Apl(\bd T) \rightarrow \Apl(S^{2k-1})$. We get a
commutative diagram of sets
$$
\xymatrix
{
[\Lambda (z),\Apl(\bd T)] \ar[d]^{I_*} \ar[r]^-{p^*} &
[\hat Q\otimes \Lambda (x,z),\Apl(\bd T)] \ar[d]^{I_*}
\ar[r]^{{\rm inc}^*} &
[\hat Q \otimes \Lambda (x),\Apl(\bd T)] \ar[d]^{I_*} \\
[\Lambda (z),\Apl(S^{2k-1})]  \ar[r]_-{p^*} &
[\hat Q\otimes \Lambda (x,z),\Apl(S^{2k-1})] \ar[r]_{{\rm inc}^*} &
[\hat Q \otimes \Lambda (x),\Apl(S^{2k-1})]
}
$$
By \cite[Chapter I Section 3 Proposition 4']{Qui} the rows are exact in the sense
that, if $inc^*f=inc^*g$ then there exists $\alpha\in [\Lambda(z),\_]$ such that
$\alpha\cdot f=g$, where $\_\cdot\_$ denotes the action of the
group $[\Lambda(z),\_]$ on $[\hat Q\otimes \Lambda (x,z), \_]$.
Note that our cofibration sequence is a model of a fibration sequence
$$
S_{\BQ}^{2k-1}=K(\BQ,2k-1)\rightarrow P\nu_{\BQ}\rightarrow
V_{\BQ}\times K(\BQ,2) \rightarrow K(\BQ,2k)
$$
and at the space level the coaction sequence is the mapping class sequence.
(see \cite[Chapter 2]{Switzer}).
For any cofibration sequence $A\rightarrow B\rightarrow C\rightarrow \Sigma A$ and
group object $G$ in a pointed model category the coaction
$\_\cdot\_\ \colon [\Sigma A,G]\times [C,G]\rightarrow [C,G]$ from the cofibration sequence
is compatible with the group action $\phi\colon [C,G]\times [C,G]\rightarrow [C,G]$ induced
by the multiplication on $G$. In particular for any $\alpha\in [\Sigma A,G]$ and
$f,g\in [C,G]$:
\begin{equation}\label{star2}
\alpha\cdot \phi(f,g)=\phi(\alpha\cdot f,g).
\end{equation}
Consider the diagram
$$
\xymatrix
{
\hat Q \ar[dd]_{\beta'} \ar[ddrr]^{\theta'|_{\hat Q}} \\
\\
\Apl(T) \ar[r]^{\sigma^*} \ar[d]_{i^*} &
\Apl(V) \ar[r]^{(\pi')^*} \ar[dl]_{(\pi|_{\bd T})^*} &
\Apl(P\nu) \ar[dll]^{(q|_{\bd T})^*} \\
\Apl(\bd T)
}
$$
The top triangle commutes up to homotopy by Lemma \ref{B2}.
The bottom left triangle commutes up to homotopy since
$\sigma\colon V \rightarrow T$ is the zero section of $\pi$ and the
bottom right triangle commutes since it is $\Apl$ of Diagram (\ref{quotie}).
Thus Diagram (\ref{hard}) commutes up to homotopy when restricted to $\hat Q$.
\par
Since $\theta'(x)=a$ and by Lemma \ref{lemma-qa0} $q^*(a)=0$, Diagram (\ref{hard}) restricted to
$\Lambda(x)$ commutes and so restricted to
$\hat Q \otimes \Lambda(x)$ commutes up to homotopy.
Thus as homotopy classes
${\rm inc}^*\mu'(proj)={\rm inc}^*q^*\theta'$, and
there exists $\alpha\in [\Lambda(z),\Apl (\bd T)]$ such that
$\alpha \cdot \mu'{\rm proj}=q^*\theta'$.
Since the action is natural in the second variable
$I\alpha \cdot I\mu'{\rm proj}=Iq^*\theta'$.
As observed above $I\mu'(proj)=Iq^*\theta'$, so $I\alpha\cdot I\mu'(proj)=I\mu'(proj)$
and so
$$
I\alpha\cdot 0=
I\alpha\cdot ((I\mu'(proj))(I\mu'(proj))^{-1})=(I\alpha\cdot I\mu'(proj))(I\mu'(proj))^{-1}
$$
$$
=(I\mu'(proj))(I\mu'(proj))^{-1}
=0
$$
with the first equality following from Equation (\ref{star2}) since $\Apl(S^{2k-1})$ is a group
object in CDGA. This implies that $p^*(I\alpha)=0$.

Since $H^{2k}(\hat Q)=0$ by looking at models we get that
$\pi_{2k-1}(S^{2k-1})\otimes \BQ \rightarrow \pi_{2k-1}(P\nu)\otimes \BQ$
is injective. Thus
$p^*:[\Lambda (z),\Apl(S^{2k-1})]  \rightarrow
[\hat Q\otimes \Lambda (x,z),\Apl(S^{2k-1})]$ is injective, and so $I\alpha=0$.

Since $Q\otimes \Lambda (z) \rightarrow \Lambda(z)$ with $dz=0$ models
$\Apl(\bd T)\rightarrow \Apl(S^{2k-1})$ and $H^{2k-1}(Q)=0$, $I$ induces an
isomorphism on $H^{2k-1}$ and so
$I_*:[\Lambda (z),\Apl(\bd T)]\rightarrow [\Lambda (z),\Apl(S^{2k-1})]$
is an isomorphism.
Thus $\alpha=0$ and so $\mu'{\rm proj}=q^*\theta'$ as homotopy classes.
\cqfd
\begin{prop}
\label{prop-modelprojxi}
Let $\nu:\BC^k\to E\to V$ be a complex vector bundle
of rank $k$ such that $2k\geq \dim V +2$ and
$Q$ be a CDGA weakly equivalent to $\Apl(V)$.
Set $\gamma_0=1$, and for $1\leq i\leq k-1$
let $\gamma_i\in Q^{2i}\cap\ker d$ be cocycle representatives of the
Chern classes $c_i(\xi)\in H^{2i}(V;\BQ)$.
Using the same notation as at the beginning of Section \ref{section-modelblowup} the diagram
$$
\xymatrix{
(Q\otimes\Lambda z,D z=0)
&&\ar[ll]_-\phi (Q\otimes\Lambda(x,z);
Dx=0,Dz=\sum_{i=0}^k\gamma_i x^{k-i})
\\&\ar@{^(->}[ul]\ar@{^(->}[ur]Q
}
$$
with $\phi|_Q={\rm inc}$, $\phi(x)=0$, $\phi(z)=z$,
$|x|=2$ and $|z|=2k-1$
is a CDGA model of the diagram
$$
\xymatrix{
E_0\ar[rr]^q\ar[rd]_{\pi_0}&&P\nu\ar[ld]^{\pi'}\\&V.
}
$$

\end{prop}
\proof
Observe that Lemma \ref{gorgle} also works in the case that
$f$ and $f'$ are identity maps. Thus the proposition follows
from Lemma \ref{(I)} by applying Lemma \ref{gorgle} twice.
\cqfd
Notice that since $2k\geq \dim V +2$, $\gamma_k=0$, and so
$Dz=\Sigma_{i=0}^{k-1} \gamma_ix^{k-i}$.
Note that in this section we never used the hypothesis that $V$ is a manifold; only
paracompactness is needed to define the classifying
map of the tautological bundle $\lambda$ and we also had some dimension restrictions.
\section{The model of the blow-up}
\label{when?}
For all this section we have the following assumptions and notation.
Let $f\colon V\hookrightarrow W$
be an embedding of connected closed manifolds of
codimension $2k$ with a fixed complex structure on its normal bundle $\nu$.
Set $n=\dim (W)$ and $m=\dim(V)$ as before
and $r=n-m=2k$. Assume that $n\geq 2m+3$ and that $H^1(f)$ is injective.
Let $\phi\colon R\to Q$ be a CDGA-model of the embedding
$f$ such that $R^{\geq n+1}=0$ and $Q^{\geq m+2}=0$.
Suppose we have fixed some shriek map $\phi^!\colon
s^{-2k}Q\to R$ of $R$-dgmodules.
Recall that all this can be done
using Proposition \ref{prop-shriekexists}.
Set $\gamma_0=1$ and for $1\leq i\leq k-1$ let $\gamma_i\in Q^{2i}\cap\ker d$
be representatives of the Chern classes $c_i(\nu)\in H^{2i}(V;\BZ)$.
We will also use the notation of Diagram (\ref{last cube}) of Section \ref{homo-type} below.
\par
In this section we pull together what we have done up to now and give our model of the
blow-up. In Section \ref{homo-type} we recall a definition of the
blow-up $\widetilde W$ and of the
map $\tilde \pi:\widetilde W\rightarrow W$
which are suitable for studying homotopy theoretic questions.
In \ref{desc of model} we define a CDGA model ${\cal B}(R,Q)$ which depends on the model
$\phi\colon R\rightarrow Q$ of the embedding $f\colon V\rightarrow W$, the shriek map
$\phi^!\colon s^{-2k}Q \rightarrow R$ and the representatives $\gamma_i$ of the Chern classes.
In Section \ref{2 equiv diags} we prove the equivalence of two diagrams of CDGAs
(Lemma \ref{lemma-2 equiv}). The first is constructed from our models
$R$ and $Q$, $\phi$, $\phi^!$ and the $\gamma_i$ and the second
comes from taking $\Apl$ of Diagram
(\ref{last cube}). In Section \ref{two pairs} we take pullbacks derived
from the equivalent diagrams of Lemma \ref{lemma-2 equiv} and show that they are equivalent
to ${\cal B}(R,Q)$ (Lemma \ref{almost}) and $\Apl(\widetilde W)$(Lemma \ref{pen}).
Putting these together we prove our main theorem
(Theorem \ref{lemma-modelblowup}) that
${\cal B}(R,Q)$ is a model for
$\Apl(\widetilde W)$. Under a certain nilpotence condition this
also implies that the rational homotopy type of the blow-up
along $f:V \rightarrow W$ depends only on the rational homotopy class of the map $f$ and
the Chern classes of the normal bundle of $V$ in $W$ (Corollary \ref{UBC}).
\subsection{The homotopy type of the blow-up}
\label{homo-type}
Consider the following cubical diagram with a triangle added to the bottom face.
As in Section \ref{sec-compl}, $T$ is a tubular neighborhood of $V$ in $W$ and
$B=\overline {W \setminus T}$. The projective bundle $P\nu$ was defined in
Section \ref{Chern} as were the projection maps $q$, $\pi_0$ and $\pi'$.
The maps $i$, $k$, $l$ and $\pi$ were first seen in Section \ref{sec-Thom}.
So we have
come across all the maps in the diagram except $\tilde \pi$.
The space $\widetilde W$ is defined as the pushout of the top face.
\begin{equation}
\label{last cube}
\xymatrix
{
\bd T\ar[rr]^k\ar[rd]^q\ar@{=}[dd]
&&B\ar[rd]\ar@{=}'[d] [dd]|{\,}\\
&P\nu\ar[rr]\ar[dd]_(0.3){\pi'}&&\widetilde{ W}\ar[dd]^{\tilde\pi}\\
\bd T\ar@{^{(}->}[d]_i\ar'[r][rr]^k|{\,}\ar[rd]_{\pi_0}&&B\ar[rd]^l\\
T\ar[r]_\pi&V\ar[rr]^f&&W
}
\end{equation}
We know that when we replace $f\pi$ by $j$ the outside bottom quadrilateral is a pushout
(see Section \ref{sec-Thom}).
In fact the bottom face of the cube is also a pushout
since $\pi$ is a deformation retraction.
However the map between these pushouts induced by $\pi$ is a homotopy equivalence
and not a homeomorphism.
\begin{defin}
The {\em blow-up of $W$ along $V$} is the pushout $\widetilde W$ and the map
$\tilde \pi:\widetilde W\rightarrow W$ is the map induced by $\pi'$ between
pushouts comprising the top and bottom faces of the above cube.
\end{defin}
This definition of blow-up is equivalent to those of
\cite{GriffithsHarris} and \cite[7.1]{McDuffSalamon}.
\subsection{Description of the model ${\cal B}(R,Q)$ for the blow-up}
\label{desc of model}

We now construct a CDGA,
${\cal B}(R,Q)$, which will be a CDGA model of the blow-up $\widetilde{W}$ of $W$
along $V$. We also
define a morphism $\iota(R,Q)\colon R\to {\cal B}(R,Q)$ that will model the projection
$\tilde\pi\colon \widetilde{W}\to W$.
\par

Let $x$ and $z$ be generators such that
$|x|=2$ and $|z|=2k-1$ and denote by $\Lambda^+(x,z)$ the augmentation ideal of the
free graded commutative algebra $\Lambda(x,z)$. The CDGA ${\cal B}(R,Q)$ is of the form:

$$
{\cal B}(R,Q)=\left(R\oplus Q \otimes \Lambda^+(x,z),D\right).
$$

The graded commutative algebra structure
on ${\cal B}(R,Q)$ is induced by the multiplications
on $R$ and $Q\otimes \Lambda^+(x,z)$, and by the $R$-module
structure on the free $Q$-module $Q\otimes\Lambda^+(x,z)$ induced by the algebra map
$\phi:R\rightarrow Q$. More explicitly for $r\in R$, $q\in Q$ and
$w\in \Lambda^+(x,z)$,

$$
r\cdot (q\otimes w)=(\phi(r)\cdot q)\otimes w,
$$

$$
(q\otimes w)\cdot r=(-1)^{(|w|+|q|)|r|}(\phi(r)\cdot q)\otimes w.
$$

Let $d_R$ and $d_Q$ denote the differentials on $R$ and $Q$ respectively,
$r\in R$ and $q\in Q$.
The differential $D$ on ${\cal B}(R,Q)$ is determined by the Leibnitz law and the formulas

\begin{eqnarray*}
D(r)&=&d_R r\\
D(q\otimes x)&=& d_Q q \otimes x\\
D(q\otimes z)&=& d_Q q\otimes z+(-1)^{|q|}\left(\phi^!(s^{-2k}q)+\sum_{i=0}^{k-1}
(q \cdot \gamma_i) \otimes x^{k-i}\right) \\
\end{eqnarray*}

There is an obvious inclusion morphism

$$
\iota(R,Q):R\to {\cal B}(R,Q).
$$

Notice that ${\cal B}(R,Q)$ actually depends
not only on $R$ and $Q$ but also on the $\gamma_i$,
$\phi$ and $\phi^!$.
These are implicit and not included in the notation.

\begin{lemma}
${\cal B}(R,Q)$ is a CDGA and $\iota(R,Q):R\rightarrow {\cal B}(R,Q)$
is a CDGA morphism.
\end{lemma}

\proof
For $q,q'\in Q$ we have $(q\otimes z)\cdot(q'\otimes z)=0$, therefore we need to check that
the Leibnitz law applied to $D((q\otimes z)\cdot(q'\otimes z))$ gives zero.
This follows because, for degree reasons, $\phi^!\phi=0$.
To check the rest of the definition of CDGA is straightforward.
\cqfd
\subsection{Two equivalent diagrams}
\label{2 equiv diags}
\begin{lemma}
\label{lemma-2 equiv}
Consider the map $e$ from Lemma \ref{esc}
and the relative Sullivan algebra $(Q\otimes \Lambda(x,z),D)$ from
Equation (\ref{equ-Qxz}).
The CDGA diagram
$$
\xymatrix
{
Q \ar[r]^-{\rm inc} & Q\otimes \Lambda(x,z) \ar[r]^{e {\rm proj}}
& Q\oplus ss^{-r}Q &
R\oplus_{\phi^!} ss^{-r} Q \ar[l]_{\phi\oplus id}
}
$$
is weakly equivalent to the diagram
$$
\xymatrix
{
\Apl(V) \ar[r]^{(\pi')^*} & \Apl(P\nu) \ar[r]^{q^*}
& \Apl(\bd T) & \Apl(B) \ar[l]_{k^*}
}
$$
\end{lemma}

\proof
We begin by fixing a common model
$\hat \phi\colon \hat R \rightarrow \hat Q$ as in Section \ref{sec-supps}.
All the notation is given in the lemmas we refer to.
Lemma \ref{igog} gives us a commutative diagram of CDGAs

\begin{equation}\label{gleep}
\xymatrix
{
Q\oplus ss^{-r}Q \ar[d]_{e^{-1}} & R\oplus_{\phi^!} ss^{-r} Q
\ar[l]_{\phi\oplus id} \ar[d]^{\xi}\\
Q\otimes \Lambda(z) & A \ar[l]_{\kappa} \\
\hat Q\otimes \Lambda(z)\ar[u]^{\beta\otimes id}\ar[d]_{\mu'}
& \hat A \ar[u]_{\zeta}\ar[d]^{\zeta'}\ar[l]_{\psi}\\
\Apl(\bd T) & \Apl (B) \ar[l]_{k^*}
}
\end{equation}

with all vertical arrows being weak equivalences.
Choose representatives
$\hat \gamma_i\in \hat Q$ of $c_i(\nu)$ such that
$\beta(\hat\gamma_i)=\gamma_i$. This can be done since
$\beta$ is an acyclic fibration.
Next Lemmas \ref{B2} and \ref{(I)} imply that we have a homotopy
commutative diagram of CDGA
\begin{equation}\label{glop}
\xymatrix{
Q \ar[d]_= \ar[r]^-{\rm inc} & Q \otimes \Lambda(x,z) \ar[d]_= \ar[r]^{e {\rm proj}}
& Q \oplus_{\phi\phi^!}ss^{-r} Q \ar[d]_{e^{-1}}\\
Q \ar[r]^-{\rm inc} & Q \otimes \Lambda(x,z)  \ar[r]^{\rm proj}
& Q \otimes \Lambda(z) \\
\hat Q \ar[u]^{\beta} \ar[r]^-{\rm inc} \ar[d]_{\sigma^*\beta'} & \hat Q\otimes \Lambda(x,z)
\ar[u]^{\beta\otimes id} \ar[r]^{\rm proj} \ar[d]_{\theta'}
& \hat Q\otimes \Lambda(z)\ar[d]_{\mu'}\ar[u]^{\beta\otimes id}\\
\Apl(V) \ar[r]_{(\pi')^*} & \Apl(P\nu) \ar[r]_{q^*} & \Apl(\bd T)
}
\end{equation}
with all vertical arrows being weak equivalences.
Note that $Q \oplus_{\phi\phi^!}ss^{-r} Q=Q \oplus ss^{-r} Q$ since $Q^{\geq m+2}=0$
entails $\phi\phi^!=0$, so we can glue Diagram
(\ref{gleep}) to the right of Diagram (\ref{glop})
and apply Lemma \ref{gorgle} three times to
get the desired result.
\cqfd
\subsection{Two pairs of pullbacks and the proof the the main theorem}
\label{two pairs}
\begin{lemma}
\label{almost}
In CDGA, the pullback of
\begin{equation}
\label{Russ}
\xymatrix
{
& Q \ar[d]^{\rm inc}\\
R\oplus_{\phi^!} ss^{-r} Q \ar[r]_{\phi\oplus id} & Q\oplus ss^{-r} Q
}
\end{equation}
is $R$ and the pullback of
\begin{equation}
\label{Tom}
\xymatrix
{
& Q \otimes \Lambda(x,z) \ar[d]^{e{\rm proj}}\\
R\oplus_{\phi^!} ss^{-r} Q \ar[r]_{\phi\oplus id} & Q\oplus ss^{-r} Q
}
\end{equation}
is ${\cal B}(R,Q)$.
\par
The map between the pullbacks induced by the inclusion
$Q \rightarrow Q\otimes \Lambda (x,z)$ is
$\iota(R,Q): R\rightarrow {\cal B}(R,Q)$.
For both of the diagrams above the map from the homotopy pullback to the
pullback is a weak equivalence. Thus the map between the homotopy pullbacks induced
by the inclusion $Q\rightarrow Q\otimes \Lambda(x,z)$ is equivalent to
$\iota(R,Q):R\rightarrow {\cal B}(R,Q)$.
\end{lemma}
\proof
Consider the maps
$$
g\colon {\cal B}(R,Q) \rightarrow R\oplus_{\phi^!} ss^{-2k} Q
$$
determined by the equations
$$
g(r,0)=(r,0)
$$
$$
g(0,q\otimes z)=(0,(-1)^{|q|}ss^{-2k}q)
$$
$$
g(0,q\otimes x^i\otimes z^{\epsilon})=(0,0)
\ {\rm if} \ i>0
$$
and
$$
g'\colon {\cal B}(R,Q) \rightarrow Q\otimes \Lambda(x,z)
$$
$$
(r,q\otimes x^i\otimes z^{\epsilon}) \mapsto
\phi(r) + q\otimes x^i\otimes z^{\epsilon}
$$
Since $(\phi\oplus {\rm id})g=e ({\rm proj}) g'$, $g$ and $g'$
determine a map
$$
h\colon {\cal B}(R,Q) \rightarrow {\rm pullback \ of \
Diagram} \ (\ref{Tom}).
$$
 Because
the forgetful functor from CDGA to
graded modules commutes with taking pullbacks
it is easy to check that $h$ is an isomorphism.
Similarly $R$ is the pullback of Diagram \ref{Russ} and
$\iota(R,Q):R\rightarrow {\cal B}(R,Q)$ is the induced map
between the pullbacks.
\par
Next we will show that these pullbacks are indeed homotopy pullbacks.
The short exact sequence of differential graded modules
$$
0\rightarrow R\rightarrow Q \oplus (R\oplus_{\phi^!} ss^{-r} Q)
\rightarrow Q\oplus ss^{-r} Q\rightarrow 0
$$
gives rise to a Mayer-Vietoris long exact sequence on homology. This maps into
the corresponding long exact sequence for the homotopy pullback. Thus the map from $R$
to the homotopy pullback of Diagram \ref{Russ}
is an equivalence by the five lemma. The map from
${\cal B}(R,Q)$ into the homotopy pullback of
Diagram \ref{Tom} is similarly a weak equivalence. The fact that the
induced map between the homotopy pullbacks is equivalent to
$\iota(R,Q): R\rightarrow {\cal B}(R,Q)$ follows by naturality.
\cqfd
\begin{lemma}\label{pen}
Recalling the cube (\ref{last cube}),
the homotopy pullback of
\begin{equation}
\label{first}
\xymatrix
{
& \Apl(V) \ar[d]^{(\pi i)^*}\\
\Apl(B)\ar[r]_{k^*} & \Apl(\bd T)
}
\end{equation}
is quasi-isomorphic to $\Apl(W)$ and the homotopy pullback of
\begin{equation}
\label{second}
\xymatrix
{
& \Apl(P\nu) \ar[d]^{q^*}\\
\Apl(B)\ar[r]_{k^*} & \Apl(\bd T)
}
\end{equation}
is quasi-isomorphic to $\Apl(\widetilde W)$.
The map between the homotopy pullbacks induced by
${\pi'}^*:\Apl(V)\rightarrow \Apl(P\nu)$
is weakly equivalent to
$\tilde\pi^*:\Apl(W)\rightarrow \Apl(\widetilde W)$.
\end{lemma}
\proof
Recall (\ref{last cube}) the following cube in which the bottom
and the top faces are push-outs :
$$
\xymatrix
{
\bd T\ar[rr]\ar[rd]^q\ar@{=}[dd]
&&B\ar[rd]\ar@{=}'[d] [dd]\\
&P\nu\ar[rr]\ar[dd]_(0.3){\pi'}&&\widetilde{ W}\ar[dd]^{\tilde\pi}\\
\bd T\ar '[r][rr]\ar[rd]_{\pi_0}&&B\ar[rd]\\
&V\ar[rr]&&W
}
$$
Let $U$ be the homotopy pullback of (\ref{first}) and
$f:\Apl (W)\rightarrow U$ be the induced map. We have maps
between Mayer-Vietoris sequences
$$
\xymatrix
{
H^{*-1}(\Apl (\bd T)) \ar[r] \ar[d]^=  & H^*(\Apl (W)) \ar[d]^{H(f)} \ar[r]
& H^*(\Apl (V)) \oplus H^*(\Apl (B)) \ar[d]^=\ar[r] & H^*(\Apl (\bd T))\ar[d]^=\\
H^{*-1}(\Apl (\bd T)) \ar[r]  & H^*(U) \ar[r]
& H^*(\Apl (V)) \oplus H^*(\Apl (B))\ar[r] & H^*(\Apl (\bd T)),
}
$$
hence $f$ is a weak equivalence by the five lemma.
Similarly $\Apl (\tilde W)$ is weakly equivalent to
the homotopy pullback of (\ref{second}). The fact that the induced map is weakly equivalent to
$\tilde \pi^*$ follows by naturality.
\cqfd
\subsection{Proof of main theorem}\label{mainthing}
Here is the main result of the paper.
\begin{thm}
\label{lemma-modelblowup}
Let $f\colon V\hookrightarrow W$
be an embedding of connected closed manifolds of
codimension $2k$ with a fixed complex structure on its normal bundle $\nu$.
Set $n=\dim (W)$ and $m=\dim(V)$ as before
and $r=n-m=2k$. Assume that $n\geq 2m+3$ and that $H^1(f)$ is injective.
Let $\phi\colon R\to Q$ be a CDGA-model of the embedding
$f$ such that $R^{\geq n+1}=0$ and $Q^{\geq m+2}=0$.
Let $\phi^!\colon s^{-2k}Q\to R$ be a shriek map of $R$-dgmodules.
Set $\gamma_0=1$ and for $1\leq i\leq k-1$ let $\gamma_i\in Q^{2i}\cap\ker d$
be representatives of the Chern classes $c_i(\nu)\in H^{2i}(V;\BZ)$.

The CDGA
$$
{\cal B}(R,Q)=\left(R\oplus Q \otimes \Lambda^+(x,z),D\right)
$$
defined in Section \ref{desc of model} is a CDGA model of $\Apl(\widetilde{W})$
where $\widetilde W$ is the blow-up
of $W$ along $V$. Also
$\iota(R,Q)\colon R\hookrightarrow {\cal B}(R,Q)$
is a CDGA model of $\Apl(\tilde\pi)\colon \Apl(W)\to\Apl(\widetilde{W})$.
\end{thm}
\proof
The theorem follows directly from Lemmas \ref{inducted},
\ref{lemma-2 equiv}, \ref{almost} and \ref{pen}.
\cqfd
Note that if $n\geq 2m+3$ and $H^1(f)$ is injective,
by Proposition \ref{prop-shriekexists}
any model of $f$ can be replaced by
one satisfying the hypotheses of the theorem.
\begin{cor}\label{UBC}
With the hypotheses of Theorem \ref{lemma-modelblowup}, if we assume that
$V$, $W$ and the blow-up $\widetilde W$ are nilpotent spaces then the rational homotopy
type of $\widetilde W$ is determined by the rational homotopy class of $f$ and by
the rational Chern classes $c_i(\nu)\in H^{2i}(V;\BQ)$ of the normal bundle.
\end{cor}
\proof This is an immediate application of Sullivan's theory \cite{Sullivan}
to Theorem \ref{lemma-modelblowup}. In particular (see \cite[Section 9]{BG}) since
$V$ and $W$ are nilpotent if $f,g\colon V\rightarrow W$ are rationally homotopic then
$\Apl (f)$ and $\Apl (g)$ have the same models $R\rightarrow Q$. Since the Chern classes
are also the same, ${\cal B}(R,Q)$ with model $\Apl(\widetilde W)$ along both embeddings.
Since $\widetilde W$ is nilpotent $\Apl(\widetilde W)$ determines the rational homotopy
type of $\widetilde W$.
\cqfd
We have a homotopy pushout
$$
\xymatrix
{
P\mu \ar[d] \ar[r] & \widetilde W \ar[d] \\
V \ar[r] & W
}
$$
and so by the Van Kampen theorem $\pi_1(\widetilde W) \rightarrow \pi_1(W)$ is
an isomorphism. So if $W$ is nilpotent then $\pi_1(\widetilde W)=\pi_1(W)$ is
a nilpotent group but the action of $\pi_1$ on the homotopy groups of the universal
cover of $\widetilde W$ may not be nilpotent. In certain cases Rao \cite{Rao} has
determined when a homotopy pushout is nilpotent, however it seems difficult to see if
they apply in our situation.
However notice that the nilpotence condition in the corollary is automatically satisfied if
$V$ is nilpotent and $W$ is simply connected because then $\widetilde W$ is also
simply connected. Hence we have proved our first
theorem from the introduction:
\begin{thm}
\label{blorgy}
Let $f:V \rightarrow W$ be an embedding of smooth closed orientable manifolds such that
$W$ is simply connected and $V$ is nilpotent and oriented.
Suppose that the normal bundle $\nu$ is equipped with the structure of a complex
vector bundle and assume that $\dim W\geq 2\dim V +3$.
Then the rational homotopy type of the blow-up
of $W$ along $V$, $\widetilde W$ can be explicitly determined
from the rational homotopy type
of $f$ and from the Chern classes $c_i(\nu)\in H^{2i}(V;\BQ)$.
\end{thm}
Even without the nilpotence condition we can still
determine a model for $\Apl(\widetilde W)$.
\begin{cor}
With the hypotheses in the first paragraph of Section \ref{when?},
a CDGA model of $\Apl(\widetilde W)$ is determined by any model $\Apl(f)$ and by
the rational Chern classes $c_i(\nu)\in H^{2i}(V;\BQ)$ of the normal bundle.
\end{cor}
As we will see in the next section this is enough to determine $H^*(\widetilde W)$
(Theorem \ref{algcoho}).
\section{Applications}
\label{section-application}
In this section we apply our model of the blow-up to three situations.
First, after some preliminaries in Section \ref{prelim} we describe
in Section \ref{blow in cp} the model of the blow-up of $\CP(n)$
along a submanifold.
This is somewhat simpler than our general model
since here the shriek map is more easily
described. Secondly, Section \ref{McDuff}
examines the model of the example of McDuff of the blow-up of $\CP(n)$
along the Kodaira-Thurston manifold.
We recover directly the Babenko-Taimanov result that this $\widetilde \CP(n)$ has a
nontrivial Massey product and is thus nonformal. Finally in Section \ref{coho alg} we describe
the cohomology algebra of certain blow-ups. Some other special cases have been studied by
Gitler \cite{Gitler}.
In Section \ref{vlast} we calculate the rational homotopy type of blow-ups of
$\CP(5)$ along $\CP(1)$. It turns out that there are infinitely many rationally
inequivalent ones.
\subsection{Preliminaries on symplectic manifolds}
\label{prelim}
A {\em symplectic form} on a $2n$-dimensional manifold $M$ is a nondegenerate
closed differential 2-form $\omega$. Thus $\omega^n$ is a volume form.
The pair $(M,\omega)$ is called a {\em symplectic manifold}
(see \cite{McDuffSalamon}).
The form $\omega$ is called {\em integral} if it belongs to the image
$H^2(M,\BZ) \rightarrow H^2(M,\BR)\cong H^2(\Omega^*(M))$
where $\Omega^*(M)$ is the de Rham complex of differential forms.
\par
It is well known that to any symplectic form $\omega$ we can associate an almost complex
structure on the tangent bundle of $M$. This almost
complex structure in turn induces a preferred
orientation on $M$ and hence
a generator $u_M\in H^{2n}(M,\BZ)$.
\begin{defin}\label{def-mult}
Define the real number $l_M$ by the Equation $[\omega^n]=l_M\cdot u_M\in H^{2n}(M;\BR)$.
\end{defin}
Since $\omega^n$ is a volume form and $u_M$ is induced from the almost
complex structure associated to $\omega$
we know that $l_M$ is a positive
real number. If $\omega$ is integral then $l_M$ is a positive integer.
We cannot always choose an integral $\omega$ so that $l_M=1$. For example if
$[\omega]\in H^2(S^2\times S^2,\BZ)$ then $[\omega]^2\in H^4(S^2\times S^2,\BZ)$
is always divisible by $2$.
The number $l_M$ will appear in the formula when we blow-up of $\CP(n)$ along $M$.
\par
An example of a symplectic manifold is given by the complex projective space
$\CP(n)$ equipped with the 2-form $\omega_0$ associated to the Fubini-Study metric
(\cite[Example 4.21]{McDuffSalamon} or \cite[page 31]{GriffithsHarris}). It is classical
(see for example \cite[pages 30-32 and 144-150]{GriffithsHarris}) that
$$
[\omega_0]=-c_1(\gamma^1)\in H^2(\CP(n);\BZ)
$$
In particular $\omega_0$ is integral and $l_{\CP(n)}=1$.
\subsection{Blow-ups of \CP(n)}
\label{blow in cp}
Let $(M,\omega)$ be a symplectic manifold. By a {\em symplectic embedding}
\begin{equation}\label{equ-symplemb}
f\colon (M,\omega) \rightarrow (\CP(n),\omega_0)
\end{equation}
we mean a smooth embedding such that $f^*(\omega_0)=\omega$.
It can be proved that $f$ respects the almost complex structures associated to
$\omega$ and $\omega_0$, therefore the normal bundle of $f$ admits
a natural complex structure
and we can consider the blow-up $\widetilde \CP(n)$ of $\CP(n)$ along $M$.
Notice also that, since $\omega_0$ is integral, if a symplectic manifold
$(M,\omega)$ symplectically embeds
in $(\CP(n),\omega_0)$ then $\omega$ is integral. Moreover the converse
is almost true thanks to a theorem
of Tischler and Gromov \cite{Gromov}, \cite{Tischler}:
if $M$ is a manifold equipped with an integral symplectic form $\omega$ then,
for $n$ large enough,
there exists a symplectic embedding as in Equation (\ref{equ-symplemb}).
\par
The first application of our model is to the blow-up of $\CP(n)$
along a symplectically embedded
submanifold. To simplify the description we will use:
\begin{lemma}
\label{lemma-complomega}
Let $Q$ be a CDGA-model of some closed symplectic manifold of dimension
$2m$ with symplectic form $\omega\in Q^2\cap\ker d$
and suppose that $Q^{\geq 2m+2}=0$.
Consider the
$\left(\Lambda(\omega)/(\omega^{m+1}),0\right)$-dgmodule structure on
$Q$ induced by multiplication by $\omega$.
Then there exists a sub-$\left(\Lambda(\omega)/(\omega^{m+1}),0\right)$-dgmodule
$I\subset Q$
such that $Q=I\oplus\Lambda(\omega)/(\omega^{m+1})$ as
$\left(\Lambda(\omega)/(\omega^{m+1}),0\right)$-dgmodules.
\\Moreover if $Q^{2m}=\BQ\cdot \omega^m$ then this
 sub-$\left(\Lambda(\omega)/(\omega^{m+1}),0\right)$-dgmodule is unique.
\end{lemma}
\proof
Take a complementary vector subspace $S$ of $\BQ \cdot \omega^m\oplus dQ^{2m-1}$ in $Q^{2m}$
and set $I_{(2m)}=S\oplus Q^{2m+1}$.
We know $I_{(2m)}$ is a sub-$\left(\Lambda(\omega)/(\omega^{m+1}),0\right)$-dgmodule of
$Q$ since $Q^{\geq 2m+2}=0$.
Suppose that for some $k\leq m$ and
for each $j>k$ we have already defined
a differential submodule $I_{(2j)}\subset Q$ such that
$Q^{\geq 2j}=\BQ\{\omega^j,\omega^{j+1},\ldots,\omega^m\}\oplus I_{(2j)}$.
We define $I_{(2k)}$ as follows.
Consider the morphism
$$
\lambda_k\colon Q^{2k}\stackrel{\cdot\omega}\to Q^{2k+2}
\stackrel{\pr}\to Q^{2k+2}/I_{(2k+2)}\cong\BQ\cdot\omega^{k+1}.
$$
Set $I_{(2k)}=\ker(\lambda_k)\oplus Q^{2k+1}\oplus I_{(2k+2)}$.
Since $\alpha\in \ker \lambda_k$ implies that $\alpha\cdot\omega\in I_{(2k+2)}$,
it is straightforward to check that $I_{(2k)}$ is
a sub-$\left(\Lambda(\omega)/(\omega^{m+1}),0\right)$-dgmodule of $Q$ such that
$Q^{\geq 2j}=\BQ\{\omega^j,\omega^{j+1},\ldots,\omega^m\}\oplus I_{(2j)}$.
Finally $I=I_{(0)}$ is the desired submodule of $Q$.
\par
If $Q=I\oplus\Lambda(\omega)/(\omega^{m+1})$ as
$\left(\Lambda(\omega)/(\omega^{m+1}),0\right)$-dgmodules then
$\ker\lambda_j$ must be included in $I$. Since $Q^{\geq 2m+2}=0$,
we must have that $I_{(2m)}=S\oplus Q^{2m+2}$ and then we can show by induction
that $I^{\geq 2k}=I_{(2k)}=\ker (\lambda_k) \oplus Q^{2k+1} \oplus I_{2k+2}$ for all
$k<m$. This implies that
$I$ is completely determined by the choice of $S$.
Therefore when $Q^{2m}=\BQ\cdot\omega^m$ the ideal $I$ is unique.
\cqfd
Let $(M,\omega)$ be a symplectic manifold of dimension $2m$,
$f\colon M \rightarrow \CP(n)$ be a symplectic
embedding and $\nu$ be the normal bundle of the embedding.
Let $(Q,d)$ be a CDGA-model of $M$ such that $Q^{\geq 2m+2}=0$.
Let $\omega\in Q^2\cap\ker d$ be a representative of the
symplectic form and $\gamma_i\in Q^{2i} \cap \ker d$ be
representatives of the Chern classes $c_i(\nu)\in H^{2i}(M)$
of the normal bundle of $M$ in $\CP(n)$.
Then $R=(\Lambda(a)/(a^{n+1}),0)$ is a CDGA model of $\CP(n)$
with $[a]=[\omega_0]\in H^2(\CP(n))$
and $a^n$ represents the orientation class of $\CP(n)$.
A model of the embedding
$f\colon M\subset \CP(n)$ is given by the map
$\phi\colon R\to Q$ defined by $\phi(a^j)=\omega^j$
which is indeed a CDGA morphism
since $Q^{\geq 2m+2}=0$.
This induces an $R$-dgmodule structure on $Q$.
Let $I$ be a differential submodule of $Q$ complementary to the
image of $\phi$ which exists by Lemma \ref{lemma-complomega}.
Let $2k=2n-2m$, the codimension of $M$ inside $\CP(n)$.
\begin{lemma}\label{lem-shch}
The map $\phi^!\colon s^{-2k}Q\to R$ defined
by
$$
\phi^!(\omega^j)=l_M a^{j+k} \ {\rm and} \ \phi^!(I)=0
$$
is a shriek map, where $l_M$ is given by $\omega^m=l_Mu_M$ for the orientation class
$u_M$ determined by the almost complex structure. (See Definition \ref{def-mult}.)
\end{lemma}
\proof
This follows directly from the definition of a shriek map
(Definition \ref{def-shriek}).
\cqfd
From Section \ref{desc of model}
we get
\begin{equation}\label{modcp}
{\cal B}(\Lambda(a)/(a^{n+1}), Q)=\left(\frac{\Lambda(a)}{(a^{n+1})}\oplus
Q\otimes \Lambda^+(x,z),D\right)
\end{equation}
a CDGA with $|x|=2$, $|z|=2k-1$. The algebra structure
extends the algebra structure on $Q$ and the $\Lambda(a)/(a^{n+1})-$module
structure and the differential is determined by the Leibnitz law
and the following equations
\begin{eqnarray*}
D(a)&=&0\\
D(q \otimes x^j)&=&dq\otimes x^j\\
D(q\otimes z)&=& dq\otimes z + (-1)^{|q|}\left(\phi^!(s^{-2k}q)
+\sum_{i=0}^{k-1}
q\gamma_i \otimes x^{k-i}\right)
\end{eqnarray*}
Notice that taking the product of $a$ with an element of $Q$ gives multiplication
by $\omega$.
Note that for the next theorem our standard dimension hypothesis would
be $2n\geq 4m+3$ or $n\geq 2n+3/2$ but since $n$ is an integer this is
equivalent to $n\geq 2m+2$.
\begin{thm}
\label{thm-blowupCPn}
Let $(M,\omega)$ be a symplectic manifold of dimension $2m$ and suppose we have been given a
symplectic embedding
$f:(M,\omega)\rightarrow (\CP(n),\omega_0)$. Let
$Q$ be a CDGA model of M such that $Q$ is connected and $Q^{\geq 2m+1}=0$. Let
$\omega\in Q^2\cap \ker d$ be a representative of the symplectic form.
If $n\geq 2m+2$ then the CDGA ${\cal B}(\Lambda(a)/(a^{n+1}), Q)$ of Equation (\ref{modcp})
is a model of the blow-up $\widetilde \CP(n)$.
\end{thm}
\proof
Since $R=\Lambda(a)/(a^{n+1})$ and $Q^{\geq 2n+1}=0$, $\phi$
defined above is the only homotopy class
of a CDGA-morphism such that $\phi(a)=\omega$.

Therefore it is a model of the embedding $f$.
The morphism $\phi^!$ is clearly a shriek map.
Therefore all the hypotheses of Theorem \ref{lemma-modelblowup}
are fulfilled and the theorem follows.
\cqfd
\subsection{McDuff's example}
\label{McDuff}
Next we look at the model of the example of McDuff (\xcite{McDuff} or
\cite[Exercise 6.55]{McDuffSalamon}) that we now review.
We start with the Kodaira-Thurston manifold
(\xcite{Thurston} or \cite[Example II.2.1]{OpreaTralle} ) which is
a closed symplectic nilmanifold $V$ of dimension $4$ and is
defined as the product of the circle $S^1$ with the orbit space $\BR^3/\Gamma$
where $\Gamma$ is the uniform lattice of the integral upper triangular $3\times 3-$matrices
in the Heisenberg group.
A CDGA-model of that manifold
is given by the following exterior algebra on four generators of degree $1$
(see \cite[example II.1.7 (2)]{OpreaTralle}):
$$
(Q,d)=(\Lambda(u_1,y_1,v_1,t_1),du=dy=dt=0,dv=uy).
$$
The symplectic form is represented in this model by $\omega=uv+yt$.
By the symplectic embedding theorem of Tischler and Gromov (\xcite{Tischler},
\cite[3.4.2]{Gromov}) for $n\geq5$
there exists a symplectic embedding $f\colon V\hookrightarrow \CP(n)$
such that $f^*(a)=[\omega]$. To fulfill the hypotheses of Theorem \ref{thm-blowupCPn}
we will suppose that $n\geq 6$. Then $m=2$ and $k=n-2\geq4$.
In \cite[p 271]{McDuff} we see that the Chern classes of $V$ are trivial
and (\cite[Theorem 14.10]{MS}) the Chern classes of $\CP(n)$
are $c(\CP(n))=\sum_{i=0}^nc_i(\CP(n))=(1+a)^{n+1}$
where $a\in H^2(\CP(n),\BZ)$ represents the form $\omega_0$ described in Section \ref{prelim}.
Therefore the Chern classes of
the normal bundle $\nu$ of the embedding are given by the Equation
$c(\nu)\cdot c(V)=f^*(c(\CP(n)))$ which yields
$$
c(\nu)=1+(n+1)\omega+\frac{n(n+1)}{2}\omega^2.
$$
The morphism $\phi^!\colon s^{4-2n}Q\to \Lambda(a)/(a^{n+1})$ is
characterized by Lemma \ref{lem-shch} and satisifies,
$\phi^!(s^{4-2n}1)=2a^{n-2}$,
$\phi^!(s^{4-2n}uv)=\phi^!(s^{4-2n}yt)=a^{n-1}$,
$\phi^!(s^{4-2n}uvyt)=a^n$, and $\phi^!(s^{4-2n}\zeta)=0$
for any other monomial $\zeta$ in
$Q=\Lambda(u,y,v,t)$. Using these data in the definition
of the CDGA
$$
B={\cal B}(\Lambda(a)/(a^{n+1}), Q)=\left(\Lambda(a)/(a^{n+1})\oplus
\Lambda(u,y,v,t\otimes \Lambda^+(x,z),D\right)
$$
above
gives a model of the McDuff example.

From this CDGA-model of the McDuff example we recover
the Babenko-Taimanov \cite{BabenkoTaimanov1}
result.
\begin{thm}
There exist non-trivial Massey products in $\widetilde\CP(n)$.
\end{thm}
\proof
In $B$ we have
$D(v\otimes x^2)=(u\otimes x)\cdot (y \otimes x)$ and
$(y\otimes x)\cdot(y\otimes x)=0$, so
$D(0)=(y\otimes x)\cdot (y\otimes x)$.
Thus $vy\otimes x^3$ is in
the triple Massey product
$\langle [u\otimes x],[y\otimes x],[y\otimes x]\rangle$.
Also $[vy\otimes x^3]$ is not in the ideal of $H^*(\widetilde\CP(n))$ generated
by $[u\otimes x]$ and $[y\otimes x]$ so
$0\not\in \langle [u\otimes x],[y\otimes x],[y\otimes x]\rangle$.
\cqfd
In this example the Massey product $\langle u,y,y\rangle$ in $V$ became
$\langle [u\otimes x],[y\otimes x],[y\otimes x]\rangle$ in $\widetilde\CP(n)$.
This is the general way in which Massey products propagate to the blow-up.
More generally we show in \cite{LS1} that any obstruction to
formality in any manifold propagates to an obstruction in its blow-up in $\CP(n)$.
\subsection{The cohomology algebra of a blow-up}
\label{coho alg}
In the next theorem $f^!$ is the classical cohomological shriek map
(see Section \ref{sec-shr}). Also the algebra structure on
$H^*(W)\oplus H^*(V)\otimes \Lambda^+(x)$ is determined by
the formula on basic tensors
$$
(w,v\otimes x^i)\cdot(w',v'\otimes x^j)
=(ww',wv'\otimes x^j+(-1)^{|w'||v|}w'v\otimes x^i+vv'\otimes x^{i+j})
$$
where $w,w'\in H^*(W)$ and $v,v'\in H^*(V)$.
\begin{thm}\label{algcoho}
Let $f:V\rightarrow W$ be an embedding of closed oriented
manifolds of codimension $2k$
with a fixed complex structure on the normal bundle $\nu$.
Assume that $\dim W\geq 2 \dim V +3$.
Let $c_i(\nu)\in H^{2i}(V)$ be the
Chern classes.
Let $I$ be the ideal in $H^*(W)\oplus H^*(V)\otimes \Lambda^+(x)$
generated by the set
$$
\{ f^!(v)+\Sigma_{i=0}^{k-1}v \cdot c_i(\nu)\otimes x^{k-i}\colon v\in H^*(V) \}.
$$
Then we have an isomorphism of algebras
$$
H^*(\widetilde W)\cong (H^*(W)\oplus H^*(V)\otimes \Lambda^+(x))/I.
$$
\end{thm}
\proof
Recall the model of $\Apl(\widetilde{W})$ obtained
in Theorem \ref{lemma-modelblowup}:
$$
B={\cal B}(R,Q)=\left(R\oplus Q\otimes \Lambda^+(x,z),D\right)
$$
We define an increasing filtration on that CDGA by
\begin{eqnarray*}
F^0B&=&R\oplus Q\otimes \Lambda^+(x)\\
F^pB&=&B\mbox{\,\,\,\,\,\,\,\,\,for\,\,}p\geq1
\end{eqnarray*}
This filtration is compatible with the CDGA structure.
The $E_1-$term of the associated spectral sequence is
\begin{eqnarray*}
E_1^{0,*}&=&H(R)\oplus H(Q)\otimes \Lambda^+(x) \\
E_1^{1,*}&=& H(Q)\otimes \Lambda(x)\cdot z\\
E_1^{p,*}&=&0\mbox{\,\,for\,\,}p\not=0,1
\end{eqnarray*}
and the $d_1-$differential is non trivial only on $E_1^{1,*}$, being there
$$
d_1([q]\otimes x^r\otimes z)
=(-1)^{|q|}\left( f^!([q])+\sum_{i=0}^{k-1}[q]\cdot c_i(\nu)\otimes x^{r+k-i}\right)
$$
hence
$$
d_1([q]\otimes x^r \otimes z)=\pm [q]\otimes x^{r+k}
+({\rm terms \ with \ lower \ powers \ of \ } x).
$$
Thus $d_1$ is injective and therefore
$$
E_2=H(E_1,d_1)=\frac{H(R)\oplus H(Q)\otimes (\Lambda^+(x))}{\im\,d_1}
$$
with
$$
{\im} \ d_1=\left(\{f^!(v)+\sum_{i=0}^{k-1}v\cdot c_i(\nu)\otimes x^{k-i} :v\in H^*(V)\}\right)
$$
The $E_2-$term is concentrated in column $p=0$, so the spectral sequence collapses
at $E_2$ where $E_2=E^{0,*}$ as algebras.
Therefore there is no possibility for extensions and
$H^*(\widetilde{W})\cong E_\infty=E_2$ as algebras.
\cqfd
Theorem \ref{algcoho} determines the algebra structure of the cohomology of the blow-up
under the ``stable" condition: $\dim(W)\geq 2 \dim(V)+3$. This result is complementary to a
theorem of Gitler \cite[Theorem 3.11]{Gitler} that determines the cohomology algebra
$H^*(\widetilde W)$ under the hypothesis that
$f^*\colon H^*(W) \rightarrow H^*(V)$ is surjective.
\subsection{A Simple Example}
\label{vlast}
As noted in Section \ref{prelim} the Fubini-Study metric $\omega_0\in H^2(\CP(n))$ is an integral
symplectic form. Thus for any
$l\in \BZ$ with $l>0$, $l\omega_0\in H^2(\CP(n))$ is also one.
So by \cite{Tischler} there exists a symplectic embedding
$$
f_l\colon \CP(1) \rightarrow \CP(5)
$$
such that $f^*_l(\omega_0)=l \omega'_0$, where $\omega_0$ denotes the Fubini-Sudy metric
in $H^2(\CP(5))$ and $\omega'_0$ denotes it in $H^2(\CP(1))$.
Identifying $a$ with $\omega_0\in H^2(\CP(5))$ and
$a'$ with $\omega_0\in H^2(\CP(1))$, we get isomorphisms
$H^*(\CP(5))\cong \Lambda(a)/(a^6)$ and
$H^*(\CP(1))\cong \Lambda(a')/((a')^2)$.
Also $f^*_l(a)=la'$ and
$$
f^!_l:s^{-8}\Lambda(a')/((a')^2) \rightarrow \Lambda(a)/(a^6)
$$
is
given by the formulas
$$
\phi^!(s^{-8}1)=la^4
$$
$$
\phi^!(s^{-8}a')=a^5.
$$
Next we calculate $c_1(\nu)$, the first Chern class of the normal bundle.
Using the following three formulas
$$
f^*(c(\CP(5))=1+6f^*(a)=1+6la'
$$
$$
c(\CP(1))=1+2a'
$$
and
$$
c(\nu)c(\CP(1))=f^*(c(\CP(5)))\in H^*(\CP(1)),
$$
a simple calculation gives us that
$$
c(\nu)=1+(6l-2)a'
$$
and hence
$$
c_1(\nu)=(6l-2)a'.
$$
Theorem \ref{algcoho} tells us that the
cohomology of the blow-up $\widetilde\CP_l(5)$ of
$\CP(5)$ along $f_l$ is
\begin{equation}
\label{firstmodel}
H^*(\widetilde\CP_l(5))\cong
(\Lambda(a)/(a^6)\oplus (\Lambda(a')/((a')^2)\otimes \Lambda^+(x)))/
\left( la^4+(6l-2)a'x^3+x^4   \right)
\end{equation}
Since $f_l(a)=la'$, $a\cdot(1\otimes x)=la'\otimes x$ so we can write
$ax=la'x$ or $\frac{1}{l} ax=a'x$. Thus $a^2x=l^2(a')^2x=0$, and also
$la^4+(6l-2)a'x^3+x^4=\frac{1}{l}(l^2a^4+(6l-2)ax^3+lx^4)$. So
it is straightforward to check that we also get an isomorphism
\begin{equation}
\label{secondmodel}
H^*(\widetilde\CP_l(5))\cong
\Lambda(a,x)/\left( a^6,a^2x,l^2a^4+(6l-2)ax^3+lx^4\right)
\end{equation}
with the isomorphism in (\ref{secondmodel}) realized by mapping
$a$ to $a$ and $x$ to $x$ on the right hand side of
(\ref{firstmodel}).
\begin{prop}
If  $\widetilde\CP_l(5)$ and
$\widetilde\CP_r(5)$ are rationally homotopy equivalent then $l=r$.
\end{prop}
\proof
We have just seen that
$$
H^*(\widetilde\CP_l(5))\cong
\Lambda(a,x)/\left( a^6,a^2x,l^2a^4+(6l-2)ax^3+lx^4\right)
$$
and
$$
H^*(\widetilde\CP_r(5))\cong
\Lambda(b,y)/\left( b^6,b^2y,r^2b^4+(6r-2)by^3+ry^4\right)
$$
Assume that
$\widetilde\CP_l(5)$ and
$\widetilde\CP_r(5)$ are rationally homotopy equivalent. So
there is an isomorphism
$$
g\colon H^*(\widetilde\CP_l(5))\rightarrow
H^*(\widetilde\CP_r(5)).
$$
Write
$$
g(a)=\alpha_1 b+\beta_1 y
$$
and
$$
g(x)=\alpha_2 b +\beta_2 y.
$$
We will first show that $\beta_1=\alpha_2=0$.
We know that $a^2x=0$ in $H^*(\widetilde\CP_l(5))$,
so $g(a^2x)=0$. This implies that
$\alpha_1^2\alpha_2=0$, $\beta_1^2\beta_2=0$ and
$2\alpha_1\beta_1\beta_2+\alpha_2\beta_1^2=0$. Since
$g$ is an isomorphism it is straightforward to check that the only possibility is
$\beta_1=0=\alpha_2$.
\par
We know that $l^2a^4+(6l-2)ax^3+lx^4=0\in H^*(\widetilde \CP_l(5))$.
So looking at the coefficients of $b^4$, $by^3$ and $y^4$ in
$g(l^2a^4+(6l-2)ax^3+lx^4)=0$ we get that for some fixed
$\delta\in \BQ$
\begin{equation}
\label{name1}
l^2\alpha_1^4=\delta r^2 \ {\rm so } \
\delta=\frac{l^2}{r^2}\alpha_1^4
\end{equation}
\begin{equation}
\label{name2}
(6l-2)\alpha_1\beta_2^3=\delta(6r-2)
\end{equation}
and
\begin{equation}
\label{name3}
l\beta_2^4=\delta r.
\end{equation}
Observe that none of $l$, $r$, $\delta$, $\beta_2$ and $\alpha_1$ can be $0$ so we
can divide by them at will.
Equations (\ref{name1}) and (\ref{name3}) imply that
\begin{equation}
\label{name4}
\frac{l}{r}=\left( \frac{\beta_2}{\alpha_1} \right)^4
\end{equation}
Subbing (\ref{name1}) into (\ref{name2}), taking fourth powers and gathering
$\alpha_1$ to one side we get
$$
(6l-2)^4\left(\frac{\beta_2}{\alpha_1}\right)^{12}=\frac{l^8}{r^8}(6r-2)^4
$$
so using (\ref{name4}) we get
\begin{equation}
\label{name5}
\frac{(6l-2)^4}{l^5}=\frac{(6r-2)^4}{r^5}
\end{equation}
Let $h(x)=\frac{(6x-2)^4}{x^5}$. Calculus tells us that
$h(x)$ is decreasing $x\geq 2$, also $h(1)\not=h(2)$,
$h(1)\not=h(3)$ and $h(1)>h(4)$. Together these
imply that $h(x)$ takes distinct values on each
positive integer. Thus we
must have $r=l$ as required.
\cqfd
Since all of the
$\widetilde\CP_l(5)$ are not rationally homotopy equivalent they are
not integrally homotopy equivalent and hence not diffeomorphic.
Using the $H^*(\CP(5);\BZ)$ module structure on $H^*(\widetilde\CP_l(5);\BZ)$
gives another way to see that all of the $\widetilde\CP_l(5)$ are
not integrally homotopy equivalent.
If we look at blow-ups of $\CP(1)$ in $\CP(4)$ we would
have difficulty proving the last theorem because two of the
relations in the cohomology
algebra would be in the same degree. This leads to some interesting algebraic equations
which seem difficult to solve. However it can still be shown using the module
structure that the blow-ups have different integral homotopy type.
This leads us to the following question which seems more likely to have
a positive answer if the codimension is large.
\vskip 0.2cm
\noindent
{\bf Question:}
Let $(M,\omega)$ be an integral symplectic manifold, $l$ a positive integer
and $f_l:M \rightarrow \CP(n)$ be an embedding such that
$f_l^*(\omega_0)=l\omega$. Let $\widetilde\CP_l(n)$ be the blowup along $f_l$.
Are all of the $\widetilde\CP_l(n)$ rationally inequivalent?
\bigbreak

%

\bigskip

\noindent\underline{Addresses :}\\
Pascal Lambrechts\\
Universit\'e d'Artois\\
Universit\'e de Louvain\\
Institut Math\'ematique\\
Chemin du Cyclotron, 2\\
B-1348 Louvain-la-Neuve, BELGIUM\\
e-mail : {\tt lambrechts@math.ucl.ac.be}
\medbreak

\noindent
Donald Stanley\\
University of Regina
\\e-mail : {\tt stanley@math.uregina.ca}
\end{document}